\documentclass[9pt]{amsart}
\usepackage[T1]{fontenc}
\usepackage{amsmath,amssymb}
\usepackage{graphicx,xcolor}
\newtheorem{proposition}{Proposition}
\newtheorem{theorem}{Theorem}
\newtheorem{lemma}{Lemma}
\newtheorem{claim}{Claim}
\newtheorem{corollary}{Corollary}
\theoremstyle{definition}
\newtheorem{definition}{Definition}
\newtheorem{remark}{Remark}
\renewcommand{\matrix}[1]{\left(\begin{array}{cc} #1\end{array}\right)}
\newcommand{\minimatrix}[1]{\left(\begin{smallmatrix}#1\end{smallmatrix}\right)}
\newcommand{\smallfrac}[2]{\mbox{$\frac{#1}{#2}$}}
\DeclareFontFamily{U}{mathx}{\hyphenchar\font45}
\DeclareFontShape{U}{mathx}{m}{n}{
      <5> <6> <7> <8> <9> <10>
      <10.95> <12> <14.4> <17.28> <20.74> <24.88>
      mathx10
      }{}
\DeclareSymbolFont{mathx}{U}{mathx}{m}{n}
\DeclareFontSubstitution{U}{mathx}{m}{n}
\DeclareMathAccent{\widecheck}{0}{mathx}{"71}
\DeclareMathAccent{\widetilde}{0}{mathx}{"72}
\DeclareMathAccent{\widebar}{0}{mathx}{"73}
\DeclareMathAccent{\widevec}{0}{mathx}{"74}
\DeclareMathAccent{\widehat}{0}{mathx}{"70}
\newcommand{\wt}[1]{\widetilde{#1}}
\newcommand{\wh}[1]{\widehat{#1}}
\newcommand{\cal}[1]{{\mathcal #1}}
\newcommand{\cv}[1]{\underline{#1}}
\newcommand{\wc}[1]{\widecheck{#1}}
\def\A{{\mathbb A}}
\def\C{{\mathbb C}}
\def\CC{\overline{\mathbb C}}
\def\D{{\mathbb D}}
\def\H{{\mathbb H}}
\def\N{{\mathbb N}}
\def\R{{\mathbb R}}
\def\S{{\mathbb S}}
\def\Z{{\mathbb Z}}
\def\Re{{\rm Re}}
\def\Im{{\rm Im}}
\def\Ker{{\rm Ker}}
\def\Res{{\rm Res}}
\def\Sym{{\rm Sym}}
\def\Nor{{\rm Nor}}
\def\Pos{{\rm Pos}}
\def\Uni{{\rm Uni}}
\def\x{{\rm x}}
\def\cqfd{\hfill$\Box$}
\def\W{{\mathcal W}}
\def\Wp{{\mathcal W}^{\geq 0}}
\def\WR{{\mathcal W}_{\R}}
\def\WRp{{\mathcal W}_{\R}^{\geq 0}}
\def\su{\mathfrak{su}}
\def\sl{\mathfrak{sl}}
\def\xiS{\xi^S}
\def\PhiS{\Phi^S}
\def\GS{G^S}
\def\xiC{\xi^C}
\def\PhiC{\Phi^C}
\def\z{z}

\begin{document}
\title{Opening nodes in the DPW method: co-planar case}
\author{Martin Traizet}
\address{Institut Denis Poisson, CNRS UMR 7350 \\
Facult\'e des Sciences et Techniques \\
Universit\'e de Tours\\
France }
\email{martin.traizet@univ-tours.fr }
\thanks{Supported by ANR-19-CE40-0014 grant}
\maketitle

{\em Abstract: we combine the DPW method and opening nodes to construct
embedded surfaces of positive constant mean curvature with Delaunay ends in euclidean
space, with no limitation to the genus or number of ends.}
\medskip

\section{Introduction}
In \cite{dorfmeister-pedit-wu}, Dorfmeister, Pedit and Wu have shown that harmonic maps
from a Riemann surface to a symmetric space admit a Weierstrass-type
representation, which means that they can be represented in terms of
holomorphic data.
In particular, surfaces with constant mean curvature one (CMC-1 for short) in
euclidean space admit such a representation, owing to the fact that the Gauss map
of a CMC-1 surface is a harmonic map to the 2-sphere.
This representation is now called the DPW method and has been widely used
to construct CMC-1 surfaces in $\R^3$ and also constant mean curvature surfaces in homogeneous spaces such as the sphere ${\mathbb S}^3$ or hyperbolic space ${\mathbb  H}^3$:
see for example \cite{dorfmeister-haak,
dorfmeister-wu,
heller1,
heller2,
heller3,
kilian-kobayashi-rossman-schmitt,
kilian-mcintosh-schmitt,
schmitt-kilian-kobayashi-rossman}.
\medskip

The input data for the DPW method is called the DPW potential. In principle, all CMC surfaces can be obtained by the DPW method. But in practice, one has to solve a Monodromy Problem,
akin to the Period Problem for the construction of minimal surfaces via the Weierstrass Representation.
So in general the topology of the constructed examples is limited
or symmetries are imposed in order to reduce the number of equations to be solved.
In contrast, Kapouleas \cite{kapouleas} has constructed embedded CMC surfaces with
no limitation on the genus or number of ends by gluing round spheres and
pieces of Delaunay surfaces, using PDE methods.
It seems an interesting question to see whether such gluing constructions can be achieved by the DPW method.
\medskip

In \cite{nnoids}, we proposed a DPW potential for CMC $n$-noids: genus zero CMC-1 surfaces with $n$ Delaunay-type ends. They look like a round sphere with $n$ half-Delaunay surfaces with small necksize attached at prescribed points. They are a particular case of the construction of Kapouleas in \cite{kapouleas}.
The potential is natural, in the sense that it is a perturbation of the standard spherical potential.
This potential has been adapted to minimal surfaces in $\H^3$ and $AdS^3$ in \cite{bobenko-heller-schmitt} and CMC>1 surfaces in $\H^3$ in \cite{raujouan2}.
\medskip

In \cite{minoids}, we proposed a DPW potential for another type of CMC $n$-noids which look like a minimal
$n$-noid (a genus zero minimal surface with $n$ catenoidal ends) whose catenoidal ends have been replaced by Delaunay ends. They had already been constructed by Mazzeo and Pacard in \cite{mazzeo-pacard} using PDE methods.
The potential is derived in a natural way from the Weierstrass data of the minimal $n$-noid.
It has also been adapted to CMC>1 surfaces in $\H^3$ in \cite{raujouan2}.
\medskip

In an unpublished paper \cite{opening-nodes}, we proposed a DPW potential for all the surfaces constructed by Kapouleas in \cite{kapouleas}. The potential was, however, quite complicated and hardly natural, and the paper was long and technical.
Our goal in this paper is to propose a much simpler and natural DPW potential in 
a particular, but still interesting case: when all the centers of the spheres to be glued together are in the same plane. The resulting CMC surface is invariant by symmetry with respect to that plane.
The symmetry allows us to take advantage of the fact that the standard holomorphic frame for Delaunay surfaces is unitary on the unit circle, which is a big asset for the resolution of the Monodromy Problem.
\medskip

The underlying Riemann surface is defined by opening nodes, which is a standard model for
Riemann surfaces with ``small necks''.
The theory of opening nodes has been used by the author to construct minimal surfaces
in euclidean space via the classical Weierstrass Representation
(see for example \cite{nosym} or \cite{triply})
or CMC-1 surfaces in hyperbolic space via Bryant Representation \cite{bryant}.
\medskip

One difficulty with the DPW method is that unlike the Weierstrass data of minimal surfaces, the
DPW potential has little geometric content so it is hard to guess a candidate for the construction of CMC surfaces with given geometric features.
The heuristic that we follow is that the DPW potential should be a perturbation of the spherical potential
where the surface is close to a round sphere and of the ``catenoidal potential'' where the surface has small catenoidal necks.
\medskip

\medskip
This paper opens up the possibility of opening nodes in the DPW method.
We hope the ideas developed in this paper will be useful to the contruction of minimal and CMC surfaces in other space forms.
\section{Main result}
Our goal is to contruct CMC surfaces by gluing spheres and half-Delaunay surfaces. The layout of these pieces is encoded by a weighted graph in the horizontal plane.
\begin{definition}
\label{definition:graph}
A horizontal weighted graph $\Gamma$ is the following data:
\begin{itemize}
\item A finite number of points $v_j\in\R^2$ for $j\in J$, called vertices. Here $J\subset\N^*$ is a finite set used to index vertices.
\item A symmetric subset $E\subset (J\times J)\setminus\Delta$ where $\Delta$ is the diagonal of 
$J\times J$, whose elements are called edges. Two vertices $v_j$ and $v_k$ are adjacent if $(j,k)\in E$.
\item A finite set of half-lines $\Delta_{jk}\subset\R^2$ for $(j,k)\in R$, called rays, such that $\Delta_{jk}$ has endpoint $v_j$. Here $R\subset J\times (\N^*\setminus J)$ is a finite set used to index rays.
\item Each edge or ray is given a non-zero weight $\tau_{jk}$, $(j,k)\in E\cup R$, with $\tau_{jk}=\tau_{kj}$
for $(j,k)\in E$.
\end{itemize}
\end{definition}
For $j\in J$, we denote $E_j=\{k\in J:(j,k)\in E\}$ the set of edges issued from the vertex $v_j$,
and $R_j=\{k\in \N^*,(j,k)\in R\}$ the indices of the rays issued from the vertex $v_j$.
Also we denote $E^+=\{(j,k)\in E:j<k\}$.

\begin{figure}
\begin{center}
\includegraphics[height=3cm]{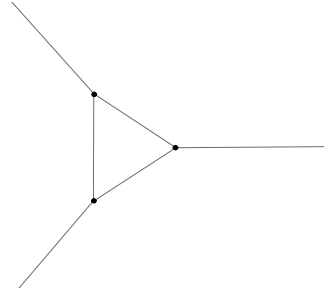}
\hspace{0.5cm}
\includegraphics[height=3cm]{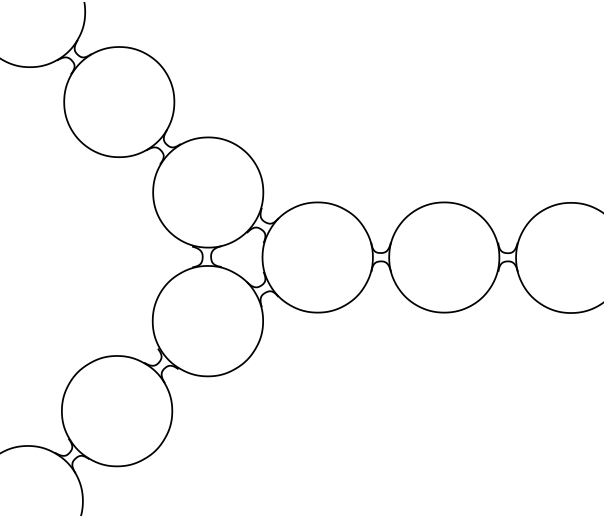}
\end{center}
\caption{A balanced graph with three edges and three rays (left) and a CMC-1 surface of genus one with three Delaunay-type ends
in the corresponding family (right).}
\label{fig1}
\end{figure}

Given a horizontal weighted graph $\Gamma$ with length-2 edges, we can construct a singular CMC-1 surface $M_0$ as follows. We identify $\R^2$ with the horizontal plane $x_3=0$.
\begin{itemize}
\item For $j\in J$, place a radius-1 sphere centered at the vertex $v_j$, so if $v_j$ and $v_k$ are adjacent, the corresponding spheres are tangent.
\item For each $(j,k)\in R$, place an infinite chain of radius-1 spheres with centers on $\Delta_{jk}$ at even distance from $v_j$.
\end{itemize}
Our goal in this paper is to construct a family of CMC-1 surfaces $(M_t)_{0<t<\epsilon}$ by desingularizing $M_0$, replacing
all tangency points between adjacent spheres by catenoidal necks of size $\simeq t\tau_{jk}$
(see Figure \ref{fig1}).
This is only a heuristic way to describe the result, and is not the way we will construct $M_t$
(although this is how Kapouleas does in \cite{kapouleas}).

For the construction to succeed, the weighted graph $\Gamma$ must satisfy a balancing condition.
For $(j,k)\in E$, we denote $\ell_{jk}=|v_j-v_k|$ and
$u_{jk}$ the unitary vector $(v_k-v_j)/\ell_{jk}$, so $u_{kj}=-u_{jk}$.
For $(j,k)\in R$, we denote $u_{jk}$ the unitary vector in the direction of the ray $\Delta_{jk}$.
\begin{definition} For $j\in J$,
we define the force $F_i$ on the vertex $v_j$ by
\begin{equation}
\label{eq:force}
F_j=\sum_{k\in E_j\cup R_j}\tau_{jk}u_{jk}.
\end{equation}
A horizontal weighted graph $\Gamma$ is balanced if $F_j=0$ for all $j\in J$.
\end{definition}
To solve our problem, we need to perturb $\Gamma$ in order
to prescribe small variations of edge-lengths and forces.
The parameters available to deform $\Gamma$ are the vertices $v_j\in\R^2$ for $j\in J$,
the unitary vectors $u_{jk}$ for $(j,k)\in R$ and the weights $\tau_{jk}\in\R$ for $(j,k)\in E^+\cup R$.
\begin{definition}
A horizontal weighted graph $\Gamma$ is non-degenerate if the jacobian of
the map
$\big( (F_j)_{j\in J},(\ell_{jk})_{(j,k)\in E^+}\big)$ with respect to the above parameters is onto.
\end{definition}
\begin{theorem}
\label{theorem:main}
Let $\Gamma$ be a balanced, non-degenerate horizontal weighted graph with length-2 edges.
There exists a smooth 1-parameter family of immersed CMC-1 surfaces $(M_t)_{0<t<\epsilon}$
with the following properties:
\begin{enumerate}
\item $M_t$ converges to $M_0$ as $t\to 0$. The convergence is for the
Hausdorf distance on compact sets of $\R^3$.
\item $M_t$ is homeomorphic to a tubular neighborhood of $\Gamma$.
\item $M_t$ is symmetric with respect to the horizontal plane.
\item For each $(j,k)\in R$, $M_t$ has a Delaunay end with weight
$\simeq 2\pi t\tau_{jk}$ and whose axis converges as $t\to 0$ to the ray $\Delta_{jk}$.
\item If all weights are positive, then $M_t$ is Alexandrov-embedded.
\item If moreover $\Gamma$ is pre-embedded, then $M_t$ is embedded.
\end{enumerate}
\end{theorem}
\begin{definition}
\label{def-preembedded}
Following Kapouleas (Definition 2.2 in \cite{kapouleas}),
we say that $\Gamma$ is pre-embedded if the distance between any two edges or rays which have no common endpoint is greater than $2$ and
the angle between any two edges or rays with a common endpoint is greater than $60^\circ$.
\end{definition}
%As already said, Theorem \ref{thm-main} was proved using a completely different method in
%\cite{kapouleas}.
%In the simplest case $N=1$, there is no need for Opening Nodes and Theorem \ref{thm-main} is proved using the DPW method in \cite{nnoids}. We follow the same strategy
%to define the DPW potential and we will use some of the results in \cite{nnoids}.
%\subsection{Reduction to length-2 edges}
%Let $\Gamma$ be a graph with $N$ vertices and even-length edges. Assume that $\Gamma$ has an edge $e_{ij}$ of length $\ell_{ij}\geq 4$. We can define a new graph $\widetilde{\Gamma}$ as follows: insert a new vertex $\v_{N+1}$
%on the edge $e_{ij}$ at distance $2$ from $\v_i$. Replace the edge $e_{ij}$ by
%the edges $[\v_i,\v_{N+1}]$ and $[\v_{N+1},\v_j]$, with respective lenghts
%$2$ and $\ell_{ij}-2$. Assign to each new edge the weight $\tau_{ij}$.
%The new graph $\widetilde{\Gamma}$ is clearly balanced.
%\begin{proposition}
%\label{prop-reduction}
%If $\Gamma$ is non-degenerate then $\widetilde{\Gamma}$ is
%non-degenerate. If $\Gamma$ is pre-embedded then $\widetilde{\Gamma}$
%is pre-embedded.
%\end{proposition}
%The proof of Proposition \ref{prop-reduction} is elementary and is omitted.
%Thanks to Proposition \ref{prop-reduction}, we can transform by induction the graph $\Gamma$ into
%a balanced, non-degenerate graph with length-2 edges. 
%Therefore, it suffices to prove Theorem \ref{thm-main} in the case where
%all edges have length 2.
%
\section{Background}
\subsection{Functional spaces}
\label{section:functional-spaces}
The DPW method uses loop groups, which are groups of smooth functions from the unit circle
$\S^1\subset\C$ to a matrix group. The circle variable is denoted $\lambda$.
The DPW method is usually formulated in the category of smooth maps, but since we plan to use the
Implicit Function Theorem, we need a Banach space.
We adopt the following choice, following \cite{nnoids,minoids}.

Fix some $\rho>1$ and let $\D_{\rho}\subset\C$ be the disk $|\lambda|<\rho$ and
$\A_{\rho}\subset\C$ the annulus $\rho^{-1}<|\lambda|<\rho$.
We decompose a smooth function $f:\S^1\to\C$ in Fourier
series
$$f(\lambda)=\sum_{i\in\Z}f_i\lambda^i$$
and define
$$\|f\|=\sum_{i\in\Z}|f_i|\rho^{|i|}$$
Let $\W$ be the space of functions $f$ with finite norm.
This is a Banach algebra, owing to the fact that the weight $\rho^{|i|}$ is submultiplicative
(see Section 4 in \cite{grochenig}).
Functions in $\W$ extend holomorphically to $\A_\rho$.

We define $\Wp$, $\W^{>0}$, $\W^{\leq 0}$ and $\W^{<0}$ as the subspaces of functions $f$ such that $f_i=0$
for $i<0$, $i\leq 0$, $i>0$ and $i\geq 0$, respectively.
Functions in $\Wp$ extend holomorphically to the disk $\D_{\rho}$ and
satisfy $|f(\lambda)|\leq \|f\|$ for all $\lambda\in\D_{\rho}$.
We write $\W^0\sim\C$ for the subspace of constant functions, so we have
a direct sum $\W=\W^{<0}\oplus\W^0\oplus\W^{>0}$.
%(The Banach algebra $\W$ is said to be decomposable, see \cite{gohberg} page 70.)
A function $f$ will be decomposed as $f=f^-+f^0+f^+$ with
$(f^-,f^0,f^+)\in\W^{<0}\times\W^0\times\W^{>0}$
(and of course $f^0=f_0$).

We define the conjugation operator by
$$\overline{f}(\lambda)=\overline{f(\overline{\lambda})}=\sum_{i\in\Z}\overline{f_i}\lambda^i.$$
We denote $\Re(f)=\frac{1}{2}(f+\overline{f})$ and $\Im(f)=\frac{1}{2i}(f-\overline{f})$
and define $\WR$ as the subspace of functions in $\W$ such that $\Im(f)=0$,
and $\WRp=\WR\cap\Wp$.

We also define the star operator by
$$f^*(\lambda)=\overline{f(1/\overline{\lambda})}=\sum_{i\in\Z}\overline{f_{-i}}\lambda^i.$$
The involution $f\mapsto f^*$ exchanges $\Wp$ and $\W^{\leq 0}$.
We have $\lambda^*=\lambda^{-1}$ and $c^*=\overline{c}$ if $c$ is a constant.
A function $f$ is real on the unit circle if and only if $f=f^*$.
Note that conjugation and star commute.

There is a theory of holomorphic functions between complex Banach space, which retain most
properties of holomorphic functions of several variables. A good reference is \cite{chae}.
\subsection{Loop groups}
\begin{itemize}
\item If $G$ is a matrix Lie group, we denote $\Lambda G$ the Banach Lie group
of maps $\Phi:\S^1\to G$ whose entries are in $\W$.
\item If $\mathfrak{g}$ is the Lie algebra of $G$, the Lie algebra of $\Lambda G$ is the set of maps
$\varphi:\S^1\to\mathfrak{g}$ whose entries are in $\W$ and is denoted $\Lambda\mathfrak{g}$.
\item $\Lambda_+ SL(2,\C)\subset\Lambda SL(2,\C)$ is the subgroup of maps $B$ whose entries are in $\Wp$, with $B\mid_{\lambda=0}$ upper triangular.
\item $\Lambda_+^{\R}SL(2,\C)\subset\Lambda_+ SL(2,\C)$ is the subgroup of maps $B$ such that
$B\mid_{\lambda=0}$ has positive entries on the diagonal.
\end{itemize}
The following result is the corner stone of the DPW method. It is usually formulated for
smooth loops \cite{presley-segal}, but adapts with no difficulty to loops with entries in $\W$
(see details in Section 3.6 of \cite{minoids}).
\begin{theorem}[Iwasawa decomposition]
\label{theorem:iwasawa}
The multiplication $\Lambda SU(2)\times \Lambda_+^{\R} SL(2,\C)\to\Lambda SL(2,\C)$
is a smooth diffeomorphism (in the sense of smooth maps between Banach manifolds). The unique splitting of an element $\Phi\in\Lambda SL(2,\C)$ as
$\Phi=FB$ with $F\in\Lambda SU(2)$ and $B\in\Lambda_+^{\R} SL(2,\C)$
is called Iwasawa decomposition.
$F$ is called the unitary factor of $\Phi$ and
denoted $\Uni(\Phi)$.
$B$ is called the positive factor and denoted $\Pos(\Phi)$.
\end{theorem}
\subsection{The DPW method}
In the DPW method, one identifies $\R^3$ 
with the Lie algebra $\su(2)$ by
$$(x_1,x_2,x_3)\in\R^3\quad\sim\quad -i\matrix{
-x_3&x_1+i x_2\\x_1-i x_2 &x_3}\in\su(2).$$
%We have $\det(X)=\|x\|^2$.
%The group $SU(2)$ acts as linear isometries on $\su(2)$ by
%conjugation: $H\cdot X=HXH^{-1}$.
The input data for the DPW method is a quadruple $(\Sigma,\xi,z_0,\phi_0)$ where
$\Sigma$ is a Riemann surface,
$\xi$ is a $\Lambda\sl(2,\C)$-valued holomorphic 1-form on $\Sigma$
of the following special form
\begin{equation}
\label{eq:xi}
\xi=\matrix{\alpha&\lambda^{-1}\beta\\ \gamma&-\alpha}
\end{equation}
where $\alpha$, $\beta$, $\gamma$
are $\Wp$-valued holomorphic 1-forms on $\Sigma$,
$z_0\in\Sigma$ is a base point
 and $\phi_0\in\Lambda SL(2,\C)$ is an initial condition.
$\xi$ is called the DPW potential.
If $\Sigma$ is simply connected, the DPW method is the following procedure:
\begin{itemize}
\item Solve the Cauchy Problem on $\Sigma$:
\begin{equation}
\label{eq:cauchy-problem}
\left\{\begin{array}{l}
d_z\Phi=\Phi\xi\\
\Phi(z_0)=\phi_0\end{array}\right.\end{equation}
to obtain a solution
$\Phi: \Sigma\to\Lambda SL(2,\C)$.
\item Compute the Iwasawa decomposition $(F(z),B(z))$ of $\Phi(z)$ for $z\in\Sigma$.
\item Define $f:\Sigma\to\su(2)\sim\R^3$ by the Sym-Bobenko formula:
\begin{equation}
\label{eq:sym-bobenko}
f(z)=\Sym(F(z))=-2 i\,\frac{\partial F(z)}{\partial\lambda}F(z)^{-1}\mid_{\lambda=1}.
\end{equation}
Then $f$ is a CMC-1 (branched) conformal immersion.
$f$ is regular at $z$ (meaning unbranched) if and only if $\beta^0(z)\neq 0$.
Its Gauss map is given by
\begin{equation}
\label{eq:normal}
N(z)=\Nor(F(z))=-i\,F(z)\matrix{-1 &0 \\ 0 & 1}F(z)^{-1}\mid_{\lambda=1}.
\end{equation}
The DPW method actually constructs a moving frame for $f$ and
the differential of $f$ is given by
\begin{equation}
\label{eq:df}
df(z)=2 i\,B_{11}^0(z)^2F(z)\matrix{0&\beta^0(z)\\\overline{\beta^0(z)}&0}F(z)^{-1}\mid_{\lambda=1}.
\end{equation}
\end{itemize}
\subsection{The Monodromy Problem}
If $\Sigma$ is not simply connected, lift the DPW potential $\xi$ to the universal cover $\wt{\Sigma}$
of $\Sigma$ and choose a point $\wt{z_0}$ in the fiber of $z_0$.
Solve the Cauchy Problem $d\Phi=\Phi\xi$ in $\wt{\Sigma}$ with initial condition $\Phi(\wt{z}_0)=\phi_0$
to define $\Phi:\wt{\Sigma}\to\Lambda SL(2,\C)$. The DPW method produces an immersion $f:\wt{\Sigma}\to\R^3$.

%We identify the group of Deck transformations of $\wt{\Sigma}\to\Sigma$ with the fundamental group
%$\pi_1(\Sigma,z_0)$. (This identification is not canoncial and depends on the choice of $\wt{z}_0$.)
For $\gamma\in\pi_1(\Sigma,z_0)$, let $\wt{\gamma}$ be the lift of $\gamma$ to $\wt{\Sigma}$ such that $\wt{\gamma}(0)=\wt{z}_0$. The monodromy of $\Phi$ with respect to $\gamma$ is
$$\cal{M}(\Phi,\gamma)=\Phi(\wt{\gamma}(1))\Phi(\wt{\gamma}(0))^{-1}$$ 
The standard condition which ensures that the immersion $f$ descends to
a well defined immersion on $\Sigma$
is the following system of equations, called the Monodromy Problem:
\begin{equation}
\label{eq:monodromy-problem}
\forall \gamma\in\pi_1(\Sigma,z_0)\quad
\left\{\begin{array}{lc}
\cal{M}(\Phi,\gamma)\in\Lambda SU(2)\quad &(i)\\
\cal{M}(\Phi,\gamma)\mid_{\lambda=1}=\pm I_2\quad&(ii)\\
\frac{\partial}{\partial\lambda}\cal{M}(\Phi,\gamma)\mid_{\lambda=1}=0\quad &(iii)
\end{array}\right.\end{equation}
We will formulate the Monodromy Problem using the notion of principal solution
(see Chapter 3.4 in \cite{teschl}).
\begin{definition}
Let $\gamma:[0,1]\to\Sigma$ be a path, not necessarily closed.
Let $Y:[0,1]\to\Lambda SL(2,\C)$ be the solution of the Cauchy Problem
$$\left\{\begin{array}{l}
Y'(s)=Y(s)\,\xi(\gamma(s))(\gamma'(s))\\
Y(0)=I_2\end{array}\right.$$
The principal solution of $\xi$ with respect to $\gamma$ is $\cal{P}(\xi,\gamma)=Y(1)$.
\end{definition}
In other words, $\cal{P}(\xi,\gamma)$
is the value at $\gamma(1)$ of the analytical continuation along $\gamma$ of the solution of the Cauchy Problem \eqref{eq:cauchy-problem} with initial condition $\Phi(\gamma(0))=I_2$.
If $p$, $q$ are two points on $\Sigma$ and the path $\gamma$ from $p$ to $q$ is clear from the context, we will sometime write $\cal{P}(\xi,p,q)$ for $\cal{P}(\xi,\gamma)$.
The principal solution has the following properties, which follow easily from its definition:
\begin{itemize}
\item $\cal{P}(\xi,\gamma)$ only depends on the homotopy class of
$\gamma$.
\item The principal solution is a morphism for the product of paths:
If $\gamma_1$ and $\gamma_2$ are two paths such that $\gamma_1(1)=\gamma_2(0)$
then
$$\cal{P}(\xi,\gamma_1\gamma_2)=\cal{P}(\xi,\gamma_1)\cal{P}(\xi,\gamma_2).$$
\item If $\psi:\Sigma_1\to\Sigma_2$ is a holomorphic map, $\xi$ is a potential on $\Sigma_2$ and $\gamma$ is a path on $\Sigma_1$, then
$$\cal{P}(\psi^*\xi,\gamma)=\cal{P}(\xi,\psi(\gamma)).$$
\item If $\sigma:\Sigma_1\to\Sigma_2$ is a anti-holomorphic map, then
$$\cal{P}(\overline{\sigma^*\xi},\gamma)=\overline{\cal{P}(\xi,\sigma(\gamma))}.$$
\end{itemize}
Back to the DPW method, if the initial condition is $\Phi(z_0)=I_2$, which will be the case in this paper,
the Monodromy Problem is equivalent to the following problem:
\begin{equation}
\label{eq:monodromy-problem-principal}
\forall \gamma\in\pi_1(\Sigma,z_0)\quad
\left\{\begin{array}{lc}
\cal{P}(\xi,\gamma)\in\Lambda SU(2)\quad &(i)\\
\cal{P}(\xi,\gamma)\mid_{\lambda=1}=\pm I_2\quad&(ii)\\
\frac{\partial}{\partial\lambda}\cal{P}(\xi,\gamma)\mid_{\lambda=1}=0\quad &(iii)
\end{array}\right.\end{equation}
\subsection{Gauging and the Regularity Problem}
\begin{definition}
A gauge on $\Sigma$ is a holomorphic map $G:\Sigma\to\Lambda_+ SL(2,\C)$.
\end{definition}
Let $\Phi$ be a solution of $d\Phi=\Phi\xi$ and $G$ be a gauge.
Let $\wh{\Phi}=\Phi G$. Then $\Phi$ and $\wh{\Phi}$
define the same immersion $f$ via the DPW method.
The gauged potential is
$$\wh{\xi}:=\wh{\Phi}^{-1}d\wh{\Phi}=G^{-1}\xi G+G^{-1}d G$$
and is denoted $\xi\cdot G$, the dot denoting the action of the gauge
group on the potential.
Gauging does not change the monodromy of $\Phi$.
\begin{definition}
We say that $\xi$ is regular at $p\in\Sigma$ if $\beta^0(p)\neq 0$. This ensures
that the immersion $f$ is unbranched at $p$.
\end{definition}
In general $\Sigma$ is a compact Riemann surface $\overline{\Sigma}$ minus a finite number of
points, and the potential $\xi$ extends meromorphically to $\overline{\Sigma}$.
\begin{definition}
We say that a pole $p$ of $\xi$ is an apparent singularity if there exists a meromorphic gauge $G$, defined in a neighborhood of $p$, such that $\xi\cdot G$ extends holomorphically at $p$ and is regular. This ensures that the immersion $f$ extends analytically at $p$.
\end{definition}
Our potential will have two kinds of poles: some of them will be ends of the immersion $f$, the others will be apparent singularities.
Note that $\xi$ must have apparent singularities at the zeros of $\beta^0$ for $f$ to be regular.
 If $\overline{\Sigma}$ has positive genus, $\beta^0$ must have zeros on $\overline{\Sigma}$ so apparent singularities cannot be avoided.
\subsection{Dressing and rigid motions}
Let $\Phi$ be a solution of the Cauchy Problem \eqref{eq:cauchy-problem}.
Let $H\in\Lambda SU(2)$ and define $\wt{\Phi}(z)=H\Phi(z)$.
Then $\wt{\Phi}$ solves $d\wt{\Phi}=\wt{\Phi}\xi$ and the Iwasawa decomposition of
$\wt{\Phi}$ is $\wt{F}=HF$ and $\wt{B}=B$.
The Sym-Bobenko formula gives
$$\wt{f}(z)=\Sym(\wt{F}(z))=\left(HfH^{-1}-2i\frac{\partial H}{\partial\lambda}H^{-1}\right)\mid_{\lambda=1}.$$
Consequently, we define a left action of $\Lambda SU(2)$ on $\su(2)$ by
\begin{equation}
\label{eq:action}
H\cdot x=\left( HxH^{-1}-2i\frac{\partial H}{\partial\lambda}H^{-1}\right)\mid_{\lambda=1}.
\end{equation}
The action is by rigid motion and $\wt{f}=H\cdot f$.
The Monodromy Problems for $\Phi$ and $H\Phi$ are equivalent because $H\in\Lambda SU(2)$.
\subsection{Spherical and catenoidal potentials}
\label{section:spherical}
Delaunay surfaces are obtained from the following standard potential on $\C^*$:
$$\xi=\matrix{0&\lambda^{-1}r+s\\\lambda r+s&0}\frac{dz}{z}$$
with initial condition $\Phi(1)=I_2$, where $r,s$ are non-zero real numbers such that $r+s=\frac{1}{2}$.
There are two limiting cases of interest to us:
\begin{itemize}
\item Spherical limit: $(r,s)=(1/2,0)$ gives
$$\xiS=\matrix{0&\lambda^{-1}/2\\\lambda/2&0}\frac{dz}{z}$$
which we call the spherical Delaunay potential.
The corresponding solution is
$$\PhiS(z)=\frac{1}{2\sqrt{z}}\matrix{z+1&\lambda^{-1}(z-1)\\\lambda(z-1)&z+1}.$$
It Iwasawa decomposition is
$$F^S(z)=\frac{1}{\sqrt{2}\sqrt{1+|z|^2}}\matrix{\overline{z}+1&\lambda^{-1}(z-1)\\
\lambda(1-\overline{z})&z+1}\matrix{e^{i\theta/2}&0\\0&e^{-i\theta/2}}$$
$$B^S(z)=\frac{1}{\sqrt{2|z|}\sqrt{1+|z|^2}}\matrix{2|z|&0\\\lambda(|z|^2-1)&1+|z|^2}$$
where $\theta=\arg(z)$.
The Sym-Bobenko formula \eqref{eq:sym-bobenko} and Equation \eqref{eq:normal} give
$$f^S(z)=\frac{-i}{1+|z|^2}\matrix{-|z-1|^2&1-|z|^2-z+\overline{z}\\1-|z|^2+z-\overline{z}&|z-1|^2}
\sim\frac{1}{1+|z|^2}\left(1-|z|^2,-2\,\Im(z),|z-1|^2\right)$$
$$N^S(z)=\frac{-i}{1+|z|^2}\matrix{-z-\overline{z}&|z|^2-1+z-\overline{z}\\ |z|^2-1+\overline{z}-z&
z+\overline{z}}
\sim\frac{1}{1+|z|^2}\left(|z|^2-1,2\,\Im(z),2\,\Re(z)\right).$$
Consider the rigid motion
\begin{equation}
\label{eq:Psi}
\Psi(x_1,x_2,x_3)=(1-x_3,-x_2,-x_1).
\end{equation}
Then
\begin{equation}
\label{eq:PsifS}
\Psi\circ f^S(z)=\frac{1}{1+|z|^2}\left(2\,\Re(z),2\,\Im(z),|z|^2-1\right)=\pi^{-1}(z)
\end{equation}
where $\pi:\C\cup\{\infty\}\to\S^2$ is the stereographic projection from the north pole.
The poles at $0$ and $\infty$ are of course apparent singularities. This is confirmed by the following gauge:
$$\GS(z)=\matrix{\frac{1+z}{\sqrt{z}}&0\\\lambda\frac{1-z}{\sqrt{z}}&\frac{\sqrt{z}}{1+z}}$$
A computation gives
$$\xiS\cdot \GS=\matrix{0&\lambda^{-1}\\0&0}\frac{dz}{2(z+1)^2}$$
which is regular at $0$ and $\infty$.
\item Catenoidal limit: $(r,s)=(0,1/2)$ gives
$$\xiC=\matrix{0&1/2\\1/2&0}\frac{dz}{z}$$
which we call the catenoidal Delaunay potential.
The corresponding solution is
$$\PhiC(z)=\frac{1}{2\sqrt{z}}\matrix{z+1&z-1\\ z-1&z+1}$$
which does not depend on $\lambda$, so the immersion degenerates into the point $0$.
A computation gives
$$N^C(z)=\frac{1}{1+|z|^2}\left(1-|z|^2,2\,\Im(z),2\,\Re(z)\right).$$
which is a conformal diffeomorphism from $\C\cup\{\infty\}$ to $\S^2$.
\end{itemize}
\subsection{Duality}
Let
$$K(\lambda)=\matrix{0&i\lambda^{-1/2}\\i \lambda^{1/2}&0}.$$
\begin{definition}
The dual potential of $\xi$ is
$$\xi^{\dagger}=K\xi K^{-1}=\matrix{-\alpha&\lambda^{-1}\gamma\\\beta&\alpha}.$$
\end{definition}
The Delaunay spherical and catenoidal potentials are dual to each other.
Note that $K$ is not a gauge.
Duality transforms the immersion in the following explicit way.
Let $\Phi^{\dagger}=K\Phi K^{-1}$ be the solution of
$d\Phi^{\dagger}=\Phi^{\dagger}\xi^{\dagger}$ with initial condition
$\Phi^{\dagger}(z_0)=K\Phi(z_0)K^{-1}$.
The Iwasawa decomposition of $\Phi^{\dagger}$ is
$F^{\dagger}=KFK^{-1}$ and $B^{\dagger}=KB K^{-1}.$
The Sym-Bobenko formula gives:
$$f^{\dagger}(z)=\matrix{0&i\\i&0}\left[f(z)+N(z)-i\matrix{1&0\\0&-1}\right]\matrix{0&-i\\-i&0}.$$
In other words, up to a rigid motion, the dual (branched) immersion $f^{\dagger}$ is the parallel surface at distance one to $f$.
\section{Strategy}
\label{section:strategy}
Fix a horizontal weighted graph $\Gamma$. Until Section \ref{section:implicit}, we do not assume that
$\Gamma$ is balanced nor has length-2 edges.
Without loss of generality, we may assume (by rotating the graph $\Gamma$) that $u_{jk}\neq\pm 1$
for all $(j,k)\in E\cup R$.
We denote $\CC$ the Riemann sphere $\C\cup\{\infty\}$.
Take a copy of the Riemann sphere $\CC_j$ for each $j\in J$, and a copy of the Riemann sphere $\CC_{jk}$ for each $(j,k)\in E^+$.
For each $(j,k)\in E^+$, identify the point $z=u_{jk}$ in $\CC_j$ with the point $z=1$ in $\CC_{jk}$,
and the point $z=u_{kj}$ in $\CC_k$ with the point $z=-1$ in $\CC_{jk}$.
This defines a compact Riemann surface with nodes $\overline{\Sigma}_0$
(the nodes are the double points created when identifying pairs of points).

Consider the meromorphic DPW potential $\xi_0$ on $\overline{\Sigma}_0$ defined by $\xi_0=\xiS$ in $\CC_j$ for $j\in J$ and
$\xi_0=\xiC$ in $\CC_{jk}$ for $(j,k)\in E^+$.
Fix an arbitrary $j_0\in J$ and take as base point $z_0$ the point $z=1$ in $\CC_{j_0}$ and the initial
condition $\Phi(z_0)=I_2$.
The fundamental group $\pi_1(\overline{\Sigma}_0,z_0)$ is generated by paths made of unit circular arcs connecting the nodes. Whenever a path $\gamma$ crosses a node, we require the fundamental solution
$\cal{P}(\xi_0,\gamma)$ to be continuous at the node.
(This seems natural and is justified by the theoretical results of Appendix \ref{appendix:neck}:
see Remark \ref{remark:neck}.)

The spherical and catenoidal potentials both take value in $\Lambda\su(2)$ when $z\in\S^1$. So if all points $u_{jk}$ are on the unit circle,
the fundamental solution $\cal{P}(\xi_0,\gamma)$ will be in $\Lambda SU(2)$ for all $\gamma\in\pi_1(\overline{\Sigma}_0,z_0)$.
Unitarization is the hard task in solving the Monodromy Problem, so this explains why we restrict to horizontal planar graphs $\Gamma$.

The strategy of the construction is the following:
for small $t\neq 0$, we define a genuine Riemann surface $\overline{\Sigma}_t$ by opening the nodes of $\overline{\Sigma}_0$.
We define a meromorphic potential $\xi_t$ on $\overline{\Sigma}_t$ as a perturbation of the above potential $\xi_0$, depending on some parameters.
These parameters are determined by solving the Regularity and Monodromy Problems by
an implicit function argument at $t=0$.
\subsection{Symmetry}
In all the paper, $\sigma(z)=1/\overline{z}$ denotes the inversion with respect to the unit circle. The potentials $\xiS$ and $\xiC$ both have the symmetry
\begin{equation}
\label{eq:symmetry}
\overline{\sigma^*\xi}=D\xi D^{-1}\quad\mbox{ with }\quad D=\matrix{i&0\\0&-i}.
\end{equation}
A potential having the symmetry \eqref{eq:symmetry} will be called $\sigma$-symmetric.
With appropriate initial condition, the solution of the Cauchy Problem \eqref{eq:cauchy-problem}
satisfies
\begin{equation}
\label{eq:symmetry-Phi}
\overline{\sigma^*\Phi}=D\Phi D^{-1}.
\end{equation}
The corresponding surface is invariant by the isometry $X\mapsto D\overline{X}D^{-1}$ in the
$\mathfrak{su}(2)$-model, which corresponds to the symmetry with respect to the plane $x_1=0$.
We keep the $\sigma$-symmetry throughout the construction, and in the end apply the rigid motion $\Psi$
so that the surface is symmetric with respect to the horizontal plane $x_3=0$.
\section{Opening nodes}
\label{section:openingnodes}
In this section, we define a family of Riemann surfaces $\overline{\Sigma}_{t,\x}$ depending on a small real parameter $t$ and a certain number of other parameters, which we denote $\x$. We start by defining the Riemann surface with nodes $\overline{\Sigma}_{0,\x}$.
We proceed as in Section \ref{section:strategy} except that the position of the nodes in $\CC_j$ become parameters. (We can fix the nodes at $1$ and $-1$ in $\CC_{jk}$ by a M\"obius transformation.)
Consider a copy $\CC_j$ of the Riemann sphere for
$j\in J$ and a copy $\CC_{jk}$ of the Riemann sphere for $(j,k)\in E^+$.
For $(j,k)\in E^+$, introduce two complex parameters $p_{jk}$ and $p_{kj}$ 
in a neighborhood of respectively $u_{jk}$ and $u_{kj}$.
It will be convenient to denote $p_{jk}'=1$ and $p_{kj}'=-1$ the nodes in $\CC_{jk}$.
Identify the point $z=p_{jk}$ in $\CC_j$ with the point $z=p_{jk}'$ in $\CC_{jk}$ and the point $z=p_{kj}$ in $\CC_{k}$ with the point $z=p_{kj}'$ in $\CC_{jk}$ to create two nodes per edge.
This defines a compact Riemann surface with nodes denoted $\overline{\Sigma}_{0,\x}$.

To open nodes for $t\neq 0$, 
we introduce local complex coordinates in a neighborhood of $p_{jk}$ and $p_{jk}'$ for $(j,k)\in E$:
$$\z_{jk}=-2i\frac{z-p_{jk}}{z+p_{jk}} : V_{jk}\subset \C_j\stackrel{\sim}{\longrightarrow} D(0,\varepsilon).$$
$$\z_{jk}'=-2i\frac{z-p_{jk}'}{z+p_{jk}'}:V_{jk}'\subset\CC_{jk}\stackrel{\sim}{\longrightarrow}D(0,\varepsilon).$$
(These coordinates are chosen so that $\overline{\Sigma}_{t,\x}$ has the desired symmetry: see Proposition \ref{proposition:symmetrySigma}.)
We assume that $\varepsilon>0$ is small enough so that the disks $V_{jk}$ for $k\in E_j$
are disjoint.
For $(j,k)\in E$, we introduce a non-zero real parameter $r_{jk}$ in a neighborhood of ${\tau}_{jk}$ and set $t_{jk}=r_{jk}t$.
Assume that $t$ is small enough so that $|t_{jk}|<\varepsilon^2$.
Remove the disks $|\z_{jk}|\leq|t_{jk}|/\varepsilon$ and
$|\z_{jk}'|\leq|t_{jk}|/\varepsilon$.
Identify each point $z$ in the annulus $|t_{jk}|/\varepsilon<|\z_{jk}|<\varepsilon$
with the point $z'$ in the annulus $|t_{jk}|/\varepsilon<|\z_{jk}'|<\varepsilon$ such that
$$\z_{jk}(z)\z_{jk}'(z')=t_{jk}.$$
In particular, the circle $|\z_{jk}|=|t_{jk}|^{1/2}$ is identified with the circle
$|\z_{jk}'|=|t_{jk}|^{1/2}$, with the reverse orientation.
This creates two necks per edge.
The resulting compact Riemann surface is denoted $\overline{\Sigma}_{t,\x}$.
Note that it does not depend on $\lambda$.
The points $z=0$, $z=1$ and $z=\infty$ in $\CC_j$ are denoted respectively
$0_j$, $1_j$ and $\infty_j$. The points $z=0$ and $z=\infty$ in $\CC_{jk}$
are denoted $0_{jk}$ and $\infty_{jk}$.
\begin{remark}
\label{remark:epsilon}
The Riemann surface $\overline{\Sigma}_{t,\x}$ does not depend on the number
$\varepsilon>0$ used to define the domains $V_{jk}$, but the smaller $\varepsilon$, the smaller
$t$ must be since we need $|t_{jk}|<\varepsilon^2$.
\end{remark}
\subsection{Symmetry}
\begin{proposition}
\label{proposition:symmetrySigma}
Assume that $p_{jk}\in\S^1$ for all $(j,k)\in E$. Then $\overline{\Sigma}_{t,\x}$ admits an anti-holomorphic
involution $\sigma$ defined by
$\sigma(z)=1/\overline{z}$ in $\CC_j$ for $j\in J$ and
 $\CC_{jk}$ for $(j,k)\in E^+$.
 \end{proposition}

Proof: a straightforward computation gives, for $p_{jk}\in\S^1$ 
$$\overline{\z_{jk}(1/\overline{z})}=\z_{jk}(z).$$
A similar relation holds for $\z_{jk}'$. Hence since $t_{jk}$ is real,
$$\z_{jk}(z)\z_{jk}'(z')=t_{jk}\quad\Rightarrow\quad \z_{jk}(\sigma(z))\z_{jk}'(\sigma(z'))=t_{jk}.$$
So if $z\sim z'$ in $\overline{\Sigma}_{t,\x}$, then $\sigma(z)\sim\sigma(z')$ in $\overline{\Sigma}_{t,\x}$.\cqfd
\subsection{Meromorphic 1-forms on $\overline{\Sigma}_{t,\x}$.}
We denote $C(p_{jk})$ the circle $|\z_{jk}|=\varepsilon$ and 
$C(p'_{jk})$ the circle $|\z_{jk}'|=\varepsilon$.
Assume $t\neq 0$ and
let $\omega$ be a meromorphic 1-form on $\overline{\Sigma}_{t,\x}$ with poles outside of the annuli
$|t_{jk}|/\varepsilon<|\z_{jk}|<\varepsilon$.
We have
\begin{equation}
\label{eq:homology1}
\int_{C(p_{jk})}\omega=
-\int_{C(p_{jk}')}\omega
\quad\mbox{ for $(j,k)\in E$}
\end{equation}
because $C(p_{jk})$ is homologous to $-C(p_{jk}')$ in $\overline{\Sigma}_{t,\x}$.
By the Residue Theorem in $\C_j$
\begin{equation}
\label{eq:homology2}
\sum_{k\in E_j}\int_{C(p_{jk})}\omega+2\pi i\sum_{q\in \C_j}\Res_q\omega=0
\quad\mbox{ for $j\in J$}
\end{equation}
where the sum is taken on all poles $q$ of $\omega$ in $\C_j$.
In the same way,
\begin{equation}
\label{eq:homology3}
\int_{C(p_{jk}')}\omega+\int_{C(p_{kj}')}\omega+2\pi i\sum_{q\in \C_{jk}}\Res_q\omega=0
\quad\mbox{ for $(j,k)\in E^+$.}
\end{equation}
\begin{definition}[Bers]
A regular differential on the Riemann surface with nodes $\overline{\Sigma}_{0,\x}$ is a meromorphic 1-form with simple poles at the nodes
$p_{jk}$ and $p_{jk}'$ for $(j,k)\in E$, with opposite residues, and possibly poles of arbitrary order away from the nodes.
\end{definition}
\begin{theorem}
\label{theorem:openingnodes}
A meromorphic 1-form $\omega$ on $\overline{\Sigma}_{t\neq 0,\x}$
(respectively a regular differential $\omega$ on $\overline{\Sigma}_{0,\x}$) is uniquely defined by prescribing its
poles, principal parts at the poles and periods on the circles $C(p_{jk})$ and $C(p_{jk}')$ for
$(j,k)\in E$, subject only to the constraints \eqref{eq:homology1},
\eqref{eq:homology2} and \eqref{eq:homology3}.
Moreover, away from the nodes and the poles,
$\omega$ depends holomorphically on $t$ in a neighborhood of $0$ and
all parameters in the construction.
\end{theorem}
This is proved for holomorphic 1-forms in \cite{fay} and for meromorphic 1-forms
with simple poles in \cite{masur} using algebraic-geometric methods.
A proof for poles of arbitrary order is given in \cite{crelle}.
The holomorphic dependence away from the nodes and the poles means the following: for
$\epsilon>0$, let $\Omega_{\epsilon}$ be $\overline{\Sigma}_{0,\x}$ minus $\epsilon$-neighborhoods of
all nodes and poles, so $\Omega_{\epsilon}\subset \overline{\Sigma}_{t,\x}$ for $t$ small enough.
Then the restriction of $\omega$ to the fixed domain $\Omega_{\epsilon}$ depends holomorphically on $(z,t,\x)$.

\section{The potential}
\label{section:potential}
In this section, we define a meromorphic potential $\xi_{t,\x}$ on $\overline{\Sigma}_{t,\x}$, with
poles at the following points:
\begin{itemize}
\item $0_j$ and $\infty_j$ in $\CC_j$ for $j\in J$, which are to be apparent singularities,
\item $p_{jk}$ in $\CC_j$ for $(j,k)\in R$, which are to be the Delaunay ends of our surface. Here $p_{jk}$ is a $\lambda$-dependent parameter in the functional space $\Wp$ in a neighborhood of $u_{jk}$,
for $(j,k)\in R$.
\item $q_{jk}$ and $\sigma(q_{jk})$ in $\CC_{jk}$, for $(j,k)\in E^+$, which are to be apparent singularities. Here $q_{jk}$ is a
$\lambda$-dependent parameter in $\Wp$ in a neighborhood of $0$, for $(j,k)\in E^+$.
\end{itemize}
\begin{remark} All these $\lambda$-dependent parameters will be used to solve the Monodromy Problem.
The cross-ratio of $1$, $-1$, $q_{jk}$ and $\sigma(q_{jk})$ is
$$(1,-1;q_{jk},\sigma(q_{jk}))=\frac{q_{jk}\overline{q_{jk}}-1+2i\,\Im(q_{jk})}{1-q_{jk}\overline{q_{jk}}+2i\,\Im(q_{jk})}.$$
The derivative of the cross-ration with respect to $\Re(q_{jk})$ at $q_{jk}=0$ is zero,
so $\Re(q_{jk})$ serves no purpose and we restrict $q_{jk}$ to the space $i\WRp$.
We could have fixed the singularities at $0_{jk}$ and $\infty_{jk}$ and perturbed the position of the nodes at $1$ and $-1$, but then $\overline{\Sigma}_{t,\x}$ would depend on $\lambda$.
We chose to have a constant Riemann surface and moving singularities (with respect to $\lambda$), which
is more conventional than the reverse.
\end{remark}
We define the meromorphic potential $\xi_{t,\x}$ on $\overline{\Sigma}_{t,\x}$ as the sum of two terms:
$$\xi_{t,\x}=\eta_{t,\x}+t\,\chi_{t,\x}$$
where the potential $\eta_{t,\x}$ is a perturbation
of the potential $\xi_0$ described in Section \ref{section:strategy},
while the potential $\chi_{t,\x}$ prescribes periods around the nodes and suitable singularities at the Delaunay ends. These potentials are defined as follows, using Theorem \ref{theorem:openingnodes}:
\begin{itemize}
\item The potential $\eta_{t,\x}$ has simple poles at $0_j$ and $\infty_j$ for $j\in J$ with residues
$$\Res_{0_j}\eta_{t,\x}=-\Res_{\infty_j}\eta_{t,\x}=M_j=\matrix{iA_j&\lambda^{-1}B_j\\\lambda C_j&-iA_j},$$
simple poles at
$q_{jk}$ and $\sigma(q_{jk})$ for $(j,k)\in E^+$ with residues
$$\Res_{q_{jk}}\eta_{t,\x}=-\Res_{\sigma(q_{jk})}\eta_{t,\x}=M_{jk}=\matrix{iA_{jk}& B_{jk}\\C_{jk}&-iA_{jk}}$$
and has vanishing periods around the nodes:
$$\int_{C(p_{jk})}\eta_{t,\x}=\int_{
C(p_{jk}')}\eta_{t,\x}=0\quad\mbox{ for $(j,k)\in E$}.$$
Here $A_j$, $B_j$, $C_j$, $A_{jk}$, $B_{jk}$, $C_{jk}$ are parameters in
a neighborhood of respectively $0$, $1/2$, $1/2$, $0$, $1/2$, $1/2$ in
$\WRp$.
\item
The potential $\chi_{t,\x}$ has the following periods around the nodes for $(j,k)\in E$:
$$\int_{C(p_{jk})}\chi_{t,\x}=-\int_{C(p_{jk}')}
\chi_{t,\x}=2\pi i \,m_{jk}\quad\mbox{ with }\quad m_{jk}=\matrix{a_{jk}&\lambda^{-1}ib_{jk}\\ic_{jk}&-a_{jk}}$$
where $a_{jk}$, $b_{jk}$, $c_{jk}$ for $(j,k)\in E$ are parameters in $\WRp$ to be determined.
It has a double pole at $p_{jk}$ in $\C_j$ for $(j,k)\in R$ with principal part
$$\matrix{0&0\\1&0}\left(\frac{a_{jk}p_{jk}dz}{(z-p_{jk})^2}+\frac{ib_{jk}dz}{z-p_{jk}}\right).$$
Here $a_{jk}$, $b_{jk}$ are parameters in $\WRp$ to be determined,  for $(j,k)\in R$.
It is known from \cite{nnoids} that such a pole creates a Delaunay end, provided the Monodromy Problem is solved.
Finally, the potential $\chi_{t,\x}$ has simple poles with equal residues at $0_j$ and $\infty_j$ and
simples poles with equal residues at
$q_{jk}$ and $\sigma(q_{jk})$.
These residues are determined by the constraints \eqref{eq:homology2} and
\eqref{eq:homology3} which give:
\begin{equation}
\Res_{0_j}\chi_{t,\x}=\Res_{\infty_j}\chi_{t,\x}=-\frac{1}{2}\sum_{k\in E_j}m_{jk}-\frac{1}{2}\sum_{k\in R_j}\matrix{0&0\\i b_{jk}&0}
\end{equation}
\begin{equation}
\Res_{q_{jk}}\chi_{t,\x}=\Res_{\sigma(q_{jk})}\chi_{t,\x}
=\frac{1}{2}(m_{jk}+m_{kj}).
\end{equation}
\end{itemize}
\subsection{Symmetry}
The residues and periods of the entries of $\eta_{t,\x}$ and $\chi_{t,\x}$ have been chosen to be either real or imaginary so that the potential has the desired symmetry:
\begin{proposition}
\label{proposition:symmetryxi}
Assume that $p_{jk}\in\S^1$ for $(j,k)\in E$ and
$p_{jk}=e^{i\theta_{jk}}$ with $\theta_{jk}\in\WRp$ for $(j,k)\in R$.
Then the potential $\xi_{t,\x}$ has the symmetry \eqref{eq:symmetry}:
$$\overline{\sigma^*\xi_{t,\x}}=D\xi_{t,\x}D^{-1}.$$
\end{proposition}
Note that the bar denotes the conjugation operator defined in Section \ref{section:functional-spaces}
so this actually means
$\overline{\sigma^*\xi_{t,\x}(z,\overline{\lambda})}=D\xi_{t,\x}(z,\lambda)D^{-1}$.
Both sides are holomorphic with respect to $\lambda$.

Proof: if $\omega$ is a meromorphic 1-form on $\overline{\Sigma}_{t,\x}$ then $\overline{\sigma^*\omega}$ is meromorphic and
$$\Res_{\sigma(p)}\overline{\sigma^*\omega}=\overline{\Res_p\omega}.$$
Hence
$$\Res_{\infty_j}\overline{\sigma^*\eta_{t,\x}}
=\overline{\Res_{0_j}\eta_{t,\x}}
=\overline{M_j}=-DM_jD^{-1}=\Res_{\infty_j}D\eta_{t,\x}D^{-1}.$$
In the same way, $\overline{\sigma^*\eta_{t,\x}}$ and $D\eta_{t,\x}D^{-1}$ have the same
residues at $0_j$,
$q_{jk}$ and $\sigma(q_{jk})$.
 Moreover, both have vanishing periods
around the nodes, so by uniqueness in Theorem \ref{theorem:openingnodes},
$$\overline{\sigma^*\eta_{t,\x}}=D\eta_{t,\x}D^{-1}.$$
For $(j,k)\in E$, we have since $\sigma(C(p_{jk}))=-C(p_{jk})$
$$\int_{C(p_{jk})}\overline{\sigma^*\chi_{t,\x}}=-\int_{C(p_{jk})}\overline{\chi_{t,\x}}
=2\pi i\, \overline{m_{jk}}=2\pi i D m_{jk}D^{-1}$$
so $\overline{\sigma^*\chi_{t,\x}}$ and $D\chi_{t,\x}D^{-1}$ have the same periods around the nodes. 
For $(j,k)\in R$, assuming that $p_{jk}=e^{i\theta_{jk}}$ with $\theta_{jk}\in \WRp$, we have
$\sigma(p_{jk})=p_{jk}$ and
$$\overline{\sigma^*\left(\frac{a_{jk}p_{jk}dz}{(z-p_{jk})^2}+\frac{ib_{jk}dz}{z-p_{jk}}\right)}
=-\frac{a_{jk}p_{jk}dz}{(z-p_{jk})^2}-\frac{ib_{jk}dz}{z-p_{jk}}+\frac{ib_{jk}dz}{z}$$
so $\overline{\sigma^*\chi_{t,\x}}$ and $D\chi_{t,\x}D^{-1}$ have the same principal part at
$p_{jk}$.
Finally, they have the same residues at $0_j$, $\infty_j$, $q_{jk}$ and $\sigma(q_{jk})$ by computations similar to the above. By uniqueness in Theorem \ref{theorem:openingnodes}, we have
$$\overline{\sigma^*\chi_{t,\x}}=D\chi_{t,\x}D^{-1}.$$
\cqfd
\subsection{Explicit formulas at $t=0$}
It will be convenient to denote, for $(j,k)\in E^+$, $M_{kj}=M_{jk}$,
$A_{kj}=A_{jk}$, etc... so $M_{jk}$ makes sense for all $(j,k)\in E$.
Be careful however that $A_{jk}$ and $A_{kj}$ are the same parameter, whereas $a_{jk}$ and
$a_{kj}$ are distinct parameters.
For complex numbers $p$, $q$, we denote $\omega_q$
the meromorphic 1-form on the Riemann sphere with simple poles at $q$ and $\sigma(q)$ with residues $1$ and $-1$, and
$\omega_{p,q}$ the meromorphic 1-form with simple poles at $p$, $q$ and $\sigma(q)$ with residues $1$, $-1/2$ and $-1/2$.
Explicitly:
$$\omega_q=\frac{dz}{z-q}-\frac{\overline{q}dz}{\overline{q}z-1}
=\frac{(1-q\overline{q})dz}{(z-q)(1-\overline{q}z)}
\quad\mbox{ and }\quad
\omega_{p,q}=\frac{dz}{z-p}-\frac{dz}{2(z-q)}-\frac{\overline{q}\,dz}{2(\overline{q}z-1)}.$$
In particular if $q=0$:
$$\omega_0=\frac{dz}{z}\quad\mbox{ and }\quad
\omega_{p,0}=\frac{dz}{z-p}-\frac{dz}{2z}.$$
\begin{proposition}
\label{proposition:time0}
At $t=0$ and for any value of the parameter $\x$,
we have in $\CC_j$ for $j\in J$:
$$\eta_{0,\x}=M_j\omega_0$$
$$\chi_{0,\x}=\sum_{k\in E_j}m_{jk}\omega_{p_{jk},0}
+\sum_{k\in R_j}\matrix{0&0\\1&0}\left(\frac{a_{jk}p_{jk}\,dz}{(z-p_{jk})^2}+i b_{jk}\omega_{p_{jk},0}\right)$$
and in $\CC_{jk}$ for $(j,k)\in E^+$:
$$\eta_{0,\x}=M_{jk}\omega_{q_{jk}}$$
$$\chi_{0,\x}=-m_{jk}\omega_{1,q_{jk}}-m_{kj}\omega_{-1,q_{jk}}.$$
\end{proposition}
Proof: the entries of $\eta_{0,\x}$ and $\chi_{0,\x}$ are regular meromorphic differentials
on the Riemann surface with nodes $\overline{\Sigma}_{0,\x}$. Proposition \ref{proposition:time0} follows from the fact that a meromorphic 1-form on the Riemann sphere is uniquely defined by its poles and principal parts.
\cqfd

We shall need the $t$-derivative of the potential $\xi_{t,\x}$ at $t=0$. We have of course
$$\frac{\partial\xi_{t,\x}}{\partial t}\mid_{t=0}=\frac{\partial\eta_{t,\x}}{\partial t}\mid_{t=0}+\chi_{0,\x}.$$
\begin{proposition}
\label{proposition:derivee0}
The $t$-derivative of the potential $\eta_{t,\x}$ at $t=0$ is given by
$$\frac{\partial\eta_{t,\x}}{\partial t}\mid_{t=0}=\left\{\begin{array}{ll}
\displaystyle\sum_{k\in E_j}r_{jk}M_{jk}\frac{(1+q_{jk}^2)}{(1-q_{jk}^2)}\frac{p_{jk}dz}{(z-p_{jk})^2}&\mbox{ in $\CC_j$ for $j\in J$}\\
\displaystyle
r_{jk}M_j\frac{dz}{(z-1)^2}
-r_{kj}M_k\frac{dz}{(z+1)^2}&\mbox{ in $\CC_{jk}$ for $(j,k)\in E^+$.}
\end{array}\right.$$
\end{proposition}
Proof: by Lemma 3 in \cite{triply}, for $(j,k)\in E$, the derivative of $\eta_{t,\x}$ with respect to the parameter $t_{jk}$
at $t=0$ is a meromorphic differential on $\overline{\Sigma}_{0,\x}$ with two double poles at $p_{jk}$,
$p'_{jk}$ and principal parts given in term of the coordinates used to open nodes by
\begin{equation}
\label{eq:principal-parts}
\frac{\partial\eta_{t,\x}}{\partial t_{jk}}\mid_{t=0}
\simeq\left\{\begin{array}{ll}
\displaystyle\frac{-d\z_{jk}}{(\z_{jk})^2}\Res_{p'_{jk}}\frac{\eta_{0,\x}}{\z'_{jk}}&\mbox{ at $p_{jk}$}\\
\displaystyle\frac{-d\z'_{jk}}{(\z'_{jk})^2}\Res_{p_{jk}}\frac{\eta_{0,\x}}{\z_{jk}}&\mbox{ at $p'_{jk}$}
\end{array}\right.
.\end{equation}
We have
$$\frac{d\z_{jk}}{(\z_{jk})^2}=\frac{i p_{jk}\,dz}{(z-p_{jk})^2}\quad\mbox{ and }\quad
\frac{d\z_{jk}'}{(\z_{jk}')^2}=\frac{i p'_{jk}\,dz}{(z-p'_{jk})^2}.$$
Observe that these are globally defined meromorphic 1-forms on the Riemann sphere so $\simeq$
in \eqref{eq:principal-parts} becomes an equality in $\CC_j$ and $\CC_{jk}$, respectively.
By Proposition \ref{proposition:time0}:
$$\Res_{p_{jk}}\frac{\eta_{0,\x}}{\z_{jk}}=\Res_{p_{jk}}\frac{(z+p_{jk})}{-2i(z-p_{jk})}M_j\frac{dz}{z}=
i M_j$$
Recalling that $q_{jk}\in i\WRp$ so $\overline{q_{jk}}=-q_{jk}$ and that $p_{jk}'=\pm 1$:
$$\Res_{p'_{jk}}\frac{\eta_{0,\x}}{\z'_{jk}}=\Res_{p'_{jk}}\frac{(z+p'_{jk})}{-2i(z-p'_{jk})}M_{jk}\frac{(1+q_{jk}^2)dz}{(z-q_{jk})(1+q_{jk}z)}
=iM_{jk}\frac{1+q_{jk}^2}{1-q_{jk}^2}.$$
Hence for $(j,k)\in E$:
$$\frac{\partial\eta_{t,\x}}{\partial t_{jk}}\mid_{t=0}=\left\{\begin{array}{ll}
\displaystyle M_{jk}\frac{(1+q_{jk}^2)}{(1-q_{jk}^2)}\frac{p_{jk}\,dz}{(z-p_{jk})^2}&\mbox{ in $\CC_j$}\\
\displaystyle M_j\frac{p'_{jk}dz}{(z-p'_{jk})^2}&\mbox{ in $\CC_{jk}$}\\
0 &\mbox{ in all other Riemann spheres.}
\end{array}\right.$$
Proposition \ref{proposition:derivee0} follows from $t_{jk}=r_{jk}t$ and the chain rule.
\cqfd
\subsection{Central value of the parameters}
\label{section:centralvalue}
The vector of all parameters of the construction (except $t$) is denoted $\x$.
Each parameter is in a neighborhood of a central value denoted with an underscore.
The central values are tabulated below. Some of them we have already seen. The others will be computed when solving the Monodromy Problem.

Also, we have tried to define the potential in a way as general and natural as possible, but it turns out
{\em a posteriori} after solving all equations that we have too many parameters, so we can fix the value of some of them: $A_j$, $B_j$ for
$j\in J$ will not be used.
Some computations are simpler with these restrictions so we assume them from now on.
\medskip
$$\begin{array}{|c|c|c|c|}
\hline \mbox{\small parameter}&\mbox{\small range}&\mbox{\small space}&\mbox{\small central value}\\
\hline
p_{jk}&(j,k)\in E&\S^1&u_{jk}\\
p'_{jk}&(j,k)\in E&\mbox{fixed}&\pm 1\\
r_{jk}&(j,k)\in E&\R&\tau_{jk}\\
\hline
A_j&j\in J&\mbox{fixed}&0\\
B_j&j\in J&\mbox{fixed}&1/2\\
C_j&j\in J&\WRp&1/2\\
\hline
q_{jk}&(j,k)\in E^+&i\WRp&0\\
A_{jk}&(j,k)\in E^+&\WRp&0\\
 B_{jk},C_{jk}&(j,k)\in E^+&\WRp&1/2\\
\hline
 a_{jk}&(j,k)\in E&\WRp&\tau_{jk}(\lambda-1)/2\\
b_{jk},c_{jk}&(j,k)\in E&\WRp&0\\
\hline
p_{jk}&(j,k)\in R&\exp(i\WRp)&u_{jk}\\
a_{jk}&(j,k)\in R&\WRp&\tau_{jk}(\lambda-1)^2/2\\
b_{jk}&(j,k)\in R&\WRp&0\\
\hline
\end{array}$$
\section{The Regularity Problem}
\label{section:regularity}
We want $0_j$, $\infty_j$ and $q_{jk}$, $\sigma(q_{jk})$ to be apparent singularities.
In this section, the entries of the potential will be denoted
\begin{equation}
\label{eq:alphabetagamma}
\xi_{t,\x}=\matrix{\alpha&\lambda^{-1}\beta\\\gamma&-\alpha}
\end{equation}
and the dependence on the parameters $(t,\x)$ will not be written to ease notations.
\subsection{Regularity at $0_j$ and $\infty_j$}
\label{section:regularity0j}
Fix $j\in J$ and consider the gauge
$$G_j=\matrix{f&0\\ \lambda g&f^{-1}}\quad\mbox{ with }\quad
f(z)=\frac{1+z}{\sqrt{z}}\quad\mbox{ and }\quad
g(z)= x_j\frac{1-z}{\sqrt{z}}+i\,y_j\frac{1+z}{\sqrt{z}}.$$
Here $x_j$, $y_j$ are parameters in $\WRp$ to be determined.
At $(x_j,y_j)=(1,0)$ we have $G_j=\GS$.
We denote
$$\wh{\xi}=\xi_{t,\x}\cdot G_j=\matrix{\wh{\alpha}&\lambda^{-1}\wh{\beta}\\\wh{\gamma}&-\wh{\alpha}}.$$
The gauge has the symmetry $\overline{G_j\circ\sigma}=D G_j D^{-1}$ so $\wh{\xi}$ has the symmetry \eqref{eq:symmetry} and it suffices to ensure that $\wh{\xi}$ is regular at $0_j$ ; regularity at
$\infty_j$ will follow by symmetry.
\begin{proposition}
\label{proposition:regularity0j}
There exists explicit values of $x_j$, $y_j$ and $C_j$ in $\WRp$,
depending analytically on $(t,\x)$,
such that
$\wh{\alpha}$ and $\wh{\beta}$ are holomorphic at $0_j$,
$\wh{\gamma}$ has a pole of multiplicity at most 2 and
\begin{equation}
\label{eq:ReReswhgamma}
\Re\left(\Res_{0_j}(z\wh{\gamma})\right)=0.
\end{equation}
\end{proposition}
Proof: straightforward computations give
$$\wh{\alpha}=\alpha+f^{-1}g\beta+f^{-1}df$$
$$\wh{\beta}=f^{-2}\beta$$
\begin{equation}
\label{eq:whgamma}
\wh{\gamma}=-2\lambda fg\alpha-\lambda g^2\beta+f^2\gamma+\lambda(f\,dg-g\,df).
\end{equation}
Recall that $\alpha$, $\beta$, $\gamma$ have simple poles at $0_j$. Hence
$\wh{\beta}$ is holomorphic at $0_j$ and $\wh{\alpha}$ has (at most) a simple pole with residue
$$\Res_{0_j}\wh{\alpha}=\Res_{0_j}\alpha+(x_j+iy_j)\Res_{0_j}\beta-\frac{1}{2}.$$
We take
\begin{equation}
\label{eq:xjyj}
x_j+iy_j=\frac{1/2-\Res_{0_j}\alpha}{\Res_{0_j}\beta}
\end{equation}
so that $\wh{\alpha}$ is holomorphic at $0_j$.
Finally, $\wh{\gamma}$ has at most a double pole at $0_j$ and since $f\,dg-g\,df$ has a simple pole at $0$,
$$\Res_{0_j}(z\wh{\gamma})=-2\lambda(x_j+iy_j)\Res_{0_j}\alpha-\lambda(x_j+iy_j)^2\Res_{0_j}\beta
+\Res_{0_j}\gamma.$$
By definition, recalling the definition of the operator $\Re$ in Section \ref{section:functional-spaces}:
$$\Re(\Res_{0_j}\gamma)=\lambda C_j.$$
So we see that Equation \eqref{eq:ReReswhgamma} is equivalent to
$$C_j=\Re\left[2(x_j+iy_j)\Res_{0_j}\alpha+(x_j+iy_j)^2\Res_{0_j}\beta\right]$$
which using Equation \eqref{eq:xjyj} simplifies to
\begin{equation}
\label{eq:Cj}
C_j=\Re\left(\frac{1/4-(\Res_{0_j}\alpha)^2)}{\Res_{0_j}\beta}\right).
\end{equation}
Note that the residues of $\alpha$ and $\beta$ involved in Equations \eqref{eq:xjyj} and \eqref{eq:Cj}
are given, as functions of $(t,\x)$, by the definition of $\xi_{t,\x}$.
In particular, at $t=0$, we have $x_j=1$ and $y_j=0$ so $G_j=\GS$, and
\begin{equation}
\label{eq:Cjt0}C_j\mid_{t=0}=\frac{1/4+A_j^2}{B_j}=1/2.
\end{equation}
\cqfd

At this point, the Regularity Problem at $0_j$ is only partially solved since $\wh{\gamma}$ still has a pole.
By Equation \eqref{eq:whgamma}, we have
\begin{equation}
\label{eq:whgamma0}
\wh{\gamma}^0=z^{-1}(z+1)^2\gamma^0
\end{equation}
so for $\wh{\gamma}$ to be holomorphic at $0_j$, it is necessary that
$$\Res_{0_j}\left(z^{-1}(z+1)^2\gamma^0\right)=0.$$
We define for $j\in J$ and $t\neq 0$
\begin{equation}
\label{eq:defRj}
\cal{R}_j(t,\x)=t^{-1}\Res_{0_j}\left(z^{-1}(z+1)^2\gamma_{t,\x}^0\right)\in\C.
\end{equation}
\begin{proposition}
\label{proposition:Rj}
For $j\in J$, the function $\cal{R}_j(t,\x)$ extends analytically at $t=0$. Moreover, at the central value,
we have $\cal{R}_j(0,\cv{\x})=\overline{F_j}/2$, where $F_j$ is the force defined in
Equation \eqref{eq:force}.
\end{proposition}
Proof: by Proposition \ref{proposition:time0}, we have $\gamma_{0,\x}=\lambda C_j \omega_0$ in $\CC_{j}$ so $\gamma_{0,\x}^0=0$.
Hence $\cal{R}_j$ extends analytically at $t=0$ and
$$\cal{R}_j(0,\x)=\Res_{0_j}\left(z^{-1}(z+1)^2\frac{\partial\gamma_{t,\x}^0}{\partial t}\mid_{t=0}\right).$$
By Proposition \ref{proposition:derivee0}, we have
$$\frac{\partial\gamma_{t,\x}}{\partial t}\mid_{t=0}
=\sum_{k\in E_j}
\left(r_{jk}C_{jk}\frac{(1+q_{jk}^2)}{(1-q_{jk}^2)}\frac{p_{jk}\,dz}{(z-p_{jk})^2}
+i c_{jk}\omega_{p_{jk},0}\right)
+\sum_{k\in R_j}\left(\frac{a_{jk}p_{jk}\,dz}{(z-p_{jk})^2}+i b_{jk}\omega_{p_{jk},0}\right).$$
At the central value (see the table in Section \ref{section:centralvalue}) and $\lambda=0$, this simplifies to
$$\frac{\partial\gamma^0_{t,\cv{\x}}}{\partial t}\mid_{t=0}
=\sum_{k\in E_j\cup R_j}\frac{\tau_{jk}u_{jk}\,dz}{2(z-u_{jk})^2}$$
which is holomorphic at $0_j$.
Hence
$$\cal{R}_j(0,\cv{\x})=\sum_{k\in E_j\cup R_j}\frac{\tau_{jk}}{2u_{jk}}.$$
\cqfd
\begin{remark}
\label{remark:regularity0j}
Proposition \ref{proposition:Rj} explains where the balancing condition comes from.
We solve the equation $\cal{R_j}=0$ in Section \ref{section:implicit} using the non-degeneracy hypothesis.
Then after the Monodromy Problem is solved, $\wh{\gamma}$ will in fact be holomorphic at $0_j$:
see Proposition \ref{proposition:spherical}.
\end{remark}

\subsection{Regularity at $q_{jk}$ and $\sigma(q_{jk})$}
\label{section:regularity-qjk}
Fix $(j,k)\in E^+$.
Recall that $\xi_{t,\x}$ has moving singularities at $q_{jk}$ and $\sigma(q_{jk})$, which depend on $\lambda$.
We use the following M\"obius transformation as local coordinate in a neighborhood of $q_{jk}$:
$$w_{jk}(z)=\frac{z-q_{jk}}{1-\overline{q_{jk}}z}=\frac{z-q_{jk}}{1+q_{jk}z}$$
We have $\sigma\circ w_{jk}=w_{jk}\circ\sigma$.
We make the change of variable $w=w_{jk}$ and denote
$$\wt{\xi}=(w_{jk}^{-1})^*\xi_{t,\x}=\matrix{\wt{\alpha}&\lambda^{-1}\wt{\beta}\\\wt{\gamma}&-\wt{\alpha}}$$
which has fixed singularities at $w=0$ and $w=\infty$ and still has the symmetry \eqref{eq:symmetry}.
We consider a gauge $G_{jk}$ of a form dual to $G_j$:
$$G_{jk}=\matrix{f^{-1}&g\\0&f}
\quad\mbox{ with }\quad
f=\frac{1+w}{\sqrt{w}}\quad\mbox{ and }\quad
g= x_{jk}\frac{1-w}{\sqrt{w}}+i\,y_{jk}\frac{1+w}{\sqrt{w}}.$$
Let
$$\wh{\xi}=\wt{\xi}\cdot G_{jk}=\matrix{\wh{\alpha}&\lambda^{-1}\wh{\beta}\\\wh{\gamma}&-\wh{\alpha}}.$$
The gauge $G_{jk}$ has the symmetry $\overline{G_{jk}\circ\sigma}=DG_{jk}D^{-1}$
so it suffices to ensure that $\wh{\xi}$ is regular at $w=0$,
regularity at $w=\infty$ will follow by symmetry.
\begin{proposition}
\label{proposition:regularity0jk}
There exists explicit values of $x_{jk}$, $y_{jk}$ and $B_{jk}$ in $\WRp$,
depending analytically on $(t,\x)$,
such that
$\wh{\alpha}$ and $\wh{\gamma}$ are holomorphic at $w=0$,
$\wh{\beta}$ has a pole of multiplicity at most 2 and
\begin{equation}
\label{eq:ReReswhbeta}
\Re\left(\Res_0(w\wh{\beta})\right)=0.
\end{equation}
\end{proposition}
Proof: we simply dualize the proof of Proposition \ref{proposition:regularity0j} with $\wt{\xi}$ in place of $\xi_{t,\x}$ and obtain:
\begin{equation}
\label{eq:xjkyjk}
x_{jk}+iy_{jk}=\frac{1/2+\Res_0\wt{\alpha}}{\Res_0\wt{\gamma}}
=\frac{1/2+\Res_{q_{jk}}\alpha}{\Res_{q_{jk}}\gamma}
\end{equation}
\begin{equation}
\label{eq:Bjk}
%B_{jk}=\Re\left(-2(x_{jk}+iy_{jk})\Res_{0_{jk}}\alpha+(x_{jk}+iy_{jk})^2\Res_{0_{jk}}\gamma\right).
B_{jk}=\Re\left(\frac{1/4-(\Res_0\wt{\alpha})^2}{\Res_0\wt{\gamma}}\right)
=\Re\left(\frac{1/4-(\Res_{q_{jk}}\alpha)^2}{\Res_{q_{jk}}\gamma}\right)
\end{equation}
\begin{equation}
\label{eq:Bjkt0}
B_{jk}\mid_{t=0}=\frac{1/4+A_{jk}^2}{C_{jk}}.
\end{equation}
\cqfd

At this point, the Regularity Problem at $q_{jk}$ is only partially solved since $\wh{\beta}$ still has a pole.
Dualizing Equation \eqref{eq:whgamma0} we have
\begin{equation}
\label{eq:whbeta0}
\wh{\beta}^0=w^{-1}(w+1)^2\wt{\beta}^0
\end{equation}
For $\wh{\beta}$ to be holomorphic, it is necessary that
$$\Res_0\left(w^{-1}(w+1)^2\wt{\beta}^0\right)=0.$$
We define for $(j,k)\in E^+$ and $t\neq 0$:
$$\cal{R}_{jk}(t,\x)=t^{-1}\Res_0\left(w^{-1}(w+1)^2\wt{\beta}^0\right)\in\C.$$
\begin{proposition}
\label{proposition:Rjk} For $(j,k)\in E^+$, the function $\cal{R}_{jk}$ extends analytically at $t=0$ and
\begin{equation}
\label{eq:Rjk}
\cal{R}_{jk}(0,\x)
%=\frac{r_{jk}}{2}+\frac{r_{kj}}{2q_{kj}^0}+2ib_{jk}^0+ib_{kj}^0(1+1/q_{kj}^0).
=r_{jk}\frac{1+(q_{jk}^0)^2}{2(1-q_{jk}^0)^2}-r_{kj}\frac{1+(q_{jk}^0)^2}{2(1+q_{jk}^0)^2}
+\frac{2ib_{jk}^0}{1-q_{jk}^0}+\frac{2iq_{jk}^0b_{kj}^0}{1+q_{jk}^0}.
\end{equation}
In particular, $\cal{R}_{jk}(0,\cv{\x})=0$ at the central value.
\end{proposition}
Proof: by Proposition \ref{proposition:time0}, we have $\beta_{0,\x}=\lambda B_{jk}\omega_{q_{jk}}$ in $\CC_{jk}$ so $\wt{\beta}_{0,\x}^0=0$.
Hence $\cal{R}_{jk}$ extends analytically at $t=0$ and
$$\cal{R}_{jk}(0,\x)=\Res_0\left(w^{-1}(w+1)^2\frac{\partial\wt{\beta}_{t,\x}^0}{\partial t}\mid_{t=0}\right).$$
By Proposition \ref{proposition:derivee0}, remembering that we fixed $B_j=B_k=1/2$:
$$\frac{\partial\beta_{t,\x}}{\partial t}\mid_{t=0}=\frac{r_{jk}\,dz}{2(z-1)^2}
-\frac{r_{kj}\,dz}{2(z+1)^2}-ib_{jk}\omega_{1,q_{jk}}-ib_{kj}\omega_{-1,q_{kj}}.$$
The first two residues are better computed using the $z$-coordinate
$$\Res_{w=0}\left(w^{-1}(w+1)^2(w_{jk}^{-1})^*\frac{dz}{(z\pm 1)^2}\right)
=\Res_{z=q_{jk}}\left(w_{jk}^{-1}(w_{jk}+1)^2\frac{dz}{(z\pm 1)^2}\right)
=\frac{1+q_{jk}^2}{(q_{jk}\pm 1)^2}.$$
The last two residues are better computed using the $w$-coordinate:
$$(w_{jk}^{-1})^*\omega_{\pm 1,q_{jk}}=\frac{dw}{w-w_{jk}(\pm 1)}-\frac{dw}{2w}$$
$$\Res_{0}\left(w^{-1}(w+1)^2(w_{jk}^{-1})^*\omega_{1,q_{jk}}\right)=\frac{-1}{w_{jk}(1)}-1=\frac{2}{q_{jk}-1}$$
$$\Res_{0}\left(w^{-1}(w+1)^2(w_{jk}^{-1})^*\omega_{-1,q_{jk}}\right)=\frac{-1}{w_{jk}(-1)}-1=\frac{-2q_{jk}}{q_{jk}+1}.$$
Collecting all terms and setting $\lambda=0$, we obtain Equation \eqref{eq:Rjk}.\cqfd
\begin{remark}
\label{remark:regularity0jk}
We solve the equation $\cal{R_{jk}}(t,\x)=0$ using the Implicit Function Theorem in Section \ref{section:implicit}.
Then after the Monodromy Problem is solved, $\wh{\beta}$ will in fact be holomorphic at
$w=0$:
see Proposition \ref{proposition:catenoid}.
\end{remark}
\section{The Monodromy Problem}
From now on, we assume that $C_j$ is given in function of $(t,\x)$ by Equation \eqref{eq:Cj} for $j\in J$
and $B_{jk}$ is given by Equation \eqref{eq:Bjk} for $(j,k)\in E^+$.
Also, we restrict $t$ to be positive.
\subsection{Definition of various paths}
\begin{center}
\begin{figure}
\includegraphics[height=5cm]{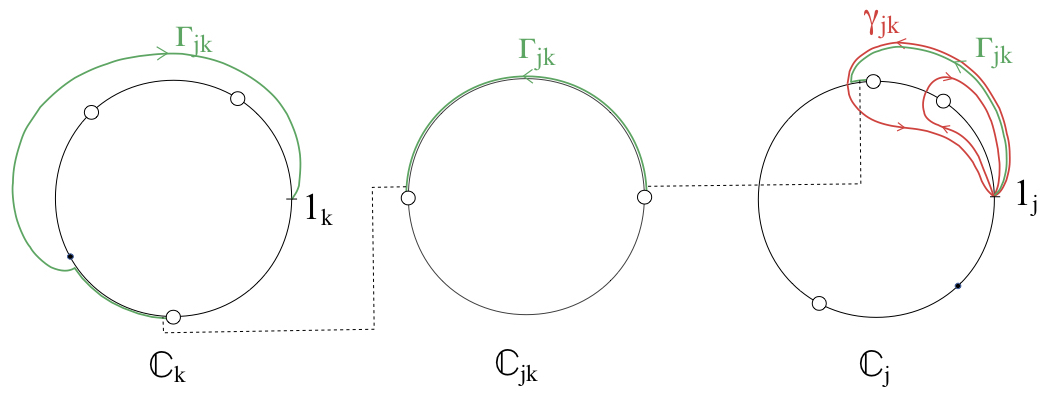}
\caption{The paths $\gamma_{jk}$ (red) and $\Gamma_{jk}$ (green), for $(j,k)\in E^+$.
The large circles represent the unit circles in $\CC_j$, $\CC_{jk}$ and $\CC_k$.
The tiny circles represent disks that are removed when opening nodes.
The bullets represent ends.
The dots connect points that are identified when opening nodes.}
\label{fig2}
\end{figure}
\end{center}

In this section, we define for $(j,k)\in E\cup R$ a loop $\gamma_{jk}$ with base point $1_j$ encircling the point $p_{jk}$, and for $(j,k)\in E^+$ a path $\Gamma_{jk}$ connecting $1_j$ to $1_k$ through the two necks
corresponding to the edge $(j,k)$
(see Figure \ref{fig2}). We study carefully how these paths transform under $\sigma$.

Fix $j\in J$.
We define an order $\prec$ on the set $E_j\cup R_j$ by
 $k\prec\ell\Leftrightarrow\arg(u_{jk})<\arg(u_{j\ell})$, where the arguments are chosen in
 $(0,2\pi)$.
For $k\in E_j\cup R_j$, we fix a curve $\alpha_{jk}$ in the domain $\{z\in\C_j:|z|>1,0<\arg(z)<2\pi\}$
from $1_j$ to $e^{i\varepsilon}u_{jk}$ and define
$\delta_{jk}=\alpha_{jk}\sigma(\alpha_{jk})^{-1}$.
The domain bounded by $\delta_{jk}$ contains the points $p_{j\ell}$ for $\ell\preceq k$.
We define inductively the loops $\gamma_{jk}$ for $k\in E_j\cup R_j$ by
\begin{equation}
\label{eq:induction-deltajk}
\delta_{jk}=\prod_{\ell\preceq k}\gamma_{j\ell}.
\end{equation}
In other words, $\gamma_{jk}=(\delta_{jk'})^{-1}\delta_{jk}$ where $k'$ is the predecessor of
$k$ for the order $\prec$.
The domain bounded by $\gamma_{jk}$ contains the point $p_{jk}$ and no other $p_{j\ell}$.
It will be convenient to also denote
$$\delta_{jk}'=\prod_{\ell\prec k}\gamma_{j\ell}$$
so $\delta_{jk}=\delta_{jk}'\gamma_{jk}$.
(An empty product means the neutral element.)
These paths transform as follows under $\sigma$:
\begin{equation}
\label{eq:symmetry-deltajk}
\sigma(\delta_{jk})=\delta_{jk}^{-1}
\end{equation}
\begin{equation}
\label{eq:symmetry-gammajk}
\sigma(\gamma_{jk})=\delta_{jk}'\gamma_{jk}^{-1}(\delta_{jk}')^{-1}.
\end{equation}
Fix $(j,k)\in E^+$. The path $\Gamma_{jk}$ is defined as follows.
Fix a number $\varepsilon'$ such that $0<\varepsilon'<\varepsilon$, where $\varepsilon$ is the
 number introduced to open nodes in section \ref{section:openingnodes}.
Recalling the definition of the coordinate $\z_{jk}$ near $p_{jk}$, we have
$$z=p_{jk}\frac{(2+i \z_{jk})}{(2-i\z_{jk})}$$
so for real $x\in[-\varepsilon,\varepsilon]$, the point $\z_{jk}=x$ is on the
unit circle and its argument is an increasing function of $x$.
First assume that $\tau_{jk}>0$ so $t_{jk}$ and $t_{kj}$ are positive.
We define the path $\beta_{jk}$ as the concatenation of the following 5 paths (taking care to avoid the disks that are removed when opening nodes):
\begin{enumerate}
\item The circular arc from $z=e^{i\varepsilon}u_{jk}$ to $\z_{jk}=\varepsilon'$.
\item The circular arc from $\z_{jk}=\varepsilon'$ to $\z_{jk}=t_{jk}/\varepsilon'$.
Its endpoint was identified with $\z_{jk}'=\varepsilon'$ when opening nodes.
\item The circular arc from $\z_{jk}'=\varepsilon'$ to $\z_{kj}'=-\varepsilon'$ on the upper half unit circle in
$\CC_{jk}$.
\item The circular arc from $\z_{kj}'=-\varepsilon'$ to $\z_{kj}'=-t_{kj}/\varepsilon'$.
Its endpoint was identified with $\z_{kj}=-\varepsilon'$ when opening nodes.
\item The circular arc from $\z_{kj}=-\varepsilon'$ to $z=e^{i\varepsilon}u_{kj'}$, where 
$j'$ is the predecessor of $j$ for the order $\prec$ on $E_k\cup R_k$
(or to $z=1_k$ in case $j$ is the minimum of $E_k\cup R_k$).
\end{enumerate}
(We could of course group paths (1) and (2) into one single arc, but it is convenient for the proof of Proposition \ref{proposition:Gamma-monodromy} to write it this way.)
If $\tau_{jk}<0$, some signs in the definition of $\beta_{jk}$ must be changed, the result being that path number (3) is now on the lower half unit circle.
All these paths are on the unit circle so $\sigma(\beta_{jk})=\beta_{jk}.$
We define the path $\Gamma_{jk}$ on $\overline{\Sigma}_{t,\x}$ from $1_j$
to $1_k$ as $\Gamma_{jk}=\alpha_{jk}\beta_{jk}\alpha_{kj'}^{-1}$
(or $\Gamma_{jk}=\alpha_{jk}\beta_{jk}$ in case $j$ is the minimum of $E_k\cup R_k$).
It transform as follows under $\sigma$:
\begin{equation}
\label{eq:symmetry-Gammajk}
\sigma(\Gamma_{jk})=\delta_{jk}^{-1}\,\Gamma_{jk}\,\delta_{kj}'.
\end{equation}
\subsection{Formulation of the Monodromy Problem}
Let $\Sigma_{t,\x}$ be the Riemann surface $\overline{\Sigma}_{t,\x}$ minus the poles of $\xi_{t,\x}$, namely
the ends $p_{jk}$ for $(j,k)\in R$, the points $0_j,\infty_j$ for $j\in J$ and the points
$q_{jk},\sigma(q_{jk})$ for $(j,k)\in E^+$.
Fix an arbitrary $j_0\in J$ and take $z_0=1_{j_0}$ as base point.
\begin{proposition}
\label{proposition:monodromy-formulation}
Assume that the Regularity Problem is solved and that
\begin{equation}
\label{eq:gamma-monodromy}
\forall (j,k)\in E\cup R,\quad
\left\{\begin{array}{l}
\cal{P}(\xi_{t,\x},\gamma_{jk})\in\Lambda SU(2)\\
\cal{P}(\xi_{t,\x},\gamma_{jk})\mid_{\lambda=1}=I_2\\
\frac{\partial}{\partial\lambda}\cal{P}(\xi_{t,\x},\gamma_{jk})\mid_{\lambda=1}=0
\end{array}\right.
\end{equation}
\begin{equation}
\label{eq:Gamma-monodromy}
\forall (j,k)\in E^+,\quad
\left\{\begin{array}{l}
\cal{P}(\xi_{t,\x},\Gamma_{jk})\in\Lambda SU(2)\\
\cal{P}(\xi_{t,\x},\Gamma_{jk})\mid_{\lambda=1}=\pm I_2\\
\cal{P}(\xi_{t,\x},\Gamma_{jk})^{-1}\frac{\partial}{\partial\lambda}\cal{P}(\xi_{t,\x},\Gamma_{jk})\mid_{\lambda=1}=i(V_k-V_j)
\end{array}\right.
\end{equation}
where $V_j$ for $j\in J$ are arbitrary matrices in $\su(2)$.
Then the Monodromy Problem \eqref{eq:monodromy-problem-principal} is solved.
\end{proposition}
Proof:
for $j\in J$, let $\gamma_{0_j}$ be a closed loop around $0_j$ in the unit disk of $\CC_j$,
with base point $1_j$.
For $(j,k)\in E^+$, let $\gamma_{q_{jk}}$ be a closed loop with base point $1_j$ defined as follows:
Items (1) and (2) in the definition of $\beta_{jk}$ from $1_j$ to $\z_{jk}'=\varepsilon$ in $\CC_{jk}$,
then a closed loop in the unit disk of $\CC_{jk}$ around $q_{jk}$, and back to $1_j$ by the same path.
Provided the Regularity Problem at $0_j$ and $q_{jk}$ are solved, the gauged potentials
$\xi_{t,\x}\cdot G_j$ and $\xi_{t,\x}\cdot G_{jk}$ have trivial monodromy around $0_j$ and
$q_{jk}$, respectively. Because the gauges have multivaluation $-I_2$ around these points, we have
$$\cal{P}(\xi_{t,\x},\gamma_{0_j})=\cal{P}(\xi_{t,\x},\gamma_{q_{jk}})=-I_2.$$
Any element of $\pi_1(\Sigma_{t,\x},z_0)$ can be written as a product of the following paths or their inverse:
\begin{enumerate}
\item $\gamma_{jk}$ for $(j,k)\in E\cup R$,
\item $\gamma_{0_j}$ for $j\in J$,
\item $\gamma_{q_{jk}}$ for $(j,k)\in E^+$,
\item $\Gamma_{jk}$ for $(j,k)\in E^+$.
\end{enumerate}
Let $c\in\pi_1(\Sigma_{t,\x},z_0)$ and decompose it as
$$c=\prod_{i=0}^{n-1} c_i$$
where each $c_i$ or $c_i^{-1}$ is a path in the above list.
Then
$$\cal{P}(\xi_{t,\x},c)=\prod_{i=0}^{n-1}\cal{P}(\xi_{t,\x},c_i)$$
so we immediately see that the first two items of the Monodromy Problem \eqref{eq:monodromy-problem} are solved.
Each path $c_i$ goes from a point $1_{j_i}$ to a point $1_{j_{i+1}}$,
with $j_{i+1}=j_i$ for paths of type (1), (2) and (3)
and $j_n=j_0$.
Then we always have
$$\cal{P}(\xi_{t,\x},c_i)^{-1}\frac{\partial}{\partial\lambda}\cal{P}(\xi_{t,\x},c_i)\mid_{\lambda=1}
=i(V_{j_{i+1}}-V_{j_i}).$$
Indeed, boths sides are zero for paths of type (1), (2), (3), and for paths of type (4) this follows from
Equation \eqref{eq:Gamma-monodromy}.
Consequently (using that $\pm I_2$ commutes with everything)
\begin{equation}
\label{eq:derivee-Gamma-monodromie}
\cal{P}(\xi_{t,\x},c)^{-1}\frac{\partial}{\partial\lambda}\cal{P}(\xi_{t,\x},c)\mid_{\lambda=1}
=\sum_{i=0}^{n-1}
\cal{P}(\xi_{t,\x},c_i)^{-1}\frac{\partial}{\partial\lambda}\cal{P}(\xi_{t,\x},c_i)\mid_{\lambda=1}
=i(V_{j_n}-V_{j_0})=0.
\end{equation}
\cqfd

We shall take the following choice for the matrices $V_j$:
\begin{equation}
\label{eq:Vj}
V_j=\frac{-i}{2}\matrix{\Re(v_j)&-i\,\Im(v_j)\\ i\,\Im(v_j)&-\Re(v_j)}
\end{equation}
where $v_j$ denotes the vertices of the given graph $\Gamma$.
Then for $(j,k)\in E^+$, we have
$v_k-v_j=\ell_{jk}u_{jk}$
so
\begin{equation}
\label{eq:Vk-Vj}
V_k-V_j=\frac{-i\, \ell_{jk}}{2}\matrix{\Re(u_{jk})&-i\,\Im(u_{jk})\\i\,\Im(u_jk)&-Re(u_{jk})}
=-\frac{\ell_{jk}}{2}N^S(u_{jk}).
\end{equation}
\begin{remark}
\label{remark:formulation}
\begin{enumerate}
\item There is geometry behind our choice for $V_j$: we are in fact requiring that the image of
$1_j$ by the immersion is $v_j$ for all $j\in J$, up to a rigid motion: see Point (2) of Proposition \ref{proposition:spherical}.
\item If the Regularity Problem at $0_j$ and $\infty_j$ is solved, then Equations \eqref{eq:gamma-monodromy} for $k\in E_j\cup R_j$ are not independent, as the fundamental group of the $n$-punctured sphere
has $n-1$ generators.
We will still solve Problems \eqref{eq:gamma-monodromy} for all $k\in E_j\cup R_j$ and infer in Point (3) of Proposition \ref{proposition:spherical}
that the Regularity Problem at $0_j$ and $\infty_j$ is solved. A similar remark holds for the Regularity Problem at
$q_{jk}$.
\end{enumerate}
\end{remark}
\subsection{The renormalized $\gamma$-Monodromy}
\label{section:gamma-monodromy}
In this section, we address the Monodromy Problem \eqref{eq:gamma-monodromy} for the curves $\gamma_{jk}$, $(j,k)\in E\cup R$.
To compensate for the lack of symmetry of $\gamma_{jk}$ (see Equation \eqref{eq:symmetry-gammajk}), we conjugate $\cal{P}(\xi_{t,\x},\gamma_{jk})$
by $\cal{P}(\xi_{t,\x},\delta_{jk}')^{1/2}$ and define
$$\wt{M}_{jk}(t,\x)=\cal{P}(\xi_{t,\x},\delta_{jk}')^{1/2}
\cal{P}(\xi_{t,\x},\gamma_{jk})
\cal{P}(\xi_{t,\x},\delta_{jk}')^{-1/2}.$$
Note that the square root is well-defined for $t$ small enough because at $t=0$, $\xi_{0,\x}$ is
holomorphic at $p_{jk}$ for all $(j,k)\in E\cup R$ so
$\cal{P}(\xi_{0,\x},\delta_{jk}')=I_2$.
As in \cite{nnoids}, we define for $t\neq 0$:
$$\wh{M}_{jk}(t,\x)=t^{-1}\log \wt{M}_{jk}(t,x).$$
\begin{proposition}
\label{proposition:gamma-monodromy}
For $(j,k)\in E\cup R$:
\begin{enumerate}
\item The renormalized monodromy $\wh{M}_{jk}(t,\x)$ extends at $t=0$ to an analytic map of $(t,\x)$ in
a neighborhood of $(0,\cv{\x})$ with value in $\Lambda sl(2,\C)$.
\item $\wh{M}_{jk}(t,\x)$ has the symmetry
\begin{equation}
\label{eq:symmetryMjk}
\overline{\wh{M}_{jk}}=-D\wh{M}_{jk}D^{-1}.
\end{equation}
\item Problem \eqref{eq:gamma-monodromy} is equivalent to the following problem for $(j,k)\in E\cup R$:
\begin{equation}
\label{eq:gamma-monodromy2}\quad\left\{\begin{array}{ll}
\wh{M}_{jk}(t,\x)\in\Lambda\su(2)\quad&(i)\\
\wh{M}_{jk}(t,\x)\mid_{\lambda=1}=0&(ii)\\
\frac{\partial}{\partial\lambda}\wh{M}_{jk}(t,\x)\mid_{\lambda=1}=0&(iii)
\end{array}\right.
\end{equation}
\item At $t=0$, we have
\begin{equation}
\label{eq:wtMjk0}
\wh{M}_{jk}(0,\x)=2\pi i\, \Res_{p_{jk}}\big[\PhiS\frac{\partial\xi_{t,\x}}{\partial t}\mid_{t=0}(\PhiS)^{-1}\big].
\end{equation}
\end{enumerate}
\end{proposition}
Proof:
\begin{enumerate}
\item By standard ODE theory, $\wt{M}_{jk}$ is an analytic map of all parameters.
At $t=0$, $\wt{M}_{jk}(0,\x)=I_2$, so $\wh{M}_{jk}$ extends analytically at $t=0$.
\item By Proposition \ref{proposition:symmetryxi} and Equations \eqref{eq:symmetry-deltajk},
\eqref{eq:symmetry-gammajk}, we have
\begin{equation}
\label{eq:symPdeltajk}
\overline{\cal{P}(\xi_{t,\x},\delta_{jk})}=D\cal{P}(\xi_{t,\x},\delta_{jk})^{-1}D^{-1}
\end{equation}
\begin{equation}
\label{eq:symPgammajk}
\overline{\cal{P}(\xi_{t,\x},\gamma_{jk})}=D\cal{P}(\xi_{t,\x},\delta_{jk}')
\cal{P}(\xi_{t,\x},\gamma_{jk})^{-1}\cal{P}(\xi_{t,\x},\delta_{jk}')^{-1}D^{-1}.
\end{equation}
Hence $\wt{M}_{jk}$ has the symmetry
$$\overline{\wt{M}_{jk}(t,\x)}=D \wt{M}_{jk}(t,\x)^{-1}D^{-1}.$$
Point (2) follows by taking the logarithm, remembering that $t\in\R$.
\item Assuming that $\cal{P}(\xi_{t,\x},\delta_{jk}')$ solves the Monodromy Problem
\eqref{eq:monodromy-problem-principal}, the Monodromy Problem for $\cal{P}(\xi_{t,\x},\gamma_{jk})$ is equivalent to
$$\left\{\begin{array}{l}
\wt{M}_{jk}(t,\x)\in\Lambda SU(2)\\
\wt{M}_{jk}(t,\x)\mid_{\lambda=1}=I_2\\
\frac{\partial}{\partial \lambda}\wt{M}_{jk}(t,\x)\mid_{\lambda=1}=0
\end{array}\right.$$
which, taking the logarithm, is equivalent to Problem \eqref{eq:gamma-monodromy2}.
Remembering the definition of $\delta_{jk}'$, Point (3) follows by induction on $k$ for
the order $\prec$ on $E_j\cup R_j$.
\item We have, since $\cal{P}(\xi_{0,\x},\gamma_{jk})=\cal{P}(\xi_{0,\x},\delta_{jk}')=I_2$:
$$\wh{M}_{jk}(0,\x)=\frac{\partial}{\partial t}\wt{M}_{jk}(t,\x)\mid_{t=0}
=\frac{\partial}{\partial t}\cal{P}(\xi_{t,\x},\gamma_{jk})\mid_{t=0}.$$
At $t=0$, we have $A_j=0$ and $B_j=C_j=1/2$ so $\xi_{0,\x}=\xiS$.
By Proposition 8 in \cite{nnoids}, we obtain
$$\frac{\partial}{\partial t}\cal{P}(\xi_{t,\x},\gamma_{jk})\mid_{t=0}
=\int_{\gamma_{jk}}\PhiS\frac{\partial\xi_{t,\x}}{\partial t}\mid_{t=0}
(\PhiS)^{-1}
=2\pi i\,\Res_{p_{jk}}\PhiS\frac{\partial\xi_{t,\x}}{\partial t}\mid_{t=0}
(\PhiS)^{-1}.$$
\cqfd
 \end{enumerate}
\subsection{The Monodromy Problem around nodes}
\label{section:monodromy-nodes}
In this section we fix $(j,k)\in E$ and solve Problem \eqref{eq:gamma-monodromy2}.
Let $U_{jk}=\Phi^S(p_{jk})$.
In view of Equation \eqref{eq:wtMjk0}, it is advantageous to conjugate $\wh{M}_{jk}$ by
the inverse of $U_{jk}$. Since $p_{jk}\in\S^1$, $U_{jk}\in\Lambda SU(2)$ and
$\overline{U_{jk}}=DU_{jk}D^{-1}$ by Equation \eqref{eq:symmetry-Phi}.
So this conjugation does not affect the Monodromy Problem \eqref{eq:gamma-monodromy2} nor the symmetry
\eqref{eq:symmetryMjk}.
We define
$$\wc{M}_{jk}(t,\x)=U_{jk}^{-1}\wh{M}_{jk}(t,\x) U_{jk}$$
\begin{equation}
\label{eq:Fjk}
\cal{F}_{jk}(t,\x)=i\left(\wc{M}_{jk;11}(t,\x)+\wc{M}_{jk;11}(t,\x)^*\right)
\end{equation}
\begin{equation}
\label{eq:Gjk}
\cal{G}_{jk}(t,\x)=\lambda\left(\wc{M}_{jk;12}(t,\x)+\wc{M}_{jk;21}(t,\x)^*\right).
\end{equation}
so that $\wc{M}_{jk}\in\Lambda\su(2)$ is equivalent to $\cal{F}_{jk}=\cal{G}_{jk}=0$.
By symmetry \eqref{eq:symmetryMjk}, $\cal{F}_{jk}(t,\x)$ and
$\cal{G}_{jk}(t,\x)$ are in $\WR$.
By definition, $\cal{F}_{jk}^*=-\cal{F}_{jk}$ so since $\cal{F}_{jk}^0\in\R$, we have $\cal{F}_{jk}^0=0$ and
we do not need to solve $\cal{F}_{jk}^-=0$.
The $\sigma$-symmetry gives us one more piece of information: if $\wc{M}_{jk}\in\Lambda\su(2)$,
then the symmetry \eqref{eq:symmetryMjk} and the definition of the conjugation and star operators give
$$\wc{M}_{jk;11}(\lambda)=-\wc{M}_{jk;11}^*(\lambda)=\overline{\wc{M}_{jk;11}^*}(\lambda)
=\wc{M}_{jk;11}(1/\lambda).$$
This implies
$$\frac{\partial}{\partial\lambda}\wc{M}_{jk;11}\mid_{\lambda=1}=0.$$
We define
\begin{equation}
\label{eq:defE1jk}
\cal{E}_{1,jk}=(\cal{E}_{1,jk,i})_{1\leq i\leq 6}=\left[
\cal{F}_{jk}^+,\cal{G}_{jk}^+,\lambda(\cal{G}_{jk}^{\leq 0})^*,i\wc{M}_{jk;11}\mid_{\lambda=1},
\wc{M}_{jk;21}\mid_{\lambda=1},\frac{\partial}{\partial\lambda}\wc{M}_{jk;21}\mid_{\lambda=1}\right]
\end{equation}
$$\x_{1,jk}=\left(a_{jk}^+,b_{jk}^+,c_{jk}^+,a_{jk}^0,b_{jk}^0,c_{jk}^0\right)$$
so Problem \eqref{eq:gamma-monodromy2} is equivalent to $\cal{E}_{1,jk}(t,\x)=0$.
\begin{proposition}
\label{proposition:monodromy-nodes}
For $(j,k)\in E$:
\begin{enumerate}
\item $\cal{E}_{1,jk}(t,\x)\in(\WR^{>0})^3\times\R^3$.
\item At the central value, $\cal{E}_{1,jk}(0,\cv{\x})=0$.
\item The partial differential of $\cal{E}_{1,jk}$
with respect to $\x_{1,jk}$ at $(0,\cv{\x})$ is an automorphism of
$(\WR^{>0})^3\times\R^3$.
\item The full differential of $\cal{E}_{1,jk}$ with respect to $\x$ at $(0,\cv{\x})$ only involves the variables
 $\x_{1,jk}$, $r_{jk}$, $A_{jk}$ and $C_{jk}$.
\item If $X\in\Ker(d_{\x}\cal{E}_{1,jk}(0,\cv{\x}))$ satisfies
 $dA_{jk}(X)=0$, then $db_{jk}(X)=dc_{jk}(X)=0.$

 \end{enumerate}
 \end{proposition}
Proof: Point (1) follows from symmetry.
By Propositions \ref{proposition:time0} and \ref{proposition:derivee0},
we have in a neighborhood of $p_{jk}$:
$$\frac{\partial}{\partial t}\xi_{t,\x}\mid_{t=0}=
m_{jk}\frac{dz}{z-p_{jk}}+r_{jk}M_{jk}
\frac{1+q_{jk}^2}{1-q_{jk}^2}
\frac{p_{jk}dz}{(z-p_{jk})^2}+\mbox{ holomorphic}.$$
A simple computation gives
$$\frac{\partial}{\partial z}\left(\PhiS(z)M_{jk}\PhiS(z)^{-1}\right)\mid_{z=p_{jk}}
=p_{jk}^{-1}\PhiS(p_{jk})[M_j,M_{jk}]\PhiS(p_{jk})^{-1}.$$
(Here $M_j$ and $M_{jk}$ have their values at $t=0$.)
Equation \eqref{eq:wtMjk0} gives
$$\wc{M}_{jk}(0,\x)=2\pi i\,m_{jk}+2\pi i\,r_{jk}
\frac{1+q_{jk}^2}{1-q_{jk}^2}\left[M_j,M_{jk}\right].$$
Observe that the partial differential of $\wc{M}_{jk}$ with respect to $q_{jk}$ is zero since ${q}_{jk}=0$ at the central value.
Point (4) follows.
Assume from now on that
$q_{jk}=0$ and
 $A_{jk}=0$.
By Equation \eqref{eq:Bjkt0}, $B_{jk}=\frac{1}{4C_{jk}}$.
We obtain
\begin{equation}
\label{eq:monodromy-nodes-proof}
\wc{M}_{jk}(0,\x)=2\pi i\matrix{a_{jk}&\lambda^{-1}ib_{jk}\\ic_{jk}&-a_{jk}}
%+2\pi i\,r_{jk}\frac{(\lambda^2-1)}{4\lambda}\matrix{-1&0\\0&1}
+\pi i\, r_{jk}\matrix{\lambda^{-1}C_{jk}-\frac{\lambda}{4C_{jk}}&0\\0&\frac{\lambda}{4C_{jk}}-\lambda^{-1}C_{jk}}.
\end{equation}
In particular at the central value, this simplifies to
\begin{equation}
\label{eq:wcMjk0}
\wc{M}_{jk}(0,\cv{\x})=2\pi i\,\tau_{jk}\frac{(\lambda-1)^2}{4\lambda}\matrix{1&0\\0&-1}
\in\Lambda\su(2).
\end{equation}
which proves Point (2).
To prove Point (3), assume that $r_{jk}=\tau_{jk}$ and $C_{jk}=1/2$ are fixed.
Differentiating Equation \eqref{eq:monodromy-nodes-proof} at $(t,\x)=(0,\cv{\x})$ we obtain:
%$$d\cal{F}_{jk}=-2\pi\left(da_{jk}-da_{jk}^*+\smallfrac{1}{2}dr_{jk}(\lambda^{-1}-\lambda)\right)$$
%$$d\cal{G}_{jk}=-2\pi(db_{jk}+\lambda dc_{jk}^*)$$
%$$d\cal{E}_{1,jk,1}=-2\pi(da_{jk}^+-\smallfrac{1}{2}\lambda dr_{jk})$$
%$$d\cal{E}_{1,jk,2}=-2\pi(db_{jk}^++\lambda dc_{jk}^0)$$
%$$d\cal{E}_{1,jk,3}=-2\pi(dc_{jk}^++\lambda db_{jk}^0)$$
%$$d\cal{E}_{1,jk,4}-d\cal{E}_{1,jk,1}\mid_{\lambda=1}=-2\pi(da_{jk}^0+\smallfrac{1}{2}dr_{jk})$$
%$$d\cal{E}_{1,jk,5}-d\cal{E}_{1,jk,3}\mid_{\lambda=1}=-2\pi(dc_{jk}^0-db_{jk}^0)$$
%$$d\cal{E}_{1,jk,6}-\frac{\partial}{\partial\lambda}d\cal{E}_{1,jk,3}\mid_{\lambda=1}=
%2\pi db_{jk}^0.$$
$$d\cal{F}_{jk}=-2\pi\left(da_{jk}-da_{jk}^*\right)$$
$$d\cal{G}_{jk}=-2\pi(db_{jk}+\lambda dc_{jk}^*)$$
$$d\cal{E}_{1,jk,1}=-2\pi\,da_{jk}^+$$
$$d\cal{E}_{1,jk,2}=-2\pi(db_{jk}^++\lambda dc_{jk}^0)$$
$$d\cal{E}_{1,jk,3}=-2\pi(dc_{jk}^++\lambda db_{jk}^0)$$
$$d\cal{E}_{1,jk,4}-d\cal{E}_{1,jk,1}\mid_{\lambda=1}=-2\pi\,da_{jk}^0$$
$$d\cal{E}_{1,jk,5}-d\cal{E}_{1,jk,3}\mid_{\lambda=1}=-2\pi(dc_{jk}^0-db_{jk}^0)$$
$$d\cal{E}_{1,jk,6}-\frac{\partial}{\partial\lambda}d\cal{E}_{1,jk,3}\mid_{\lambda=1}=
2\pi db_{jk}^0.$$
Point (3) easily follows from these formulas.
Finally, to prove Point (5), relax the hypothesis $r_{jk}=\tau_{jk}$ and $C_{jk}=1/2$. By Equation \eqref{eq:monodromy-nodes-proof}, the off-diagonal part of $\wc{M}_{jk}$ does not change, so
$d\cal{E}_{1,jk,i}$ for $i\in\{2,3,5,6\}$ do not change.
Since these equations determine $b_{jk}$ and $c_{jk}$, we obtain $db_{jk}(X)=dc_{jk}(X)=0$.
\cqfd

\begin{remark}
We will solve all equations at the same time by one single application of the Implicit Function Theorem in Section \ref{section:implicit}.
\end{remark}
\subsection{The Monodromy Problem around ends}
\label{section:monodromy-ends}
In this section we fix $(j,k)\in R$ and solve Problem \eqref{eq:gamma-monodromy2}. 
We follow closely the resolution of the same problem in \cite{nnoids}.
We cannot take $U_{jk}=\Phi^S(p_{jk})$ because $p_{jk}\not\in\S^1$ so we take
$U_{jk}=\Phi^S(u_{jk})$ and conjugate $\wh{M}_{jk}$ by the inverse of $U_{jk}$.
Observe that if $a_{jk}=b_{jk}=0$ then $\xi_{t,\x}$ is holomorphic at $p_{jk}$ so
$\wh{M}_{jk}=0$. This prompts us to take
$$a_{jk}=(\lambda-1)^2\,\wh{a}_{jk}\quad\mbox{ and }\quad b_{jk}=(\lambda-1)^2\,\wh{b}_{jk}$$
with $\wh{a}_{jk},\wh{b}_{jk}\in\WRp$.
This way, Points (ii) and (iii) of Problem \eqref{eq:gamma-monodromy2} are automatically satisfied.
We define
$$\wc{M}_{jk}(t,\x)=\frac{\lambda}{(\lambda-1)^2}U_{jk}^{-1}\wh{M}_{jk}(t,\x)U_{jk}$$
which extends at $\lambda=1$ to an analytic map of $(t,\x)$
(see details in Section 6.2 of \cite{nnoids}).
Since $(\lambda-1)^2/\lambda$ is unitary on the unit circle, Point (i) of Problem \eqref{eq:gamma-monodromy2} is equivalent to $\wc{M}_{jk}(t,\x)\in\Lambda\su(2)$.
Define $\cal{F}_{jk}$ and $\cal{G}_{jk}$ by Equations \eqref{eq:Fjk} and \eqref{eq:Gjk} and
$$\cal{E}_{2,jk}=(\cal{E}_{2,jk,i})_{1\leq i\leq 4}=\left[\cal{F}_{jk}^+,\cal{G}_{jk}^+,(\cal{G}_{jk}^-)^*,\cal{G}_{jk}^0\right].$$
Problem \eqref{eq:gamma-monodromy2} is equivalent to
$\cal{E}_{2,jk}(t,\x)=0$.
Writing $p_{jk}=e^{i\theta_{jk}}$ with $\theta_{jk}\in\WRp$, we define
$$\x_{2,jk}=(\wh{a}_{jk}^+,\wh{b}_{jk}^+,\theta_{jk}^+,\wh{b}_{jk}^0).$$

\begin{proposition}
\label{proposition:monodromy-ends}
For $(j,k)\in R$:
\begin{enumerate}
\item $\cal{E}_{2,jk}(t,\x)\in(\WR^{>0})^3\times\R$.
\item At the central value, $\cal{E}_{2,jk}(0,\cv{\x})=0$.
\item The differential of $\cal{E}_{2,jk}$ with respect to $\x$ at $(0,\cv{\x})$ only involves the variable $\x_{2,jk}$ and is an automorphism
of $(\WR^{>0})^3\times\R$.
\end{enumerate}
 \end{proposition}
Proof: Point (1) follows from symmetry. Equation \eqref{eq:wtMjk0} gives
$$\wc{M}_{jk}(0,\x)=2\pi i\lambda\, U_{jk}^{-1}\Res_{p_{jk}}\left[
\PhiS(z)\matrix{0&0\\1&0}\PhiS(z)^{-1}
\left(\frac{\wh{a}_{jk}p_{jk}}{(z-p_{jk})^2}+\frac{i\wh{b}_{jk}}{z-p_{jk}}\right)\right]
U_{jk}.$$
A simple computation gives
$$\frac{\partial}{\partial z}\PhiS(z)\matrix{0&0\\1&0}\PhiS(z)^{-1}\mid_{z=p_{jk}}
=\frac{\lambda^{-1}}{2p_{jk}}\PhiS(p_{jk})\matrix{1&0\\0&-1}\PhiS(p_{jk})^{-1}.$$
This gives
\begin{equation}
\label{eq:whMjk0}
\wc{M}_{jk}(0,\x)=2\pi i\,U_{jk}^{-1}\PhiS(p_{jk})\matrix{\wh{a}_{jk}/2&0\\\lambda i \wh{b}_{jk}&-\wh{a}_{jk}/2}\PhiS(p_{jk})^{-1}U_{jk}.
\end{equation}
At the central value, this simplifies to
$$\wc{M}_{jk}(0,\cv{\x})=2\pi i\,\matrix{\tau_{jk}/4&0\\0&-\tau_{jk}/4}\in\Lambda\su(2)$$
which proves Point (2).
Using Equation \eqref{eq:whMjk0}, we obtain at the central value
$$\frac{\partial}{\partial p_{jk}}\wc{M}_{jk}
=2\pi i\,u_{jk}^{-1}\left[\matrix{0&\lambda^{-1}/2\\\lambda/2&0},\matrix{\tau_{jk}/4&0\\0&-\tau_{jk}/4}\right]
=\frac{2\pi i\,\tau_{jk}}{4 u_{jk}}\matrix{0&-\lambda^{-1}\\\lambda &0}.$$
Hence by the chain rule, since $dp_{jk}/d\theta_{jk}=iu_{jk}$ at $\x=\cv{\x}$:
$$d_{\x}\wc{M}_{jk}(0,\cv{\x})=2\pi i \matrix{d\wh{a}_{jk}/2&-\lambda^{-1}i\tau_{jk}\,d\theta_{jk}/4\\\lambda i \,d\wh{b}_{jk}+\lambda i\tau_{jk}\,d\theta_{jk}/4&-d\wh{a}_{jk}/2}$$
which gives
$$d\cal{E}_{2,jk,1}=-\pi\,d\wh{a}_{jk}^+$$
$$d\cal{E}_{2,jk,2}=\frac{\pi}{2}\tau_{jk}\,d\theta_{jk}^+$$
$$d\cal{E}_{2,jk,3}=-2\pi d\wh{b}_{jk}^+-\frac{\pi}{2}\tau_{jk}\,d\theta_{jk}^+$$
$$d\cal{E}_{2,jk,4}=-2\pi\,d\wh{b}_{jk}^0.$$
Point (3) easily follows.\cqfd
\subsection{The $\Gamma$-Monodromy Problem}
\label{section:Gamma-monodromy}
In this section, we fix $(j,k)\in E^+$ and we solve Problem \eqref{eq:Gamma-monodromy}.
To compensate for the lack of symmetry of $\Gamma_{jk}$ (see Equation \eqref{eq:symmetry-Gammajk}), we multiply $\cal{P}(\xi_{t,\x},\Gamma_{jk})$
by suitable (different) factors on the left and right. Then we conjugate by $\PhiS(u_{jk})\in\Lambda SU(2)$ to simplify computations.
We define for $t>0$:
$$P_{jk}(t,\x)=\PhiS(u_{jk})^{-1}\cal{P}(\xi_{t,\x},\delta_{jk})^{-1/2}
\cal{P}(\xi_{t,\x},\Gamma_{jk})\cal{P}(\xi_{t,\x},\delta_{kj}')^{1/2}
\PhiS(u_{jk}).$$

 \begin{definition}
\label{def-tlogt}
Let $f(t)$ be a function of the real variable $t\geq 0$. We say that $f$ is a smooth function of $t$ and $t\log t$
if there exists a smooth function of two variables $g(t,s)$ defined in a neighborhood of
$(0,0)$ in $\R^2$ such that $f(t)=g(t,t\log t)$ for
$t>0$ and $f(0)=g(0,0)$.
\end{definition}
\begin{remark}
The function $t\log t$ extends continuously at $0$ but the extension is not differentiable at $0$ and is only of
H\"older class $C^{0,\alpha}$ for all $\alpha\in (0,1)$.
Therefore, a smooth function of $t$ and $t\log t$ is only of class $C^{0,\alpha}$ and is not differentiable at $t=0$.
\end{remark}

\begin{proposition}
\label{proposition:Gamma-monodromy}
\begin{enumerate}
\item $P_{jk}(t,\x)$ has the symmetry
$$\overline{P_{jk}}=DP_{jk}D^{-1}.$$
\item $P_{jk}(t,\x)$ extends at $t=0$ to a smooth function of $t$, $t\log t$ and $\x$. Moreover, we have
at $t=0$:
$$P_{jk}(0,\x)=\PhiS(u_{jk})^{-1}
\PhiS(p_{jk})\exp\left(M_{jk}\int_{1}^{-1}\omega_{q_{jk}}\right)\PhiS(p_{kj})^{-1}
\PhiS(u_{jk}).$$
\item At the central value
$$P_{jk}(0,\cv{\x})=\pm\matrix{\lambda&0\\0&\lambda^{-1}}\in\Lambda SU(2).$$
\item Provided Problem \eqref{eq:gamma-monodromy} is solved, Problem
\eqref{eq:Gamma-monodromy} is equivalent to
\begin{equation}
\label{eq:Gamma-monodromy2}
\left\{\begin{array}{ll}
P_{jk}(t,\x)\in\Lambda SU(2)\quad&(i)\\
P_{jk}(t,\x)\mid_{\lambda=1}=\pm I_2&(ii)\\
P_{jk}(t,\x)^{-1}\frac{\partial}{\partial\lambda}P_{jk}(t,\x)\mid_{\lambda=1}=\displaystyle\frac{\ell_{jk}}{2}\matrix{1&0\\0&-1}\quad&(iii)
\end{array}\right.
\end{equation}
\end{enumerate}

\end{proposition}
Proof:
\begin{enumerate}
\item Equation \eqref{eq:symmetry-Gammajk} and Proposition \ref{proposition:symmetryxi} give
$$\overline{\cal{P}(\xi_{t,\x},\Gamma_{jk})}=D\cal{P}(\xi_{t,\x},\delta_{jk})^{-1}
\cal{P}(\xi_{t,\x},\Gamma_{jk})\cal{P}(\xi_{t,\x},\delta_{kj}')D^{-1}.$$
Using the symmetry \eqref{eq:symPdeltajk} and $\overline{\Phi^S(u_{jk})}=D\Phi^S(u_{jk})D^{-1}$, we obtain Point (1).
\item The function $\cal{P}(\xi_{t,\x},\alpha_{jk})$ is an analytic function of all parameters
by Theorem \ref{theorem:openingnodes} because the path
$\alpha_{jk}$ stays away from the nodes.
The same holds for the paths number (1), (3), (5)
 in the definition of the path $\beta_{jk}$
and the path $\alpha_{kj'}$.
By Theorem \ref{theorem:neck} in Appendix \ref{appendix:neck} (see also Remark \ref{remark:neck}), the principal solution of $\xi_{t,\x}$ on path number (2) extends at $t=0$ to a smooth function of $t$, $t\log t$ and $\x$, with the following value at $t=0$:
$$\cal{P}(M_j\omega_0,\z_{jk}=\varepsilon',\z_{jk}=0)\cal{P}(M_{jk}\omega_{q_{jk}},\z'_{jk}=0,\z'_{jk}=\varepsilon').$$
In the same way, the principal solution on path number (4)
 extends to a smooth function of $t$, $t\log t$
and $\x$ with the following value at $t=0$:
$$\cal{P}(M_{jk}\omega_{q_{jk}},\z'_{kj}=-\varepsilon',\z'_{kj}=0)\cal{P}(M_k\omega_0,\z_{kj}=0,\z_{kj}=-\varepsilon').$$
Collecting all terms, the function $\cal{P}(\xi_{t,\x},\Gamma_{jk})$ extends to a smooth function of $t$ and $t\log t$ with the following value at $t=0$:
\begin{eqnarray*}
&&\cal{P}(M_j\omega_0,1_j,e^{i\varepsilon}u_{jk})
\cal{P}(M_j\omega_0,e^{i\varepsilon}u_{jk},\z_{jk}=\varepsilon')
\cal{P}(M_j\omega_0,\z_{jk}=\varepsilon',\z_{jk}=0)\\
&&\times\,\cal{P}(M_{jk}\omega_{q_{jk}},\z'_{jk}=0,\z'_{jk}=\varepsilon')
\cal{P}(M_{jk}\omega_{q_{jk}},\z'_{jk}=\varepsilon',\z'_{kj}=-\varepsilon')
\cal{P}(M_{jk}\omega_{q_{jk}},\z'_{kj}=-\varepsilon',\z'_{kj}=0)\\
&&\times\,\cal{P}(M_k\omega_0,\z_{kj}=0,\z_{kj}=-\varepsilon')
\cal{P}(M_k\omega_0,\z_{kj}=-\varepsilon',e^{i\varepsilon}u_{kj'})
\cal{P}(M_k\omega_0,e^{i\varepsilon}u_{kj'},1_k)\\
&=&\cal{P}(M_j\omega_0,1_j,p_{jk})
\cal{P}(M_{jk}\omega_{q_{jk}},1,-1)
\cal{P}(M_k\omega_0,p_{kj},1_k)
%&=&\PhiS(p_{jk})\exp\left(M_{jk}\int_{1}^{-1}\omega_{q_{jk}}\right)\PhiS(p_{kj})^{-1}.
\end{eqnarray*}
In the above computation, $M_j$ and $M_{jk}$ have their value at $t=0$, so
$M_j\omega_0=\xi^S$.
Point (2) follows.
\item At $(t,\x)=(0,\cv{\x})$, we have $M_{jk}\omega_{q_{jk}}=\xiC$ and $p_{jk}=-p_{kj}=u_{jk}$ so
$$P_{jk}(0,\cv{\x})=\PhiC(-1)\PhiS(-1)^{-1}=\pm\matrix{\lambda&0\\0&\lambda^{-1}}.$$
\begin{remark}
Note that $\PhiS$ and $\PhiC$ are both multivalued with multivaluation $\pm I_2$. This is why we put a $\pm$ sign in Point (3). We do not need to resolve this multivaluation.\end{remark}
\item Assuming that Problem \eqref{eq:gamma-monodromy} is solved, Items (i) and (ii) of Problems \eqref{eq:Gamma-monodromy} and \eqref{eq:Gamma-monodromy2} are clearly equivalent. Assuming Item (ii) is true, we have:
$$\frac{\partial}{\partial\lambda}P_{jk}(t,\x)\mid_{\lambda=1}=
\PhiS(u_{jk})^{-1}\frac{\partial}{\partial\lambda}\cal{P}(\xi_{t,\x},\Gamma_{jk})
\PhiS(u_{jk})\mid_{\lambda=1}.$$
On the other hand, by Equations \eqref{eq:normal} and \eqref{eq:Vk-Vj}:
$$\frac{\ell_{jk}}{2}\PhiS(u_{jk})\matrix{1&0\\0&-1}\PhiS(u_{jk})^{-1}\mid_{\lambda=1}=
\frac{-i\,\ell_{jk}}{2}N^S(u_{jk})=i(V_k-V_j)$$
So Items (iii) of Problems \eqref{eq:Gamma-monodromy} and \eqref{eq:Gamma-monodromy2} are equivalent.
\cqfd
\end{enumerate}
We define for $(t,\x)$ in a neighborhoof of $(0,\cv{\x})$:
$$\wt{P}_{jk}(t,\x)=\log(P_{jk}(t,\x)P_{jk}(0,\cv{\x})^{-1}).$$
By Point (1) of Proposition \ref{proposition:Gamma-monodromy}, $\wt{P}_{jk}$ has the symmetry
\begin{equation}
\label{eq:symmetrywtPjk}
\overline{\wt{P}_{jk}}=D\wt{P}_{jk}D^{-1}.
\end{equation}
\begin{proposition}
Problem \eqref{eq:Gamma-monodromy2}
is equivalent to
\begin{equation}
\label{eq:Gamma-monodromy3}
\left\{\begin{array}{ll}
\wt{P}_{jk}(t,\x)\in\Lambda\su(2)\quad&(i)\\
\wt{P}_{jk;12}(t,\x)\mid_{\lambda=1}=0&(ii)\\
\frac{\partial}{\partial\lambda}\wt{P}_{jk}(t,\x)\mid_{\lambda=1}=\displaystyle\frac{\ell_{jk}-2}{2}\matrix{1&0\\0&-1}\quad&(iii)
\end{array}\right.
\end{equation}
\end{proposition}
Proof:
\begin{itemize}
\item Items (i) of Problems \eqref{eq:Gamma-monodromy2} and \eqref{eq:Gamma-monodromy3}
are equivalent by Point (3) of Proposition \ref{proposition:Gamma-monodromy}.
\item Item (ii) of Problem \eqref{eq:Gamma-monodromy2} is equivalent to
$\wt{P}_{jk}\mid_{\lambda=1}=0$.
Assuming that Item (i) of Problem \eqref{eq:Gamma-monodromy3} holds, we have by symmetry
$$\wt{P}_{jk;11}(1)=-\wt{P}_{jk;11}^*(1)=-\overline{\wt{P}_{jk;11}(1)}=-\wt{P}_{jk;11}(1)$$
so $\wt{P}_{jk;11}(1)=0$. Hence Items (ii) of Problems \eqref{eq:Gamma-monodromy2} and \eqref{eq:Gamma-monodromy3} are equivalent.
\item Assuming Item (ii) of Problem \eqref{eq:Gamma-monodromy2} is satisfied, we have
$$P_{jk}(t,\x)=P_{jk}(0,\cv{\x})=\pm I_2$$
so
$$\frac{\partial}{\partial\lambda}\wt{P}_{jk}(t,\x)\mid_{\lambda=1}
=P_{jk}(t,\x)^{-1}\frac{\partial}{\partial\lambda}P_{jk}(t,\x)\mid_{\lambda=1}
-P_{jk}(0,\cv{\x})^{-1}\frac{\partial}{\partial\lambda}P_{jk}(0,\cv{\x})\mid_{\lambda=1}$$
and by Point (3) of Proposition \ref{proposition:Gamma-monodromy},
$$P_{jk}(0,\cv{\x})^{-1}\frac{\partial}{\partial\lambda}P_{jk}(0,\cv{\x})\mid_{\lambda=1}
=\matrix{1&0\\0&-1}.$$
So Items (iii) of Problem \eqref{eq:Gamma-monodromy2} and \eqref{eq:Gamma-monodromy3} are equivalent.
\cqfd
\end{itemize}
We define
$$\cal{F}_{jk}(t,\x)=\wt{P}_{jk;11}(t,\x)+\wt{P}_{jk;11}(t,\x)^*$$
$$\cal{G}_{jk}(t,\x)=i\left(\wt{P}_{jk;12}(t,\x)+\wt{P}_{jk;21}(t,\x)^*\right).$$
By symmetry \eqref{eq:symmetrywtPjk}, $\cal{F}_{jk}(t,\x)$ and
$\cal{G}_{jk}(t,\x)$ are in $\WR$.
We define
$$\cal{E}_{3,jk}=
(\cal{E}_{3,jk,i})\mid_{1\leq i \leq 7}
=
\left[
\cal{F}_{jk}^+,\cal{G}_{jk}^+,(\cal{G}_{jk}^-)^*,\cal{F}_{jk}^0,\cal{G}_{jk}^0,i\,\wt{P}_{jk;12}\mid_{\lambda=1},
i\,\frac{\partial\wt{P}_{jk;12}}{\partial\lambda}\mid_{\lambda=1}\right]$$
\begin{equation}
\label{eq:defLjk}
\cal{L}_{jk}(t,\x)=\frac{\partial\wt{P}_{jk;11}(t,\x)}{\partial\lambda}\mid_{\lambda=1}-\frac{(\ell_{jk}-2)}{2}.
\end{equation}
Problem \eqref{eq:Gamma-monodromy3} is equivalent to $\cal{E}_{3,jk}(t,\x)=0$ and
$\cal{L}_{jk}(t,\x)=0$.
We leave aside the equation $\cal{L}_{jk}(t,\x)=0$ for the moment and will solve it in Section \ref{section:implicit} using the non-degeneracy hypothesis.
Regarding the equation $\cal{E}_{3,jk}=0$, recall that $q_{jk}\in i\WRp$ and $p_{jk}=e^{i\theta_{jk}}$
with $\theta_{jk}\in\R$ and define
$$\x_{3,jk}=\left(A_{jk}^+,C_{jk}^+,\Im(q_{jk}^+),A_{jk}^0,C_{jk}^0,\theta_{jk},\theta_{kj}\right).$$
\begin{proposition}
\label{proposition:monodromy-edges}
For $(j,k)\in E^+$:
\begin{enumerate}
\item $\cal{E}_{3,jk}(t,\x)\in(\WR^{>0})^3\times\R^4$.
\item At the central value, $\cal{E}_{3,jk}(0,\cv{\x})=0$.
\item The partial differential of $\cal{E}_{3,jk}$ at $(0,\cv{\x})$ with respect to $\x_{3,jk}$ is an automorphism
of $(\WR^{>0})^3\times\R^4$.
\item The full differential of $\cal{E}_{3,jk}$ at $(0,\cv{\x})$ only involves the variables $\x_{3,jk}$ and $\Im(q_{jk}^0)$.
\item If $X\in\Ker(d\cal{E}_{3,jk}(0,\cv{\x}))$, then
$dA_{jk}(X)=dC_{jk}(X)=0.$
 \end{enumerate}
 \end{proposition}
 Proof:
 \begin{itemize}
 \item Point (1) comes from symmetry.
 \item Point (2) is clear since $\wt{P}_{jk}(0,\cv{\x})=0$ by definition.
 \item
 We have at $t=0$
 $$\wt{P}_{jk}(0,\x)=\log\left[\PhiS(u_{jk})^{-1}
 \PhiS(p_{jk})
\exp\left(M_{jk}\int_1^{-1}\omega_{q_{jk}}\right)
 \PhiS(p_{kj})^{-1}\PhiS(u_{kj})\PhiC(-1)^{-1}\right].$$
 Point (4) follows by inspection.
 \item  We set $q_{jk}=0$ to compute the partial derivatives with respect to all parameters but $q_{jk}$.
By Equation \eqref{eq:Bjkt0}, at $t=0$ we have $M_{jk}^2=\frac{1}{4}I_2$ so
 $$\exp\left(M_{jk}\int_1^{-1}\omega_0\right)=\exp(\pi i M_{jk})=2 i M_{jk}.$$
 By Equation \eqref{eq:Bjkt0} we have
 $\partial B_{jk}/\partial A_{jk}=0$ and $\partial B_{jk}/\partial C_{jk}=-1$ at $\x=\cv{\x}$.
 This gives by the chain rule
 $$\frac{\partial \wt{P}_{jk}}{\partial A_{jk}}(0,\cv{\x})=
 2i\frac{\partial M_{jk}}{\partial A_{jk}} \PhiC(-1)^{-1}
 =2i\matrix{i&0\\0&-i}\matrix{0&-i\\-i&0}
 =\matrix{0&2i\\-2i&0}$$
 $$\frac{\partial \wt{P}_{jk}}{\partial C_{jk}}(0,\cv{\x})=
2i\frac{\partial M_{jk}}{\partial C_{jk}} \PhiC(-1)^{-1}
=2i\matrix{0&-1\\1&0}\matrix{0&-i\\-i&0}
 =\matrix{-2&0\\0&2}$$
 $$\frac{\partial \wt{P}_{jk}}{\partial \theta_{jk}}(0,\cv{\x})=
\xiS(u_{jk})\frac{\partial p_{jk}}{\partial\theta_{jk}}
 =\matrix{0&\lambda^{-1}i/2\\\lambda\, i/2&0}$$
 $$\frac{\partial \wt{P}_{jk}}{\partial \theta_{kj}}(0,\cv{\x})=
-\PhiC(-1)\xiS(u_{kj})\frac{\partial p_{kj}}{\partial\theta_{kj}}\PhiC(-1)^{-1}
 =\matrix{0&-\lambda\,i/2\\-\lambda^{-1}i/2&0}.$$
\item Next we compute the partial derivative with respect to $q_{jk}$ at $(0,\cv{\x})$:
 $$\frac{\partial}{\partial q_{jk}}\int_1^{-1}\omega_{q_{jk}}=
 \frac{\partial}{\partial q_{jk}}\int_1^{-1}\frac{dz}{z-q_{jk}}-\frac{q_{jk}\,dz}{1+q_{jk}z}
 \int_1^{-1}\frac{dz}{z^2}-dz=4$$
$$\frac{\partial \wt{P}_{jk}}{\partial q_{jk}}(0,\cv{\x})=
4\,\PhiC(-1)\matrix{0&1/2\\1/2&0}\PhiC(-1)^{-1}
 =\matrix{0&2\\2&0}$$
% \begin{remark} With $q_{jk}\in\WRp$, the previous computation gives
% $\frac{\partial}{\partial q_{jk}}\int_1^{-1}\omega_{q_{jk}}=0$, so we were right
% to assume $q_{jk}\in i\WRp$.
% \end{remark}
 \item
Write $q_{jk}=i \nu_{jk}$ with $\nu_{jk}\in\WRp$.
 Remembering that $\theta_{jk},\theta_{kj}\in\R$ we obtain at $(0,\cv{x})$:
 $$d\cal{F}_{jk}=-2dC_{jk}-2dC_{jk}^*$$
  $$d\cal{G}_{jk}=-2(dA_{jk}+dA_{jk}^*)-2(d \nu_{jk}-d \nu_{jk}^*)$$
 $$d\cal{E}_{3,jk,1}=-2dC_{jk}^+$$
 $$d\cal{E}_{3,jk,2}=-2dA_{jk}^+ -2d \nu_{jk}^+$$
$$d\cal{E}_{3,jk,3}=-2dA_{jk}^+ +2d \nu_{jk}^+$$
$$d\cal{E}_{3,jk,4}=-4dC_{jk}^0$$
$$d\cal{E}_{3,jk,5}=-4dA_{jk}^0.$$
 If $X\in\Ker(d\cal{E}_{3,jk}(0,\cv{\x}))$, we obtain from these formulas 
 $dA_{jk}(X)=dC_{jk}(X)=d\nu_{jk}^+(X)=0$, so
 Point (5) is proved.
 Regarding Point (3), the partial differential of $(\cal{E}_{3,jk,i})_{1\leq i\leq 5}$
 with respect to $(A_{jk}^+,C_{jk}^+,\nu_{jk}^+,A_{jk}^0,C_{jk}^0)$ is clearly an automorphism of
 $(\WR^{>0})^3\times\R^2$.
Observe that $d\cal{E}_{3,jk,i}$ for $1\leq i\leq 5$ do not involve the real variables $\theta_{jk}$ and
 $\theta_{kj}$
so $d\cal{E}_{3,jk}$ has block-triangular form and is suffices to compute the differential of the remaining two equations with respect to these variables:
 $$d_{\theta_{jk},\theta_{kj}}\cal{E}_{3,jk,6}=\smallfrac{1}{2}(-d\theta_{jk}+d\theta_{kj})$$
 $$d_{\theta_{jk},\theta_{kj}}\cal{E}_{3,jk,7}=\smallfrac{1}{2}(d\theta_{jk}+d\theta_{kj}).$$
 Points (3) follows.\cqfd
\end{itemize}
\section{Solving all equations with the Implicit Function Theorem}
\label{section:implicit}
There remains a few parameters that we have not used yet and that we can fix, namely:
$r_{jk}=\tau_{jk}$ for $(j,k)\in E^+$, $\wh{a}_{jk}^0=\tau_{jk}/2$ and $\theta_{jk}^0=\arg(u_{jk})$
for $(j,k)\in R$. Remembering that $q_{jk}=i\nu_{jk}$, we define
$$\x_1=(\x_{1,jk})_{(j,k)\in E}\quad
\cal{E}_1=(\cal{E}_{1,jk})_{(j,k)\in E}$$
$$\x_2=(\x_{2,jk})_{(j,k)\in R}\quad
\cal{E}_2=(\cal{E}_{2,jk})_{(j,k)\in R}$$
$$\x_3=(\x_{3,jk})_{(j,k)\in E^+}\quad
\cal{E}_3=(\cal{E}_{3,jk})_{(j,k)\in E^+}$$
$$\x_4=(r_{kj},
\nu_{jk}^0)_{(j,k)\in E^+}\quad
\cal{E}_4=(\cal{R}_{jk})_{(j,k)\in E^+}$$
$$\x=(\x_1,\x_2,\x_3,\x_4)\quad\cal{E}=(\cal{E}_1,\cal{E}_2,\cal{E}_3,\cal{E}_4).$$
Recall that the central value $\cv{\x}$ depends smoothly on the graph $\Gamma$ (which has not yet been assumed to be balanced).
\begin{proposition}
\label{proposition:implicit}
The partial differential of $\cal{E}(t,\x)$ with respect to $\x$ at $(0,\cv{\x})$ is an automorphism.
By the Implicit Function Theorem, for $t\geq 0$ in a neighhorhood of $0$, there exists
$\x(t,\Gamma)$, depending smoothly on $t$, $t\log t$ and the graph $\Gamma$, such that $\cal{E}(t,\x(t,\Gamma))=0$ and $\x(0,\Gamma)=\cv{\x}(\Gamma)$.
\end{proposition}
Proof:
\begin{enumerate}
\item By Propositions \ref{proposition:monodromy-nodes}, \ref{proposition:monodromy-ends}
and \ref{proposition:monodromy-edges}, the partial differential of $(\cal{E}_1,\cal{E}_2,\cal{E}_3)$ with respect to $(\x_1,\x_2,\x_3)$ has upper-triangular $3\times 3$ block-form, with automorphisms on the diagonal, so is an automorphism.
\item Let us prove that $L=d_{\x}\cal{E}(0,\cv{\x})$ is injective. Let $X\in\Ker(L)$.
By Point (5) of Proposition \ref{proposition:monodromy-edges}, $dA_{jk}(X)=dC_{jk}(X)=0$.
By Point (5) of Proposition \ref{proposition:monodromy-nodes}, $db_{jk}(X)=0$.
Differentiating Equation \eqref{eq:Rjk} and remembering that we fixed $r_{jk}$ so $dr_{jk}=0$, we obtain
$$d_{\x}\,\cal{R}_{jk}(0,\cv{\x})=-\frac{1}{2}dr_{kj}+2i\,db_{jk}^0+2i\tau_{jk}d \nu_{jk}^0.$$
Hence
 $d r_{kj}(X)=d \nu_{jk}^0(X)=0$, so
$X_4=0$. Hence, by Point (1), $X=0$.
\item Since $\x_4,\cal{E}_4$ are in spaces of the same finite dimension, Points (1) and (2) imply that $L$ is an automorphism by
elementary linear algebra.
By Point (2) of Proposition \ref{proposition:Gamma-monodromy}, $\cal{E}(t,\x)$ is a smooth function of
$t$, $t\log t$ and $\x$, which (Definition \ref{def-tlogt}) means that there exists a smooth function
$\wt{\cal{E}}(t,s,\x)$ such that $\cal{E}(t,\x)=\wt{\cal{E}}(t,t\log t,\x)$. We apply the Implicit Function Theorem
to $\wt{\cal{E}}$ at $(t,s,\x)=(0,0,\cv{\x}(\Gamma))$ and obtain a smooth function
$\x(t,s,\Gamma)$ such that $\wt{\cal{E}}(t,s,\x(t,s,\Gamma))=0$. Specializing to $s=t\log t$, we obtain
Proposition \ref{proposition:implicit}.\cqfd
\end{enumerate}

We are not done yet ; we still have to solve the equations $\cal{R}_j=0$ and $\cal{L}_{jk}=0$, where
$\cal{R}_j$ is defined by Equation \eqref{eq:defRj} and $\cal{L}_{jk}$ is defined by Equation \eqref{eq:defLjk}.
Define
$$\cal{F}(t,\Gamma)=\left(\big(\cal{R}_j(t,\x(t,\Gamma))\big)_{j\in J},\big(\cal{L}_{jk}(t,\x(t,\Gamma))\big)_{(j,k)\in E^+}\right).$$
By Proposition \ref{proposition:Rj} and since $\wt{P}_{jk}(0,\cv{\x})=0$, we have:
$$\cal{F}(0,\Gamma)=\left(\big(\smallfrac{1}{2}\,\overline{F_j(\Gamma)}\big)_{j\in J},
\big(1-\smallfrac{1}{2}\ell_{jk}(\Gamma)\big)_{(j,k)\in E^+}\right).$$
By the Implicit Function Theorem, we obtain:
\begin{proposition}
\label{proposition:implicit2}
Assume that the central graph $\cv{\Gamma}$ has length-2 edges, is balanced and non-degenerate.
Then for $t\geq 0$ small enough, there exists a deformation $\Gamma(t)$
of $\cv{\Gamma}$, depending smoothly on $t$ and $t\log t$, such that
$\Gamma(0)=\cv{\Gamma}$ and 
$\cal{F}(t,\Gamma(t))=0$.
\end{proposition}
%%% GEOMETRY %%%

\section{Geometry of the immersion}
From now on, we assume that the parameter vector $\x$ is given by Proposition \ref{proposition:implicit}
and $\Gamma(t)$ is given by Proposition \ref{proposition:implicit2}.
We write $\x_t=\x(t,\Gamma(t))$ which is a smooth function of $t$ and $t\log t$.
To ease notation, we write $\overline{\Sigma}_t=\overline{\Sigma}_{t,\x_t}$
and $\xi_t=\xi_{t,\x_t}$. In the same way, we write
$a_{jk,t}$, $p_{jk,t}$, etc... for the value of the parameters $a_{jk}$, $p_{jk}$ at time $t$.
We also write $\tau_{jk,t}$, $u_{jk,t}$, etc... for the quantities associated to
$\Gamma(t)$, and $\tau_{jk}=\tau_{jk,0}$, $u_{jk}=u_{jk,0}$ for the quantities associated to
the given graph $\Gamma(0)$.

We denote $\Sigma_t$ the Riemann surface $\overline{\Sigma}_t$ minus the poles of $\xi_t$.
Let $p:\wt{\Sigma}_t\to\Sigma_t$ be a universal cover. Recall that we have fixed an arbitrary $j_0\in J$ and taken
$z_0=1_{j_0}$ as base point. Choose an arbitrary $\wt{z}_0$ in the fiber $p^{-1}(z_0)$.
Let $\Phi_t$ be the solution of $d\Phi_t=\Phi_t\xi_t$ on
$\wt{\Sigma}_t$ with initial condition $\Phi_t(\wt{z}_0)=I_2$,
$f_t=\Sym(\Uni(\Phi_t))$ the immersion given by the DPW method and
$\wt{f}_t=\Psi\circ f_t$ where $\Psi$ is the rigid motion given by Equation
\eqref{eq:Psi}.
Recall that $\overline{\Sigma}_t$ does not depend on $\lambda$, but $\Sigma_t$ does, which is a problem as the DPW method requires a fixed Riemann surface. We address this issue in Section \ref{section:delaunay} using the results from \cite{nnoids} where the same problem already occured.
At this point, all we know for sure is that $f_t$ is a well defined immersion on $\Sigma_t$
minus $\varepsilon$-neighborhoods of $0_{jk}$, $\infty_{jk}$ for $(j,k)\in E^+$
and $u_{jk}$ for $(j,k)\in R$.

Fix a small $\varepsilon_1$ such that $0<\varepsilon_1\leq\varepsilon/2$ and for $t>0$ small enough,
consider the following fixed compact subdomains of $\overline{\Sigma}_t$:
$$\Omega_{j,\varepsilon_1}=\CC_j\setminus\bigcup_{k\in E_j\cup R_j}D(u_{jk},\varepsilon_1)
\quad\mbox{ for $j\in J$ (spherical parts)}$$
$$\Omega_{jk,\varepsilon_1}=\CC_{jk}\setminus D(\pm 1,\varepsilon_1)\quad\mbox{ for $(j,k)\in E^+$
(catenoidal parts)}.$$
\subsection{Spherical parts}
Without loss of generality, we may assume by translating the graph that $v_{j_0}=0$ so $V_{j_0}=0$.
Recall the definition of the gauge $G_j$
in Section \ref{section:regularity0j} which we now denote $G_{j,t}$ as it depends on $t$.
\begin{proposition}
\label{proposition:spherical}
For $j\in J$ and $t>0$:
\begin{enumerate}
\item The potential $\xi_t$ restricted to $\Omega_{j,\varepsilon_1}\setminus\{0_j,\infty_j\}$
depends $C^1$ on $t$.
\item $\wt{f}_t(1_j)=v_{j,t}+e_1$.
\item $\xi_t\cdot G_{j,t}$ is regular at $0_j$ and $\infty_j$, so $f_t$ extends analytically to $0_j$ and $\infty_j$.
\item As $t\to 0$, $\wt{f}_t-v_{j,t}$ converges on $\Omega_{j,\varepsilon_1}$ to the inverse stereographic projection $\pi^{-1}:\CC\to\S^2$.
More precisely, we have
$$\|\wt{f}_t-v_{j,t}-\pi^{-1}\|_{C^1(\Omega_{j,\varepsilon_1})}\leq c\,t$$
for some uniform constant $c$ (depending on $\varepsilon_1$) and the norm is computed for the
spherical metric on the Riemann sphere.
\end{enumerate}
\end{proposition}
Proof:
\begin{enumerate}
\item
Recall that $\x(t)$ is a smooth function of $t$ and $t\log t$ so is not even differentiable at $t=0$.
However, assuming that Equation \eqref{eq:Cj} holds, we have, for all values of the parameter $\x$,
$\xi_{0,\x}=\xi^S$ in $\Omega_{j,\varepsilon_1}$, so $\xi_{0,\x}$ does not depend on $\x$.
By Proposition \ref{proposition:tlogt} in Appendix \ref{appendix:tlogt}, $\xi_t=\xi_{t,\x(t)}$, restricted to
$\Omega_{j,\varepsilon_1}\setminus\{0_j,\infty_j\}$, extends to a $C^1$ function of $t$
in a neighborhood of $0$.
\item Choose a path $c$ from $z_0$ to $1_j$ on $\Sigma_t$ and let $\wt{c}$ be the lift
of $c$ to $\wt{\Sigma}_t$ such that $\wt{c}(0)=\wt{z_0}$. Let $\wt{1}_j=\wt{c}(1)\in p^{-1}(1_j)$.
Let $\wt{\Omega}_{j,\varepsilon_1}$ be the component of
$p^{-1}(\Omega_{j,\varepsilon_1}\setminus\{0_j,\infty_j\})$
containing $\wt{1}_j$.
Since $\xi_0=\xi^S$ in $\CC_j$ we have
\begin{equation}
\label{eq:Phi0Cj}
\Phi_0=\Phi_0(\wt{1}_j)\Phi^S\quad\mbox{ in $\wt{\Omega}_{j,\varepsilon_1}$}.
\end{equation}
By Equation \eqref{eq:Gamma-monodromy},
%and \eqref{eq:derivee-Gamma-monodromie}
we have
\begin{equation}
\label{eq:Phit1j}
\left\{\begin{array}{l}
\Phi_t(\wt{1}_j)\in\Lambda SU(2)\\
\Phi_t(\wt{1}_j)\mid_{\lambda=1}=\pm I_2\\
\Phi_t(\wt{1}_j)^{-1}\frac{\partial}{\partial\lambda}\Phi_t(\wt{1}_j)\mid_{\lambda=1}=i\,V_{j,t}.
\end{array}\right.
\end{equation}
By the Sym-Bobenko formula \eqref{eq:sym-bobenko} and Equation \eqref{eq:Vj},
$$f_t(1_j)=2V_{j,t}\sim (0,-\Im(v_{j,t}),-\Re(v_{j,t})).$$
$$\wt{f}_t(1_j)=\Psi(f_t(1_j))=(1+\Re(v_{j,t}),\Im(v_{j,t},0).$$
%% REGULARITY
\item
To prove Point (3), we apply Theorem \ref{theorem:regularity} in Appendix \ref{appendix:regularity}
 to the potential $\wh{\xi}_t=\xi_t\cdot G_{j,t}$.
Let $\ell$ be the maximum of $E_j\cup R_j$ for the order $\prec$.
The path $\delta_{j\ell}$ bounds a disk-type domain in $\Omega_{j,\varepsilon_1}$ containing $0_j$ and $\infty_j$ and not containing $-1$.
The potential $\wh{\xi}_t$ satisfies Hypothesis (1) to (3) of Theorem \ref{theorem:regularity} in $\Omega$
 thanks to
Propositions \ref{proposition:regularity0j}, \ref{proposition:Rj} and \ref{proposition:implicit2}.
By Equation \eqref{eq:induction-deltajk}, $\Phi_t$
solves the Monodromy Problem on $\delta_{j\ell}$.
At $t=0$ we have by Equation \eqref{eq:xjyj} $x_j=1$ and $y_j=0$ so
$G_{j,0}=\GS$. Hence
$$\wh{\xi}_0=\xiS\cdot \GS=\matrix{0&\lambda^{-1}\\0&0}\frac{dz}{2(z+1)^2}.$$
Let $\wh{\Phi}_t=\Phi_0(\wt{1}_j)^{-1}\Phi_t G_{j,t}$.
Since $\Phi_0(\wt{1}_j)\in\Lambda SU(2)$, $\wh{\Phi}_t$
solves the Monodromy Problem on $\delta_{j\ell}$ and
$$\wh{\Phi}_0=\PhiS \GS=\matrix{2&\frac{z-1}{2\lambda(z+1)}\\0&1/2}.$$
Theorem \ref{theorem:regularity}
tells us that $\wh{\xi}_t$ is holomorphic at $0_j$.
Finally $\wh{\beta}_0(0_j)=dz/2$ so $\wh{\beta}_t^0(0_j)\neq 0$ for $t$ small enough
so $f_t$ is regular at $0_j$. Regularity at $\infty_j$ follows by $\sigma$-symmetry.
\item
Let $\wc{\Phi}_t$ be the solution of $d\wc{\Phi}_t=\wc{\Phi}_t\xi_t$ with initial condition
$\wc{\Phi}(\wt{1}_t)=I_2$ and $\wc{f}_t=\Sym(\Uni(\wc{\Phi}))$.
By Point (1) and standard ODE theory, $\wc{\Phi}_t$ is a 
$C^1$ function of $t$ in a neighborhood of $0$ and $z\in\Omega_{j,\varepsilon_1}\setminus\{0_j,\infty_j\}$.
Since Iwasawa decomposition is a diffeomorphism (Theorem \ref{theorem:iwasawa}), $\Uni(\wc{\Phi}_t)$ and $\Pos(\wc{\Phi}_t)$
are $C^1$, so by Equation \eqref{eq:df}, $d\wc{f}_t$ is $C^1$.
Let $K$ be a compact subset of $\overline{\Omega}_{j,\varepsilon_1}\setminus\{0_j,\infty_j\}$.
By the mean value inequality,
$$\|\wc{f}_t(z)-\wc{f}_t(1_j)-\wc{f}_0(z)+\wc{f}_0(1_j)\|\leq C(K)t\quad\mbox{ for $z\in K$}.$$
Since $\Phi_t(\wt{1}_j)\mid_{\lambda=1}=I_2$,
$f_t$ and $\wc{f}_t$ differ by a translation.
Also $\xi_0=\xi^S$ so $\wc{f}_0=f^S$. Hence
$$\| f_t(z)-f_t(1_j)-f^S(z)+f^S(1)\|\leq C(K)t.$$
By Point (2) and Equation \eqref{eq:PsifS} we obtain
$$\|\wt{f}_t(z)-v_{j,t}-\pi^{-1}(z)\|\leq C(K)t\quad\mbox{ for $z\in K$}.$$
This estimate is extended to neighborhoods of $0_j$ and $\infty_j$
using the gauge $G_{j,t}$.
\cqfd
\end{enumerate}
\subsection{Delaunay ends}
\label{section:delaunay}
For $p\in\C$, $D^*(p,r)$ denotes the punctured disk $0<|z-p|<r$.
\begin{proposition}
\label{proposition:delaunay}
There exists $\varepsilon_2>0$ such that for $(j,k)\in R$ and $t>0$ small enough:
\begin{enumerate}
\item $f_t$ extends analytically to $D^*(p_{jk,t}^0,\varepsilon_2)$.
\item $f_t$ has a Delaunay end of weight $\simeq 2\pi t \tau_{jk}$ at $p_{jk,t}^0$.
\item The axis of the Delaunay end of $\wt{f}_t$ at $p_{jk,t}^0$ converges to the half-line $v_j +\R^+u_{jk}$ as $t\to 0$.
\item If $\tau_{jk}>0$, then $f_t(D^*(p_{jk,t}^0,\varepsilon_2))$ is embedded.
\end{enumerate}
\end{proposition}
Proof: These facts are proved in \cite{nnoids} in a similar situation, using general results about
Delaunay ends from \cite{kilian-rossman-schmitt} and \cite{raujouan}.
The potential in \cite{nnoids} has the form
$$\matrix{0&\lambda^{-1}dz\\ t(\lambda-1)^2\omega_t&0}$$
where $\omega_t$ has double poles. We gauge our potential to a similar form so we can apply the results of \cite{nnoids}. Fix $(j,k)\in R$. Recall that $\alpha_t$, $\beta_t$ are holomorphic at $p_{jk,t}$ and
$\gamma_t$ has a double pole with principal part
$$\gamma_t=t(\lambda-1)^2\left(\frac{\wh{a}_{jk,t}p_{jk,t}dz}{(z-p_{jk,t})^2}+\frac{i\,\wh{b}_{jk,t}dz}{z-p_{jk,t}}+O(1)dz\right)$$
where $O(1)$ means a holomorphic function in a neighborhood of $p_{jk,t}$.
Define $\kappa_t\in\Wp$ by $\kappa_t=p_{jk,t}\beta_t(p_{jk,t})/dz$.
At $t=0$, we have $\beta_0=\frac{dz}{2z}$ so $\kappa_0=1/2$.
Consider the gauge
$$G_t=\matrix{\sqrt{\frac{\kappa_t}{z}}&0\\\frac{\lambda}{2\sqrt{\kappa_t z}}&\sqrt{\frac{z}{\kappa_t}}}.$$
A computation gives
$$\wh{\xi}_t
:=\xi_t\cdot G_t=\matrix{\alpha_t+\frac{\beta_t}{2\kappa_t}-\frac{dz}{2z}&
\frac{z\beta_t}{\lambda\kappa_t}\\-\frac{\lambda\alpha_t}{z}-\frac{\lambda\beta_t}{4\kappa_t z}+\frac{\kappa_t\gamma_t}{z}&-\alpha_t-\frac{\beta_t}{2\kappa_t}+\frac{dz}{2z}}.$$
Thanks to our choice of $\kappa_t$ and given the principal part of $\gamma_t$, $\wh{\xi}_t$ has the form
\begin{equation}
\label{eq:delaunay-whxi}
\wh{\xi}_t=\matrix{0&\lambda^{-1}dz\\ t(\lambda-1)^2\omega_t&0}+\matrix{O(1)& O(z-p_{jk,t})\\O(1)&O(1)}
\end{equation}
with
$$\omega_t=\kappa_t\left(\frac{\wh{a}_{jk,t}dz}{(z-p_{jk,t})^2}+\frac{(i\wh{b}_{jk,t}-\wh{a}_{jk,t})dz}{p_{jk,t}(z-p_{jk,t})}\right).$$
The gauged potential $\wh{\xi}_t$ now has the same form as in \cite{nnoids} up to a holomorphic term which is of no consequence (see Remark \ref{remark:delaunay} below).
By Proposition 4 in \cite{nnoids}, $f_t$ extends analytically to $D^*(p_{jk,t}^0,\varepsilon_2)$,
$\kappa_t\wh{a}_{jk,t}$ is a real constant and $f_t$
has a Delaunay end of weight $8\pi t\kappa_t\wh{a}_{jk,t}$ at $p_{jk,t}^0$.
Since $\kappa_0\wh{a}_{jk,0}=\tau_{jk}/4$, Point (2) follows.
Let $\wh{\Phi}_t=\Phi_0(\wt{1}_j)^{-1}\Phi_t G_t$.
At $t=0$ we have by Equation \eqref{eq:Phi0Cj}
$$\wh{\Phi}_0(z)=\PhiS(z)G_0(z)=
\frac{1}{\sqrt{2}}\matrix{1&\lambda^{-1}(z-1)\\\lambda &z+1}
=H\matrix{1&\lambda^{-1}z\\0&1},\quad H=\frac{1}{\sqrt{2}}\matrix{1&-\lambda^{-1}\\\lambda &1}.$$
Let $\wc{\Phi}_t=H^{-1}\wh{\Phi}_t$ and $\wc{f}_t$ be the corresponding immersion.
Then $\wc{\Phi}_0$ has the same value as $\Phi_0$ in \cite{nnoids}.
By Proposition 5 in \cite{nnoids}, the axis of the Delaunay end of $\wc{f}_t$ at $p_{jk,t}^0$ converges
to the half-line through $(0,0,1)$ spanned by $-u_{jk}$.
(The signs in Proposition 5 are actually opposite, but this is because we have the opposite Sym-Bobenko formula in \cite{nnoids}.)
Applying the isometries represented by $H$, $\Phi_0(\wt{1}_j)$ and the rigid motion $\Psi$,
we obtain Point (3).
Point (4) is proved in Proposition 6 in \cite{nnoids}.
\begin{remark}
\label{remark:delaunay}
The proof of Proposition 4 in \cite{nnoids} uses a gauge of the form
$$G=\matrix{\frac{\sqrt{w}}{k}&0\\\frac{-\lambda}{2k\sqrt{w}}&\frac{k}{\sqrt{w}}}\quad\mbox{ with $w=z-p_{jk,t}$ and $k\in\Wp$}.$$
Then
$$\wh{\xi}_t\cdot G=\matrix{\wh{\alpha}_t-\frac{\wh{\beta}_t}{2w}+\frac{dw}{2w}&
\frac{k^2\wh{\beta}_t}{\lambda w}\\ \frac{\lambda\wh{\alpha}_t}{k^2}-\frac{\lambda\wh{\beta}_t}{4k^2w}+\frac{w\wh{\gamma}_t}{k^2}+\frac{\lambda\,dw}{2k^2w}&-\wh{\alpha}_t+\frac{\wh{\beta}_t}{2w}-\frac{dw}{2w}}.$$
What only matters is the residue of $\wh{\xi}_t\cdot G$ at $w=0$.
So the second term in the right-hand side of
\eqref{eq:delaunay-whxi} can be neglected because its $(1,2)$ entry has a zero at $w=0$ and the other entries are holomorphic.
\end{remark}
\subsection{Catenoidal parts}
Recall from Section \ref{section:regularity-qjk} the definition of the complex coordinate $w_{jk}$ on $\CC_{jk}$ and
the gauge $G_{jk}$, which we write respectively $w_{jk,t}$ and $G_{jk,t}$ as they now depend on $t$.
We denote $\wh{\xi}_{jk,t}=(w_{jk,t}^{-1})^*\xi_t\cdot G_{jk,t}$.
We cannot use $z=1$ as base point in $\CC_{jk}$ so we use instead the point
$i_{jk}$ defined as $z=i$ if $\tau_{jk}>0$ and
$z=-i$ if $\tau_{jk}<0$.
This point lies on the path $\Gamma_{jk}$.
\begin{proposition}
\label{proposition:catenoid}
For $(j,k)\in E^+$ and $t>0$ small enough:
\begin{enumerate}
\item The potential $\wh{\xi}_{jk,t}$ is regular at $w=0$. Consequently,
the immersion $f_t$ extends analytically to a neighborhood of $0_{jk}$ and $\infty_{jk}$.
\item The potential $\xi_t$ is regular on $\Sigma_t$, so $f_t$ is a regular immersion.
\item The blow-up $t^{-1}(\wt{f}_t -\wt{f}_t(i_{jk}))$) converges on 
$\Omega_{jk,\varepsilon_1}$ as $t\to 0$
to a minimal catenoidal immersion from $\C\setminus\{\pm 1\}$ to $\R^3$.
The limit catenoid has waist radius $|\tau_{jk}|$ and its axis, oriented from the end
at $z=1$ to the end at $z=-1$, is a line parallel to $u_{jk}$ and oriented by
$\tau_{jk}u_{jk}$.
The convergence is for the $C^1$ norm.
\end{enumerate}
\end{proposition}
Proof: fix $(j,k)\in E^+$.
\begin{enumerate}
%% REGULARITY
\item We start by computing $\Phi_0$ in $\Omega_{jk,\varepsilon_1}$.
Split the path $\Gamma_{jk}$ as $\Gamma_{jk}=\Gamma_{jk1}\Gamma_{jk2}$ with $\Gamma_{jk1}(1)=\Gamma_{jk2}(0)=i_{jk}$.
Consider the lift of $\Gamma_{jk1}$ to $\wt{\Sigma}_t$ starting at $\wt{1}_j$
and let $\wt{\imath}_{jk}$ be its endpoint.
Consider the lift of $\Gamma_{jk2}$ starting at
$\wt{\imath}_{jk}$ and let $\wt{1}_k$ be its endpoint.
Let $\wt{\Omega}_{jk,\varepsilon_1}$ be the component of $p^{-1}(\Omega_{jk,\varepsilon_1}\cap\Sigma_t)$ which contains
$\wt{\imath}_{jk}$.
By Theorem \ref{theorem:neck},
$\Phi_t(\wt{\imath}_{jk})$ extends to a smooth function of $t$ and $t\log t$.
Moreover, since $\xi_0=\xi^C$ in $\Omega_{jk,\varepsilon_1}$, we have
$\Phi_0=M\PhiC$ for some matrix $M$ which is determined by the fact that $\Phi_0$ is continuous at the nodes
(see Remark \ref{remark:neck}). This gives by Equation \eqref{eq:Phi0Cj}
\begin{equation}
\label{eq:Phi0Cjk}
\Phi_0(z)=\Phi_0(\wt{1}_j)\PhiS(u_{jk})\Phi^C(z)
=\Phi_0(\wt{1}_k)\PhiS(u_{kj})\Phi^C(-1)^{-1}\PhiC(z)
\quad\mbox{ in $\wt{\Omega}_{jk,\varepsilon_1}$}.
\end{equation}
\item The proof of Point (1) is essentially the same as the proof of Point (3) of
Proposition \ref{proposition:spherical}.
We apply the dual version of Theorem \ref{theorem:regularity}, Corollary \ref{corollary:regularity} in Appendix \ref{appendix:regularity}.
Observe that
\begin{equation}
\label{eq:path1}
\Gamma_{jk1}^{-1}\gamma_{jk}\Gamma_{jk1}\in\pi_1(\Sigma_t,i_{jk})
\end{equation}
is homotopic to a loop $\delta_1$ contained in $\Omega_{jk,\varepsilon_1}$
going around $1$ in the clockwise direction,
and
\begin{equation}
\label{eq:path2}
\Gamma_{jk2}\gamma_{kj}\Gamma_{jk2}^{-1}\in\pi_1(\Sigma_t,i_{jk})
\end{equation}
is homotopic to a loop $\delta_2$ contained in $\Omega_{jk,\varepsilon_1}$ going around $-1$ in the clockwise direction.
The product of the loops \eqref{eq:path1} and \eqref{eq:path2} is a reparametrization (changing the base point) of
\begin{equation}
\label{eq:path3}
\gamma_{jk}\Gamma_{jk}\gamma_{kj}\Gamma_{jk}^{-1}.
\end{equation}
The Monodromy Problem for $\Phi_t$ on the loop \eqref{eq:path3} is solved so it is also solved on
$\delta_1\delta_2$.
We now make the change of variable $w=w_{jk,t}(z)$.
The path $w_{jk,t}(\delta_1\delta_2)$ bounds a disk-type domain in $\CC\setminus\{\pm 1\}$ containing $0$ and $\infty$.
The potential $\wh{\xi}_t$ satisfies
Hypothesis (1) to (3) of Corollary \ref{corollary:regularity} thanks to
Propositions \ref{proposition:regularity0jk}, \ref{proposition:Rjk} and \ref{proposition:implicit}.
Let $\wh{\Phi}_t=(w_{jk,t}^{-1})^*\Phi_t G_{jk,t}$.
At $t=0$, we have $w_{jk,0}(z)=z$ and by Equation \eqref{eq:xjkyjk} $x_{jk}=1$, $y_{jk}=0$ so
$$G_{jk,0}=\matrix{\frac{\sqrt{z}}{1+z}&\frac{1-z}{\sqrt{z}}\\0&\frac{1+z}{\sqrt{z}}}$$
$$\wh{\xi}_0=\xiC\cdot G_{jk,0}=\matrix{0&0\\1&0}\frac{dz}{2(z+1)^2}.$$
Using Equation \eqref{eq:Phi0Cjk}, we have$$\wh{\Phi}_0=\Phi_0(\wt{1}_j)\PhiS(u_{jk})\PhiC G_{jk,0}
=\Phi_0(\wt{1}_j)\PhiS(u_{jk})\matrix{1/2&0\\\frac{z-1}{2(z+1)}&2}$$
where $\Phi_0(\wt{1}_j)\PhiS(u_{jk})\in\Lambda SU(2)$.
By Corollary \ref{corollary:regularity}, $\wh{\xi}_t$ is holomorphic at $0$.
\begin{remark} To deal with the fact that $\wh{\xi}_t$ is not $C^1$ with respect to $t$, we write
$t=\exp(-1/s^2)$, so $\wh{\xi}_{t(s)}$ extends to a smooth function of $s$ in a neighborhood of $0$,
and use $s$ as the time parameter when applying Corollary \ref{corollary:regularity}.
\end{remark}
\item 
By Equation \eqref{eq:whbeta0}, since $\wh{\beta}_t$ is holomorphic at $w=0$, $\wt{\beta}_t^0$ has a zero of multiplicity at least one at $w=0$. So $\beta_t^0$ has a zero of multiplicity
at least one at $z=q_{jk,t}^0$ and at $z=\sigma(q_{jk,t}^0)$ by symmetry, for a total of
$2\,\mbox{card}(E^+)$ zeros.
It has simple poles at $0_j$ and $\infty_j$ for $j\in J$. By elementary topology, the genus of $\overline{\Sigma}_t$ is
$g=\mbox{card}(E^+)-\mbox{card}(J)+1$. Hence the number of zeros of $\beta_t^0$, counting multiplicities, is equal
to
$$\#\mbox{poles}+2g-2=2\,\mbox{card}(J)+2g-2=2\,\mbox{card}(E^+).$$
So the zeros at $q_{jk,t}^0$ and $\sigma(q_{jk,t}^0)$ are simple and $\beta_t^0$ has no other zero. This proves Point (2), and yields that $\wh{\beta}_t^0$ does not vanish at $w=0$,
so completes the proof of Point (1).
\item To prove Point (3), we use Theorem 4 in \cite{minoids}.
One technical issue is that this theorem requires a $C^1$ family of potentials $\xi_t$ and we do not have that regularity.
This problem is solved as follows.
Forget for a moment that the parameter $\x$ has been determined as a smooth function of $t$ and $\log t$ and consider the potential $\xi_{t,\x}$, only assuming that the parameter
$B_{jk}$ is given by Equation \eqref{eq:Bjk}.
Consider the gauged potential
$$\wc{\xi}_{t,\x}=\xi_{t,\x}\cdot\wc{G}_{\x}\quad\mbox{ with }\quad
\wc{G}_{\x}=\frac{1}{\sqrt{2C_{jk}}}\matrix{1&2i\,A_{jk}\\0&2C_{jk}}.$$
Then at $t=0$ we have in $\Omega_{jk,\varepsilon_1}$, using Proposition \ref{proposition:time0}
and Equation \eqref{eq:Bjk}:
$$\wc{\xi}_{0,\x}=\eta_{0,\x}\cdot \wc{G}_{\x}
=\matrix{0&2(B_{jk}C_{jk}-A_{jk}^2)\\1/2&0}\frac{dz}{z}
=\matrix{0&1/2\\1/2&}\frac{dz}{z}=\xiC.$$
Since this does not depend on $\x$, Proposition \ref{proposition:tlogt} in Appendix \ref{appendix:tlogt}
ensures that $\wc{\xi}_t=\wc{\xi}_{t,\x(t)}$ extends to a $C^1$ function in a neighborhood of $t=0$.
Moreover
$$\frac{d}{dt}\wc{\xi}_t=\frac{\partial}{\partial t}\wc{\xi}_{t,\x}\mid_{(t,\x)=(0,\x(0))}
=\frac{\partial}{\partial t}\xi_{t,\cv{\x}}\mid_{t=0}.$$
Define in $\wt{\Omega}_{jk,\varepsilon_1}$
$$\wc{\Phi}_t=\wc{H}_{jk,t}\Phi_t\wc{G}_{\x(t)}\quad\mbox{ with }\quad
\wc{H}_{jk,t}=\PhiC(\wt{\imath}_{jk})\Uni(\Phi_t(\wt{\imath}_{jk}))^{-1}.$$
At $t=0$, we have $\wc{G}_{\x(0)}=I_2$ and by Equation \eqref{eq:Phi0Cjk}
\begin{equation}
\label{eq:wcHjk0}
\wc{H}_{jk,0}=\PhiC(\wt{\imath}_{jk})\Phi_0(\wt{\imath}_{jk})^{-1}=
\left(\Phi_0(\wt{1}_j)\PhiS(u_{jk})\right)^{-1}.
\end{equation}
Hence $\wc{\Phi}_0=\PhiC$ in $\wt{\Omega}_{jk,\varepsilon_1}$.
Let $\wc{f}_t=\Sym(\Uni(\wc{\Phi}_t))$.
By Theorem 4 in \cite{minoids}, $t^{-1}\wc{f}_t$ converges to a minimal immersion with Weierstrass data
$$g=-\frac{\wc{\Phi}_{0;11}}{\wc{\Phi}_{0;21}}=\frac{z+1}{1-z}$$
$$\omega=4(\wc{\Phi}_{0;21})^2\frac{\partial}{\partial t}\wc{\beta}_t^0\mid_{t=0}
%=\frac{(z+1)^2}{z}\frac{\partial}{\partial t}\beta_{t,\cv{\x}}\mid_{t=0}
=\frac{(z-1)^2}{z}\frac{\tau_{jk}}{2}\left(\frac{dz}{(z-1)^2}-\frac{dz}{(z+1)^2}\right)=\frac{2\tau_{jk}\,dz}{(z+1)^2}$$
using Proposition \ref{proposition:derivee0}.
With the change of variable $w=(z+1)/(1-z)$ we obtain $g=w$ and $\omega=\tau_{jk}dw/w^2$.
This is the Weierstrass data of a catenoid with neck-size $|\tau_{jk}|$ and
vertical axis (from the end at $w=\infty$ to $w=0$) oriented by $-\tau_{jk}e_3$.
Let $h_{jk,t}$ be the rigid motion represented by
$\wc{H}_{jk,t}$ and $\vec{h}_{jk,t}$ its linear part, where the action is given by \eqref{eq:action}.
We have $\wc{f}_t=h_{jk,t}\circ f_t$.
At $t=0$, we have by Equation \eqref{eq:wcHjk0}
$\wc{H}_{jk,0}\mid_{\lambda=1}=\PhiS(u_{jk})^{-1}\mid_{\lambda=1}$.
So by Equation \eqref{eq:normal}, $\vec{h}_{jk,0}$ maps 
$N^S(u_{jk})$ to $e_3$. This means that $t^{-1}(f_t-f_t(i_{jk}))$ converges to a catenoid
with axis (from the end at $z=1$ to $z=-1$) oriented by $-\tau_{jk}N^S(u_{jk})$.
We have
$$\vec{\Psi}(N^S(u_{jk}))=(-\Re(u_{jk}),-\Im(u_{jk}),0)\sim -u_{jk}$$
so $t^{-1}(\wt{f}_t-\wt{f}_t(i_{jk}))$ converges to a catenoid with axis oriented by
$\tau_{jk}u_{jk}$.
The convergence is on compact subsets of
$\Omega_{jk,\varepsilon_1}\setminus\{0_{jk},\infty_{jk}\}$.
It is extended to neighborhoods of $0_{jk}$ and $\infty_{jk}$ using the
gauge $G_{jk,t}$.
\cqfd
\end{enumerate}
\subsection{Edge-length estimate}
Recall that $\ell_{jk,t}=\|v_{k,t}-v_{j,t}\|$ is the length of the edge $(j,k)$ on $\Gamma_t$.
\begin{proposition}
\label{proposition:length-estimate}
As $t\to 0$, we have for $(j,k)\in E^+$
\begin{equation}
\label{eq:length-estimate}
\ell_{jk,t}=2-2\, \tau_{jk}\,t\log t+O(t).
\end{equation}
\begin{equation}
\label{eq:vjkt}
\wt{f}_t(i_{jk})=\smallfrac{1}{2}(v_{j,t}+v_{k,t})+O(t)
\end{equation}
\end{proposition}
\begin{remark} \begin{enumerate}
\item Equation \eqref{eq:length-estimate} estimates how much the spheres centered at $v_j$ and $v_k$
move away from each other if $\tau_{jk}>0$ (or toward each other if $\tau_{jk}<0$) to fit in a catenoidal neck
of size $\simeq \tau_{jk}t$.
It is in agreement with the half-period of a Delaunay surface of necksize
$\tau_{jk}t$ which is known to have the asymptotic \eqref{eq:length-estimate} as $t\to 0$
(see for example Proposition 7 in \cite{mazzeo-pacard} -- a scaling of $1/2$ must be applied because the mean curvature is the trace of the fundamental form in that paper).
\item Equation \eqref{eq:vjkt} tells us that the waist of the catenoidal neck is centered at the middle of
$v_{j,t}$, $v_{k,t}$, up to an $O(t)$ term.
\end{enumerate}
\end{remark}
Proof: forget for a moment that $\x$ is determined as a function of $t$.
We first compute the term of order $t\log t$ in $\cal{P}(\xi_{t,\x},\Gamma_{jk})$.
Recall that the only terms where a $t\log t$ appears
are those corresponding to path numbers (2) and (4) in the definition of $\beta_{jk}$.
To estimate the term corresponding to path number (2), we use Point (3) of Theorem
\ref{theorem:neck} where $\gamma$ denotes the circle $|\z_{jk}|=\varepsilon$. We have
$$\cal{P}(\xi_{t,\x},\gamma)=\cal{P}(\xi_{t,\x},1_j,\z_{jk}=\varepsilon')^{-1}
\cal{P}(\xi_{t,\x},\gamma_{jk})\cal{P}(\xi_{t,\x},1_j,\z_{jk}=\varepsilon')$$
Using Equation \eqref{eq:wcMjk0}
\begin{eqnarray*}
\lefteqn{\frac{\partial}{\partial t}\cal{P}(\xi_{t,\x},\gamma)\mid_{t=0}
=\PhiS(\z_{jk}=\varepsilon')^{-1}\frac{\partial}{\partial t}\cal{P}(\xi_{t,\x},\gamma_{jk})\mid_{t=0}\PhiS(\z_{jk}=\varepsilon')}\\
&=&\PhiS(\z_{jk}=\varepsilon')^{-1}
\PhiS(p_{jk})\wc{M}_{jk}(0,\x)\PhiS(p_{jk})^{-1}\PhiS(\z_{jk}=\varepsilon')\\
&=&2\pi i\tau_{jk}\frac{(\lambda-1)^2}{4\lambda}\PhiS(\z_{jk}=\varepsilon')^{-1}
\PhiS(u_{jk})\matrix{1&0\\0&-1}\PhiS(u_{jk})^{-1}\PhiS(\z_{jk}=\varepsilon')
+O(\x-\cv{\x}).
\end{eqnarray*}
By Theorem \ref{theorem:neck}, the principal solution of $\xi_{t,\x}$ on path number (2) is equal to
\begin{eqnarray*}
&&\left[I_2+\tau_{jk}t\log t\frac{(\lambda-1)^2}{4\lambda}\PhiS(\z_{jk}=\varepsilon')^{-1}
\PhiS(u_{jk})\matrix{1&0\\0&-1}\PhiS(u_{jk})^{-1}\PhiS(\z_{jk}=\varepsilon')\right]\\
&\times&\cal{P}(\xi_{0,\x},\z_{jk}=\varepsilon',\z_{jk}=t_{jk}/\varepsilon')+O(t)+t\log t\,O(\x-\cv{\x}).
\end{eqnarray*}
By the same argument, the principal solution of $\xi_{t,\x}$ on path number (4) is equal to
\begin{eqnarray*}
&&\left[I_2-\tau_{jk}t\log t\frac{(\lambda-1)^2}{4\lambda}\PhiC(\z'_{kj}=-\varepsilon')^{-1}\PhiC(-1)\matrix{1&0\\0&-1}\PhiC(-1)^{-1}\PhiC(\z'_{kj}=-\varepsilon')\right]\\
&\times&\cal{P}(\xi_{0,\x},\z'_{kj}=-\varepsilon',\z'_{kj}=-t_{kj}/\varepsilon')+O(t)
+t\log t\,O(\x-\cv{\x}).
\end{eqnarray*}
The computation in the proof of Point (2) of Proposition \ref{proposition:Gamma-monodromy} gives
after simplification
$$\cal{P}(\xi_{t,\x},\Gamma_{jk})=
\cal{P}(\xi_{0,\x},\Gamma_{jk})+\tau_{jk}t\log t\frac{(\lambda-1)^2}{4\lambda}
\PhiS(u_{jk})\left[\matrix{1&0\\0&-1},\PhiC(-1)\right]\PhiS(u_{kj})^{-1}
+O(t)+t\log t\,O(\x-\cv{\x}).$$
Recalling the definition of $P_{jk}$ and $\wt{P}_{jk}$ from Section \ref{section:Gamma-monodromy}, we obtain
$$P_{jk}(t,\x)=P_{jk}(0,\x)+\tau_{jk}t\log t \frac{(\lambda-1)^2}{4\lambda}
\left[\matrix{1&0\\0&-1},\PhiC(-1)\right]\PhiS(u_{kj})^{-1}\PhiS(u_{jk})+O(t)+t\log t\,O(\x-\cv{\x}).$$
$$\wt{P}_{jk}(t,\x)=\wt{P}_{jk}(0,\x)+\tau_{jk}t\log t\frac{(\lambda-1)^2}{4\lambda}
\matrix{2&0\\0&-2}+O(t)+t\log t\,O(\x-\cv{\x})+O((\x-\cv{\x})^2).$$
We substitute the value $\x(t)=\x(t,\Gamma)$ given by Proposition \ref{proposition:implicit}.
(At this point, the graph $\Gamma$ is fixed.)
Recalling that $\x(t)=\cv{\x}+O(t\log t)$ and $\wt{P}_{jk}(0,\cv{\x})=0$, we obtain
\begin{equation}
\label{eq:wtPjk-tlogt}\wt{P}_{jk}(t,\x(t))=d_{\x}\wt{P}_{jk}(0,\cv{\x})(\x(t)-\cv{\x})
+\tau_{jk}t\log t\frac{(\lambda-1)^2}{2\lambda}\matrix{1&0\\0&-1}+O(t).
\end{equation}
Write $\delta\x(t)=\x(t)-\cv{\x}$ and extend this notation to all parameters.
Recalling from the proof of Proposition \ref{proposition:monodromy-edges} the formula
for $d_{\x}\wt{P}_{jk}(0,\cv{\x})$ and the definition of $\cal{F}_{jk}$,
we obtain
$$\cal{F}_{jk}^+(t,\x(t))=-2\,\delta C_{jk}^+(t) +\lambda\tau_{jk}t\log t+O(t)=0$$
$$\cal{F}_{jk}^0(t,\x(t))=-4\,\delta C_{jk}^0(t) -2\tau_{jk}t\log t+O(t)=0$$
which gives
\begin{equation}
\label{eq:deltaCjkt}
\delta C_{jk}(t)=\frac{\tau_{jk}}{2}t\log t\, (\lambda-1)+O(t).
\end{equation}
The definition of $\cal{L}_{jk}$ and Equation \eqref{eq:wtPjk-tlogt} give
(recalling that $\tau_{jk}$ and $\ell_{jk}$ depend on the graph $\Gamma$)
$$\cal{L}_{jk}(t,\x(t,\Gamma))=\frac{\partial}{\partial\lambda}(-2\delta C_{jk}(t,\Gamma))\mid_{\lambda=1}-\frac{\ell_{jk}(\Gamma)-2}{2}=-\tau_{jk}(\Gamma)t\log t +O(t)-\frac{\ell_{jk}(\Gamma)-2}{2}.$$
By Proposition \ref{proposition:implicit2}, the graph $\Gamma_t$ satisfies
$\cal{L}_{jk}(t,\x(t,\Gamma_t))=0$, and this gives Point (1) of Proposition \ref{proposition:length-estimate}.
\medskip

Since the $t\log t$ factor in Equation \eqref{eq:wtPjk-tlogt} is diagonal, the resolution of the remaining equations of the system $\cal{E}_{3,jk}(t,\x(t))=0$, which only involve the off-diagonal part of $\wt{P}_{jk}(t,\x(t))$, gives
$$\delta A_{jk}(t)=O(t),\quad
\delta q_{jk}^+(t)=O(t),\quad
\delta \theta_{jk}(t)=O(t)\quad\mbox{ and }\quad
\delta\theta_{kj}(t)=O(t).$$
By Point (5) of Proposition \ref{proposition:monodromy-nodes}, we obtain
$\delta b_{jk}(t)=O(t)$. Finally, the resolution of $\cal{R}_{jk}(t,\x(t))=0$
gives $\delta r_{kj}(t)=O(t)$ and $\delta q_{jk}^0(t)=O(t)$ so $\delta q_{jk}(t)=O(t)$.
\medskip

Recall that $\Gamma_{jk1}$ denotes the first half of the path $\Gamma_{jk}$, from $1_j$
to $i_{jk}$. By a computation similar to the above we have
$$\cal{P}(\xi_{t,\x},\Gamma_{jk1})=\cal{P}(\xi_{0,\x},\Gamma_{jk1})+
\tau_{jk}t\log t\frac{(\lambda-1)^2}{4\lambda}\PhiS(u_{jk})\matrix{1&0\\0&-1}\PhiC(i_{jk})+O(t)
+t\log t \,O(\x-\cv{\x})$$
$$\cal{P}(\xi_{0,\x},\Gamma_{jk1})=\PhiS(p_{jk})\exp\left(M_{jk}\int_1^{i_{jk}}\omega_{q_{jk}}\right)$$
$$\frac{\partial}{\partial C_{jk}}\cal{P}(\xi_{0,\x},\Gamma_{jk1})
=\PhiS(u_{jk})\frac{\partial}{\partial C_{jk}}\exp\left(M_{jk}\pi\,i_{jk}/2\right)
=\sqrt{2}\,i_{jk}\PhiS(u_{jk})\matrix{0&-1\\1&0}.$$
We substitute $\x=\x(t)$. Using Equation \eqref{eq:deltaCjkt} and that $\delta q_{jk}$,
$\delta A_{jk}$, $\delta \theta_{jk}$ are $O(t)$, we obtain
$$\cal{P}(\xi_{t,\x(t)},\Gamma_{jk1})=\PhiS(u_{jk})\PhiC(i_{jk})\big(I_2+\tau_{jk}t\log t \,Q_{jk}+O(t)\big)$$
with
\begin{eqnarray*}
Q_{jk}&=&\frac{\sqrt{2}}{2}i_{jk}(\lambda-1)\PhiC(i_{jk})^{-1}\matrix{0&-1\\1&0}+\frac{(\lambda-1)^2}{4\lambda}\PhiC(i_{jk})^{-1}\matrix{1&0\\0&-1}\PhiC(i_{jk})\\
&=&\frac{\lambda-1}{2}\matrix{1&-i_{jk}\\i_{jk}&-1}+\frac{(\lambda-1)^2}{4}\matrix{0&i_{jk}\\-i_{jk}&0}
=\matrix{\frac{1}{2}(\lambda-1)&\frac{i_{jk}}{4}(\lambda^{-1}-\lambda)\\
\frac{i_{jk}}{4}(\lambda-\lambda^{-1})&\frac{1}{2}(1-\lambda)}\end{eqnarray*}
The differential of Iwasawa decomposition at the identity is the projection on the factors of the decomposition of the Lie algebra $\Lambda\sl(2,\C)$ as $\Lambda\su(2)\oplus\Lambda\sl^+_{\R}(2,\C)$.
The matrix $Q_{jk}$ decomposes as
$$Q_{jk}=\frac{i_{jk}}{4}\matrix{0&\lambda^{-1}-\lambda\\ \lambda-\lambda^{-1}&0}
+\frac{1}{2}\matrix{\lambda-1&0\\0&1-\lambda} \in \Lambda\su(2)\oplus\Lambda\sl^+_{\R}(2,\C).$$
Hence
$$\Uni(\cal{P}(\xi_{t,\x(t)},\Gamma_{jk1}))=\PhiS(u_{jk})\PhiC(i_{jk})\left[
I_2+\frac{\tau_{jk}\,i_{jk}}{4}t\log t\matrix{0&\lambda^{-1}-\lambda\\ \lambda-\lambda^{-1}&0}
+O(t)\right].$$
Finally the Sym Bobenko formula \eqref{eq:sym-bobenko} gives
\begin{eqnarray*}
\lefteqn{\Sym(\Uni(\cal{P}(\xi_{t,\x(t)},\Gamma_{jk1})))}\\
&=&
f^S(u_{jk})-2i\frac{\tau_{jk}\,i_{jk}}{4}t\log t\,\PhiS(u_{jk})\PhiC(i_{jk})
\matrix{0&-2\\2&0}\PhiC(i_{jk})^{-1}\PhiS(u_{jk})^{-1}\mid_{\lambda=1}+O(t)\\
&=&f^S(u_{jk})-i\tau_{jk}t\log t\,\PhiS(u_{jk})\matrix{-1&0\\0&1}\PhiS(u_{jk})^{-1}\mid_{\lambda=1}+O(t)\\
&=&f^S(u_{jk})+\tau_{jk}t\log t\,N^S(u_{jk})+O(t)
\end{eqnarray*}
and Point (2) follows.\cqfd
\subsection{Transition annuli}
For $(j,k)\in E$ and $t>0$,
let ${\cal A}_{jk,t}$ be the annulus $|t_{jk}|/\varepsilon<|\z_{jk,t}|<\varepsilon$ which is identified with the annulus $|t_{jk}|/\varepsilon<|\z_{jk,t}'|<\varepsilon$
when opening nodes.
We have for $|z_{jk,t}|\leq 1$
\begin{equation}
\label{eq:estimee-zjkt}
\frac{1}{2}\leq\frac{|z_{jk,t}|}{|z-p_{jk,t}|}\leq\frac{3}{2}.
\end{equation}
So provided $|p_{jk,t}-u_{jk}|\leq \frac{\varepsilon}{6}$, which is true for $t$ small enough,
the outer boundary component of $\cal{A}_{jk,t}$ (namely the circle $|z_{jk,t}|=\varepsilon$)
is included in $\Omega_{j,\varepsilon_1}$.
Likewise, the inner boundary component of $\cal{A}_{jk,t}$ (namely the circle $|z'_{jk,t}|=\varepsilon$) is included in
$\Omega_{jk,\varepsilon_1}$.
\begin{proposition}
\label{proposition:transition} For $t>0$ small enough and $(j,k)\in E^+$:
\begin{enumerate}
\item The images of $\cal{A}_{jk,t}$ and $\cal{A}_{kj,t}$ by $\wt{f}_t$ are graphs over annuli in the plane orthogonal to $u_{jk}$.
\item If $\tau_{jk}>0$, the image of the annulus $\cal{A}_{jk,t}\cup\Omega_{jk,\varepsilon_1}\cup\cal{A}_{kj,t}$ is embedded.
\end{enumerate}
\end{proposition}
Proof:
\begin{enumerate}
\item
We may think of the universal covering $\wt{\cal A}_{jk,t}$ of $\cal{A}_{jk,t}$
as the Riemann surface on which $\log \z_{jk,t}$
is well defined.
Let $c>0$ such that for $u\in\S^1$ and $z\in D(u,\frac{1}{2})$
$$\|\PhiS(z)-\PhiS(u)\|\leq c|z-u|\quad\mbox{ and }\quad
\|\PhiC(z)-\PhiC(u)\|\leq c|z-u|.$$
Then for $t$ small enough we have, by Equations \eqref{eq:Phi0Cj}, \eqref{eq:Phi0Cjk}
 and \eqref{eq:estimee-zjkt}
\begin{equation}
\label{eq:boundary-estimate}
\|\Phi_t-\Phi_0(\wt{1}_j)\PhiS(u_{jk})\|\leq 4c\varepsilon\quad\mbox{ on 
$\partial\cal{A}_{jk,t}$}
\end{equation}
We would like to apply the maximum principle to conclude that
the same estimate holds inside $\cal{A}_{jk,t}$.
This is of course not possible because $\Phi_t$ is not well defined on $\cal{A}_{jk,t}$, but
this problem is easily solved as follows.
Define on $\wt{A}_{jk,t}$
$$G_t=\exp\left(-\frac{\log \z_{jk,t}}{2\pi i}\log \cal{M}(\Phi_t,\gamma_{jk})\right).$$
Then $G_t\Phi_t$ descends to a well defined holomorphic function on $\cal{A}_{jk,t}$.
Also, we have $G_t=I_2+O(t)$, so by Equation \eqref{eq:boundary-estimate},
for $t$ small enough
$$\|G_t\Phi_t-\Phi_0(\wt{1}_j)\PhiS(u_{jk})\|\leq 5c\varepsilon\quad\mbox{ on 
$\partial{\cal A}_{jk,t}$}.$$
By the maximum principle
$$\|G_t\Phi_t-\Phi_0(\wt{1}_j)\PhiS(u_{jk})\|\leq 5c\varepsilon\quad\mbox{ in 
${\cal A}_{jk,t}$}.$$
(The maximum principle for Banach valued holomorphic functions states that if $\|f\|$ has an interior maximum then $\|f\|$ is constant, and is an easy consequence of the Gauss mean value formula.)
Hence for $t$ small enough
$$\|\Phi_t-\Phi_0(\wt{1}_j)\PhiS(u_{jk})\|\leq 6c\varepsilon\quad\mbox{ in
${\cal A}_{jk,t}$}.$$
Fix a positive $\alpha<1/4$.
Using that Iwasawa decomposition is differentiable, we have, provided $\varepsilon$ is chosen small enough
(observe that $c$ is a universal constant)
$$\|\Uni(\Phi_t)-\Phi_0(\wt{1}_j)\PhiS(u_{jk})\|\leq \alpha\quad\mbox{ in $\cal{A}_{jk,t}$}.$$
Let $N_t$ be the Gauss map of $f_t$. By Equation \eqref{eq:normal}, we obtain since
$\Phi_0(\wt{1}_j)\mid_{\lambda=1}=I_2$
$$\|N_t-N^S(u_{jk})\|\leq 2\alpha\quad\mbox{ in $\cal{A}_{jk,t}$}$$
Recall that $\vec{\Psi}\circ N^S(u_{jk})=-u_{jk}$ and
let $\pi_{u_{jk}}^{\perp}$ be the projection on the plane orthogonal to $u_{jk}$.
Then
$\pi_{u_{jk}}^{\perp}\circ\wt{f}_t$ is a local diffeomorphism on $\cal{A}_{jk,t}$.
By Proposition \ref{proposition:spherical}, the image by $\pi_{u_{jk}}^{\perp}\circ\wt{f}_t$
of the outer boundary component $|\z_{jk,t}|=\varepsilon$ is close to a circle of center
$\pi_{u_{jk}}^{\perp}(v_j)$ and radius of order $\varepsilon_1$.
By Proposition \ref{proposition:catenoid},
the image by $\pi_{u_{jk}}^{\perp}\circ\wt{f}_t$ of the inner boundary component $|\z_{jk,t}'|=\varepsilon$ is close to a circle of center $\pi_{u_{jk}}^{\perp}(v_j)$ and radius of order
$t$. Hence the projection of the inner boundary component is inside the projection of the outer
boundary component, so $\pi_{u_{jk}}^{\perp}\circ\wt{f}_t$ is a diffeomorphism onto its image
by a standard covering argument. This proves Point (1).
\item
Introduce a coordinate system $(x,y,z)$ with origin $\wt{f}_t(i_{jk})$ and $x$-axis
parallel to the line $(v_{j,t},v_{k,t})$. In the following, left and right refer to the $x$-axis
(so $v_{j,t}$ is on the left of $v_{k,t}$).
Let $S_t$ be the hemisphere $-1\leq x\leq 0$ of the unit sphere centered at $(-1,0,0)$.
Assume $\tau_{jk}>0$.
 By Proposition \ref{proposition:catenoid}, the right boundary
of $\wt{f}_t(\cal{A}_{jk,t})$ is on the left of $S_t$. By Proposition \ref{proposition:spherical},
the left boundary of $\wt{f}_t(\cal{A}_{jk,t})$ is at distance $O(t)$ from the radius 1-sphere centered at
$v_{j,t}$ so is on the left of $S_t$ by Point (2) of Proposition \ref{proposition:length-estimate}.
Moreover, the mean curvature vector on $\wt{f}_t(\cal{A}_{jk,t})$ points to the left (because it does so
on the left boundary).
By the maximum principle, $\wt{f}_t(\cal{A}_{jk,t})$ is on the left of $S_t$ so in particular lies in the half-space $x<0$.
By the same argument, $\wt{f}_t(\cal{A}_{kj,t})$ lies in the half-space
$x>0$. Hence they are disjoint and it is now clear, from Proposition \ref{proposition:catenoid}
and Point (1), that the image of $\cal{A}_{jk,t}\cup \Omega_{jk,\varepsilon_1}\cup\cal{A}_{kj,t}$
is embedded.\cqfd
\end{enumerate}
\subsection{Embeddedness}
Let $M_t$ be the image of $\wt{f}_t$.
\begin{proposition}
\label{proposition:alexandrov}
If all weights $\tau_{jk}$ are positive, then for $t>0$ small enough, $M_t$ is Alexandrov-embedded.
If moreover the graph $\Gamma$ is pre-embedded, then for $t>0$ small enough, $M_t$ is embedded.
\end{proposition}
Proof: we follow closely the proof of Proposition 7 in \cite{minoids}.
Assume that all $\tau_{jk}$ are positive.
By Proposition \ref{proposition:spherical} and taking $\varepsilon_1>0$ small enough, we may find,
for $(j,k)\in E\cup R$, a Jordan curve $\gamma'_{jk,t}$, freely homotopic to $\gamma_{jk}$, whose image is in a plane $\Pi_{jk,t}$ orthogonal to $(v_{j,t},v_{k,t})$, and moreover:
\begin{itemize}
\item If $(j,k)\in E$, $\gamma'_{jk,t}$ lies in $\Omega_{j,\varepsilon_1}\cap\cal{A}_{jk,t}$,
\item If $(j,k)\in R$, $\gamma'_{jk,t}$ lies in $\Omega_{j,\varepsilon_1}\cap D^*(p_{jk,t}^0,\varepsilon_2)$.
\end{itemize}
Let $\Delta_{jk,t}\subset\Pi_{jk,t}$ be the flat disk bounded by $\wt{f}_t(\gamma'_{jk,t})$.
\begin{itemize}
\item For $j\in J$, let $\Omega'_{j,t}\subset\Omega_{j,\varepsilon_1}$ be the domain bounded by the curves $\gamma'_{jk,t}$ for $k\in E_j\cup R_j$. By Proposition \ref{proposition:spherical},
$\wt{f}_t(\Omega'_{j,t})$ is embedded and does not intersect the disks $\Delta_{jk,t}$ for
$k\in E_j\cup R_j$. Hence the union of $\wt{f}_t(\Omega'_{j,t})$ and $\overline{\Delta}_{jk,t}$
for $k\in E_j\cup R_j$ is the image of a continuous injection of the 2-sphere. 
By the Jordan Brouwer Theorem, it is the boundary of a bounded domain $W_{j,t}$.
\item For $(j,k)\in R$, let $D_{jk,t}$ be the disk bounded by $\gamma'_{jk,t}$ and
$D_{jk,t}^*=D_{jk,t}\setminus\{p_{jk,t}^0\}$. By Proposition \ref{proposition:delaunay},
$\wt{f}_t(D_{jk,t}^*)$ is embedded.
By the proof of Claim 3 in \cite{minoids}, its reunion with $\overline{\Delta}_{jk,t}$ bounds a cylindrically
bounded domain $W_{jk,t}$.
\item For $(j,k)\in E^+$, let $\cal{A}'_{jk,t}\subset\cal{A}_{jk,t}\cup\Omega_{jk,\varepsilon_1}\cup\cal{A}_{kj,t}$ be the annulus bounded by $\gamma'_{jk,t}$ and $\gamma'_{kj,t}$.
By Proposition \ref{proposition:transition}, $\wt{f}_t(\cal{A}'_{jk,t})$ is embedded.
By Claim \ref{claim:alexandrov} below and the Jordan Brouwer Theorem,
$\wt{f}_t(\cal{A}'_{jk,t})\cup \overline{\Delta}_{jk,t}\cup \overline{\Delta}_{kj,t}$ is the boundary of
a bounded domain $W_{jk,t}$.
\end{itemize}
Let $W_t$ be the closed manifold with boundary obtained as the disjoint union of all $\overline{W}_{j,t}$ for $j\in J$
and $\overline{W}_{jk,t}$ for $(j,k)\in E\cup R$, identifying $\overline{W}_{j,t}$ and $\overline{W}_{jk,t}$ for $k\in E_j\cup R_j$
along their common boundary $\Delta_{jk,t}$.
Let $F_t:W_t\to\R^3$ be the canonical injection on each $\overline{W}_{j,t}$ and $\overline{W}_{jk,t}$.
Note that $F_t$ is a priori not injective, since the domains may overlap
(its image $F_t(W_t)$ is what is called an immersed domain.)
But $F_t$ is a proper local diffeomorphism whose boundary restriction parametrizes $M_t$.
Moreover, we may compactify $W_t$ by adding one point per domain $\overline{W}_{jk,t}$ for $(j,k)\in R$. This proves that $M_t$ is
Alexandrov embedded.
\medskip

Assume now that $\Gamma$ is pre-embedded.
Then the domains $W_{j,t}$ for $j\in J$ and
$W_{jk,t}$ for $(j,k)\in E\cup R$ are disjoint, and their closures intersect only along the disks
$\Delta_{jk,t}$.
Hence the map $F_t$ is an embedding so $M_t$ is embedded.
\cqfd
\begin{claim}
\label{claim:alexandrov}
We may choose the curves $\gamma'_{jk,t}$ and $\gamma'_{kj,t}$ so that
$\wt{f}_t(\cal{A}'_{jk,t})$ does not intersect the disks $\Delta_{jk,t}$ and $\Delta_{kj,t}$.
\end{claim}
Proof: we continue with the coordinate system $(x,y,z)$ introduced in the proof of Point (2) of Proposition \ref{proposition:transition}.
By Proposition \ref{proposition:spherical}, we may find a Jordan curve
$\gamma''_{jk,t}$ in $\Omega_{j,\varepsilon_1}\cap\cal{A}_{jk,t}$ whose image is at constant distance from the $x$-axis. Let $\cal{A}_{jk,t}''$ be the annulus bounded by $\gamma_{jk,t}''\cup
\gamma_{kj,t}'$ and $A''_t=\wt{f}_t(\cal{A}_{jk,t}'')$.
Consider half a period of a Delaunay surface $D_t$ with axis $Ox$ and necksize $\tau_{jk}t/2$,
bounded on the left by a circle of maximum radius and on the right by a circle of radius $\tau_{jk}t/2$.
Translate the Delaunay surface $D_t$ from the left until a first contact point $p''_t$ with $A''_t$ occurs. By Propositions \ref{proposition:catenoid} and \ref{proposition:transition}, $p''_t$ cannot be on the right boundary of $D_t$ (which is too small) nor on the left boundary of $D_t$
(which is too big). By the maximum principle, $p''_t$ must be on the left boundary of $A''_t$ and
has minimum $x$-coordinate.
Choose the curve $\gamma'_{jk,t}$ so that $\Pi_{jk,t}$ is the plane orthogonal to the $x$-axis and containing
$p''_t$. Then $A''_t$, being on the right of $\Pi_{jk,t}$, does not intersect $\Delta_{jk,t}$.
The annulus bounded by $\gamma'_{jk,t}$ and $\gamma''_{jk,t}$ is inside $\Omega_{j,\varepsilon_1}$
so its image does not intersect $\Delta_{jk,t}$ by Proposition \ref{proposition:spherical}.
Hence $\wt{f}(\cal{A}'_{jk,t})$ does not intersect $\Delta_{jk,t}$, and in the same way, it does not
intersect $\Delta_{kj,t}$.\cqfd
\medskip

This concludes the proof of Theorem \ref{theorem:main}.
\appendix
\section{A regularity result}
\label{appendix:regularity}
In this section we prove a regularity result in the spirit of Theorem 5 in \cite{opening-nodes} or Theorem 6 in \cite{lawson}.
The philosophy of these results is to idenfity which part of the Regularity Problem is solved when the Monodromy Problem around a singularity is solved.
For use in other papers, we consider the Monodromy Problem associated to the general Sym-Bobenko formula in space forms with Sym-points at $\lambda_1,\lambda_2$, with either
\begin{enumerate}
\item $\lambda_1=\lambda_2=1$ ($\R^3$ case)
\item $\lambda_1=e^{i\theta},\lambda_2=e^{-i\theta}$ with $0<\theta<\pi$ ($\S^3$ case)
\item $\lambda_1=e^q$, $\lambda_2=e^{-q}$ with $q>0$ ($\H^3$ case).
\end{enumerate}
The Monodromy Problem in cases (2) and (3) is
\begin{equation}
\label{eq:monodromy-problem-spaceforms}
\left\{\begin{array}{ll}
\cal{M}(\Phi,\gamma)\in\Lambda SU(2)\\
\cal{M}(\Phi,\gamma)\mid_{\lambda_1}=\cal{M}(\Phi,\gamma)\mid_{\lambda_2}=\pm I_2
\end{array}\right.
\end{equation}
\begin{theorem}
\label{theorem:regularity}
Let $\Omega\subset\CC$ be a $\sigma$-symmetric domain containing the points $0$, $\infty$ and not containing $-1$.
Let $\xi_t=\minimatrix{\alpha_t&\lambda^{-1}\beta_t\\\gamma_t&-\alpha_t}$ be a $C^1$ family of
$\sigma$-symmetric DPW potentials with the following properties:
\begin{enumerate}
\item $\alpha_t$, $\beta_t$ are holomorphic in $\Omega$ and $\gamma_t$ has at most a double pole
at $0$ and $\infty$,
\item $\Re(\Res_0(z\gamma_t))=0$,
\item $\Res_0(\gamma_t^0)=0$.
\end{enumerate}
Assume that there exists a continuous family of $\sigma$-symmetric solutions $\Phi_t$ of $d\Phi_t=\Phi_t\xi_t$ in the universal covering of $\Omega\setminus\{0,\infty\}$
and a $\sigma$-symmetric curve $\delta\subset\Omega$ bounding a disk-type domain containing $0$ and $\infty$,
such that the Monodromy Problem \eqref{eq:monodromy-problem} or \eqref{eq:monodromy-problem-spaceforms} for $\cal{M}(\Phi_t,\delta)$ is solved. Further assume that at $t=0$
$$\xi_0=\matrix{0&\lambda^{-1}\\0&0}\frac{k\,dz}{(z+1)^2}$$
with $k\in\R^*$ and $\Phi_0(1)$ is diagonal.
Then for $t$ in a neighborhood of $0$, $\xi_t$ is holomorphic at $0$ and $\infty$.
\end{theorem}
Proof: Let $(F,B)$ be the Iwasawa decomposition of $\Phi_0(1)$
(both factors are diagonal).
Replacing $\Phi_t$ by $F^{-1}\Phi_t$, we may assume that $\Phi_0(1)$ is a diagonal matrix
in $\Lambda^+_{\R}SL(2,\C)$ so
$$\Phi_0(z)=\matrix{\rho&0\\0&\frac{1}{\rho}}\matrix{1&\frac{k(z-1)}{2\lambda(z+1)}\\0&1}$$
with $\rho\in\WRp$.
By Hypothesis (2) and (3), we may write
$$\Res_0(z\gamma_t)=ia_t\quad\mbox{ and } \quad
\Res_0(\gamma_t)=\lambda( b_t+ic_t)$$
with $a_t,b_t,c_t\in\WRp$. For $\x=(a,b,c)\in(\WRp)^3$, define
$$\omega_{\x}=ia(1-z^2)\frac{dz}{z^2}+\lambda\big(b+ic\frac{1-z}{1+z}\big)\frac{dz}{z}.$$
Then
$$\overline{\sigma^*\omega_{\x}}=-\omega_{\x},\quad
\Res_0(z\omega_{\x})=ia\quad\mbox{ and }\quad
\Res_0\omega_x=\lambda(b+ic).$$
Writing $\x_t=(a_t,b_t,c_t)$, we see that $\gamma_t-\omega_{\x_t}$ is holomorphic at $0$ and
$\infty$ by symmetry.
Define
$$\xi_{t,\x}=\matrix{\alpha_t&\lambda^{-1}\beta_t\\ \gamma_t-\omega_{\x_t}&-\alpha_t}+\matrix{0&0\\\ \omega_{\x}&0}.$$
(The first term is holomorphic in $\Omega$).
Let $\Phi_{t,\x}$ be the solution of $d\Phi_{t,\x}=\Phi_{t,\x}\xi_{t,\x}$ with initial condition
$\Phi_{t,\x}(z_0)=\Phi_t(z_0)$, where $z_0$ is an arbitrary base point.
Then $\xi_{t,\x_t}=\xi_t$ and $\Phi_{t,\x_t}=\Phi_t$.
Since $\xi_0$ is holomorphic at $0$, we have $\x_0=0$.
Theorem \ref{theorem:regularity} follows from the following
\begin{lemma}
For $(t,\x)$ in a neighborhood of $(0,0)$, the only solution to the Monodromy Problem \eqref{eq:monodromy-problem} or \eqref{eq:monodromy-problem-spaceforms} for
$\cal{M}(\Phi_{t,\x},\delta)$ is $\x=0$.
\end{lemma}
Proof: let
$$M(t,\x)=H\log\cal{M}(\Phi_{t,\x},\delta)H^{-1}
\quad\mbox{ with }\quad 
H=\matrix{\lambda^{1/2}&0\\0&\lambda^{-1/2}}\in\Lambda SU(2).$$
(The reason to conjugate by $H$ will be clear in a moment.)
By Proposition 8 in \cite{nnoids}, the partial differential of $M$ with respect to $\x$ at $(0,0)$, applied to the vector $\x=(a,b,c)$, is
given by
$$d_{\x}M(0,0)\cdot\x=\int_{\delta}N\omega_{\x}$$
where
$$N=H\Phi_0\matrix{0&0\\1&0}\Phi_0^{-1}H^{-1}
=\matrix{\frac{k(z-1)}{2\lambda(z+1)}&\frac{-k^2\rho^2(z-1)^2}{4\lambda(z+1)^2}\\\frac{1}{\lambda \rho^2}&-\frac{k(z-1)}{2\lambda(z+1)}}.$$
Since $N\omega_x$ has only poles at $0$, $-1$ and $\infty$, we have by the Residue Theorem
$$d_{\x}M(0,0)\cdot\x=2\pi i\big(\Res_0(N\omega_{\x})+\Res_{\infty}(N\omega_{\x})\big)
=-2\pi i\,\Res_{-1}(N\omega_{\x}).$$
Computing the residue at $z=-1$, we obtain
$$d_{\x}M(0,0)=2\pi i\matrix{-k\,db&\lambda^{-1}ik^2\rho^2(2da-\lambda dc/2)\\
2i\,dc/\rho^2&k\,db}.$$
Define
$$\wt{a}=-kb,\quad
\wt{b}=k^2\rho^2(2a-\lambda c/2)\quad\mbox{ and }\quad
\wt{c}=2c/\rho^2.$$
It is clear that $(a,b,c)\mapsto (\wt{a},\wt{b},\wt{c})$ is an automorphism of $(\WRp)^3$.
The point of this change of variables (and the conjugation by $H$) is that we now have
$$d_{\x}M(0,0)=2\pi i\matrix{d\wt{a}&\lambda^{-1}i\,d\wt{b}\\ i\,d\wt{c}&-d\wt{a}}.$$
This is precisely Equation \eqref{eq:monodromy-nodes-proof} with $a_{jk}$, $b_{jk}$, $c_{jk}$ 
replaced by $\wt{a}$, $\wt{b}$, $\wt{c}$ and $r_{jk}=0$.
So in the $\R^3$ case, the proof of Point (3) of Proposition
\ref{proposition:monodromy-nodes} yields that for $t$ in a neighborhood of $0$, the Monodromy Problem
\eqref{eq:monodromy-problem} uniquely determines $(\wt{a},\wt{b},\wt{c})$, hence $\x$, as a function of 
$t$.
Now $\x=0$ is a trivial solution (since $\xi_{t,0}$ is holomorphic in $\Omega$) so $\x=0$ is the unique solution. In the $\S^3$ case (respectively the $\H^3$ case), the Monodromy Problem is equivalent,
using the $\rho$-symmetry, to 
$$\left\{\begin{array}{l}
M\in\Lambda\su(2)\\
\Im(M_{11}\mid_{\lambda=e^{i\theta}})=0\\
M_{21}\mid_{\lambda=e^{i\theta}}=0
\end{array}\right.
\quad\mbox{ respectively }\quad
\left\{\begin{array}{l}
M\in\Lambda\su(2)\\
\Im(M_{11}\mid_{\lambda=e^q})=0\\
\Re(M_{12}\mid_{\lambda=e^q})=0\\
\Re(M_{21}\mid_{\lambda=e^q})=0\end{array}\right.
$$
In the $\S^3$ case, define $\cal{F}$, $\cal{G}$ by Equations \eqref{eq:Fjk} and \eqref{eq:Gjk} with $M$ in place
of $\wc{M}_{jk}$ and
$$\cal{E}=(\cal{E}_i)_{1\leq i\leq 5}=\big(\cal{F}^+,\cal{G}^+,\lambda(\cal{G}^-)^*,\Im(M_{11}\mid_{\lambda=e^{i\theta}}),
M_{21}\mid_{\lambda=e^{i\theta}}\big).$$
Then
$$d\cal{E}_1=-2\pi d\wt{a}^+$$
$$d\cal{E}_2=-2\pi (d\wt{b}^++\lambda d\wt{c}^0)$$
$$d\cal{E}_3=-2\pi(d\wt{c}^+ +\lambda d\wt{b}^0)$$
$$d\cal{E}_4+\Re(d\cal{E}_1\mid_{\lambda=e^{i\theta}})=2\pi d\wt{a}^0$$
$$d\cal{E}_5-d\cal{E}_3\mid_{\lambda=e^{i\theta}}=2\pi(e^{i\theta}d\wt{b}^0-d\wt{c}^0).$$
It easily follows (since $e^{i\theta}\not\in\R$) that $d\cal{E}(0,0)$ is an isomorphism from
$\WR^{>0}\times\R^3$ to $\WR^{>0}\times\R\times\C$. Again, we conclude with the Implicit Function Theorem that the only solution of the Monodromy Problem is $\x=0$.
We omit the proof in the $\H^3$ case which is similar.
\cqfd
\medskip

By duality, we obtain the following result (with the same hypothesis on $\Omega$ and $\delta$):
\begin{corollary}
\label{corollary:regularity}
Let $\xi_t$ be a $C^1$ family of $\sigma$-symmetric DPW potentials on $\Omega$ with the following properties:
\begin{enumerate}
\item $\alpha_t$, $\gamma_t$ are holomorphic in $\Omega$ and $\beta_t$ has at most a double pole
at $0$ and $\infty$,
\item $\Re(\Res_0(z\beta_t))=0$,
\item $\Res_0(\beta_t^0)=0$.
\end{enumerate}
Assume that there exists a continuous family of $\sigma$-symmetric solutions $\Phi_t$ of $d\Phi_t=\Phi_t\xi_t$ such that the Monodromy Problem \eqref{eq:monodromy-problem} or \eqref{eq:monodromy-problem-spaceforms} for $\cal{M}(\Phi_t,\delta)$ is solved. Further assume that at $t=0$
$$\xi_0=\matrix{0&0\\1&0}\frac{k\,dz}{(z+1)^2}$$
and $\Phi_0(1)$ is diagonal.
Then for $t$ in a neighborhood of $0$, $\xi_t$ is holomorphic at $0$ and $\infty$.
\end{corollary}
\section{Principal solution through a neck}
\label{appendix:neck}
Fix some numbers $0<\varepsilon'<\varepsilon$. For $t\in\C$ such that $0<|t|<\varepsilon^2$,
let $\cal{A}_t\subset\C$ be the annulus $|t|/\varepsilon<|z|<\varepsilon$ and
$\psi_t:\cal{A}_t\to\cal{A}_t$ be the involution defined by $\psi_t(z)=t/z$.
We see an element of the universal cover $\wt{\C^*}$ of $\C^*$ as a complex number $t\in\C^*$ with a determination of its argument (which we do not write), so the function $\log t$ is well defined on
$\wt{\C^*}$.
We denote $t\mapsto e^{2\pi i}t$ the Deck transformation of $\wt{\C^*}$ which increases the argument of
$t$ by $2\pi$.
For $t\in\wt{\C^*}$, let $\beta_t$ be the curve from $\varepsilon'$ to $t/\varepsilon'$ 
parametrized for $s\in[0,1]$ by
$$\beta_t(s)=(\varepsilon')^{1-2s}t^s=(\varepsilon')^{1-2s}e^{s\log t}.$$
Our goal is to understand the limit behavior of
$\cal{P}(\xi_t,\beta_t)$ as $t\to 0$, under suitable hypothesis on the potential $\xi_t$.
Let $\gamma$ be the circle parametrized by $\gamma(s)=\varepsilon'e^{2\pi is}$.
\begin{theorem}
\label{theorem:neck}
Let $\xi_t$ be a family of $\Lambda\sl(n,\C)$ valued holomorphic 1-forms on $\cal{A}_t$, depending holomorphically
on $t\in D^*(0,\varepsilon^2)$, and let $\wh{\xi}_t=\psi_t^*\xi_t$. Assume that
$$\lim_{t\to 0}\xi_t=\xi_0\quad\mbox{ and }\quad\lim_{t\to 0}\wh{\xi}_t=\wh{\xi}_0$$
where $\xi_0$ and $\wh{\xi}_0$ are holomorphic in $D(0,\varepsilon)$ and the limit is uniform on
compact subsets of $D^*(0,\varepsilon)$.
Define for $t\in\wt{\C^*}$ small enough
$$\wt{F}(t)=\cal{P}(\xi_t,\gamma)^{-\frac{\log t}{2\pi i}}\cal{P}(\xi_t,\beta_t).$$
Then
\begin{enumerate}
\item The function $\wt{F}$ satisfies $\wt{F}(e^{2\pi i}t)=\wt{F}(t)$ so descends to a well defined holomorphic function $F(t)$ defined in a punctured neighborhood of $0$.
\item The function $F$ extends holomorphically at $t=0$ with
$$F(0)=\cal{P}(\xi_0,\varepsilon',0)\cal{P}(\wh{\xi}_0,0,\varepsilon').$$
\item If $t>0$, the function $\cal{P}(\xi_t,\beta_t)$ extends to a smooth function of 
$t$ and $t\log t$ with value $F(0)$ at $t=0$. Moreover we have as $t\to 0$
$$\cal{P}(\xi_t,\beta_t)=\left(I_2+\frac{t\log t}{2\pi i}\frac{\partial}{\partial t}\cal{P}(\xi_t,\gamma)\mid_{t=0}\right)F(0)+O(t).$$
\end{enumerate}
\end{theorem}
\begin{remark}
\label{remark:neck}
We apply Theorem \ref{theorem:neck} in the proof of Proposition \ref{proposition:Gamma-monodromy}
with $\xi_t=(\z_{jk}^{-1})^*\xi_{t,\x}$ and $t=t_{jk}$.
Then $\z'_{jk}=\psi_t\circ \z_{jk}$ so
$\wh{\xi}_t=((\z'_{jk})^{-1})^*\xi_{t,\x}$.
By Proposition \ref{proposition:time0}, $\xi_0$ and $\wh{\xi}_0$ are both holomorphic
in $D(0,\varepsilon)$, with $\xi_0=(\z_{jk}^{-1})^*(M_j\omega_0)$
and $\wh{\xi}_0=((\z'_{jk})^{-1})^*(M_{jk}\omega_{q_{jk}})$.
Theorem \ref{theorem:neck} says that
$\cal{P}((\z_{jk}^{-1})^*\xi_{t,\x},\varepsilon',t_{jk}/\varepsilon')$ extends at $t=0$ to
a smooth function of $t_{jk}$ and $t_{jk}\log t_{jk}$ with value at $t_{jk}=0$
$$\cal{P}((\z_{jk}^{-1})^*(M_j\omega_0),\varepsilon',0)\,\cal{P}(((\z'_{jk})^{-1})^*(M_{jk}\omega_{q_{jk}}),0,\varepsilon').$$
To justify that the extension is a smooth function of $t,t\log t$ and $\x$ we use Hartog Theorem on 
separate holomorphy to ensure that the function $F$ depends holomorphically on $(t,\x)$.
In other words,
$\cal{P}(\xi_{t,\x},\z_{jk}=\varepsilon',\z'_{jk}=\varepsilon')$ extends to a smooth function of 
$t$, $t\log t$ and $\x$ with value at $t=0$
$$\cal{P}(M_j\omega_0,\z_{jk}=\varepsilon',\z_{jk}=0)\,\cal{P}(M_{jk}\omega_{q_{jk}},\z'_{jk}=0,\z'_{jk}=\varepsilon')
=\cal{P}(\xi_{0,\x},\z_{jk}=\varepsilon',\z'_{jk}=\varepsilon')$$
where in the last expression, it is understood that the principal solution is continuous at the node.
This gives some theoretical ground for the heuristic explained in Section \ref{section:strategy}.
\end{remark}
Proof: first of all, by the change of variables $z'=z/\varepsilon'$ and $t'=t/(\varepsilon')^2$, we may
assume without loss of generality that $\varepsilon'=1$ (so $\varepsilon>1$). The expression of $\beta_t$ simplifies
to $\beta_t(s)=t^s$.
The restriction of $\xi_t$ to the unit circle $\gamma$ extends holomorphically at $t=0$, with value
$\xi_0$. Since $\xi_0$ is holomorphic in $D(0,\varepsilon)$, $\cal{P}(\xi_0,\gamma)=I_2$.
Hence $\log \cal{P}(\xi_t,\gamma)$ and $\wt{F}(t)$ are well defined for $t$ small enough.
Point (1) follows from the fact that the path $\beta_{e^{2\pi i}t}$ is homotopic to $\gamma\beta_t$.
To prove Point (2), we split the path $\beta_t$ into $\beta_t=\alpha_t\wh{\alpha}_t^{-1}$
where
$$\alpha_t(s)=\beta_t(s/2)=t^{s/2}\quad\mbox{ and }\quad
\wh{\alpha}_t(s)=\beta_t(1-s/2)=\psi_t(\alpha_t(s)).$$
Since $F(t)$ is well-defined, we may assume that $|\arg t|\leq \pi$. Provided $|t|\leq e^{-\pi}$, we have
$|\log t|\leq 2|\log |t|\,|$ so
\begin{equation}
\label{eq:neck-alphau}
|\alpha'_t(s)|=\smallfrac{1}{2}|t|^{s/2}|\log t|\leq |t|^{s/2}|\log |t|\,|.
\end{equation}
Integrating the estimate \eqref{eq:neck-alphau}, we see that the length of the spiral $\alpha_t$ is bounded by $2$.
\begin{lemma}
\label{lemma:neck}
There exists a uniform constant $C$ such that for $t$ small enough enough:
\begin{equation}
\label{eq:neck-estimate1}
\int_0^1\left\|\,\left[\xi_t(\alpha_t(s))-\xi_0(\alpha_t(s))\right]\alpha'_t(s)\right\|ds\leq C|t|^{1/2}.\end{equation}
\end{lemma}
Proof: we use the letter $C$ to denote various uniform constants.
Fix some $\varepsilon_2\in (1,\varepsilon)$.
On the circle $C(0,\varepsilon_2)$, $\xi_t$ depends holomorphically on $t$ in a neighborhood of $0$
so 
\begin{equation}
\label{eq:neck-estimate2}
\int_{C(0,\varepsilon_2)}\|\xi_t-\xi_0\|\leq C |t|.
\end{equation}
By the change of variable formula, the convergence of $\wh{\xi}_t$ to $\wh{\xi}_0$ and the holomorphicity
of $\xi_0$ and $\wh{\xi}_0$ in $D(0,\varepsilon)$:
\begin{equation}
\label{eq:neck-estimate3}
\int_{C(0,|t|/\varepsilon_2)}\|\xi_t-\xi_0\|
\leq\int_{C(0,\varepsilon_2)}\|\wh{\xi}_t\|+\int_{C(0,|t|/\varepsilon_2)}\|\xi_0\|
\leq C.
\end{equation}
We expand $\xi_t-\xi_0$ in Laurent series in the annulus $|t|/\varepsilon_2\leq|z|\leq\varepsilon_2$ as
$$\xi_t(z)-\xi_0(z)=\sum_{k\in\Z}A_k(t)z^kdz$$
where the matrices $A_k(t)$ are given by
$$A_k(t)=\frac{1}{2\pi i}\int_{C(0,\varepsilon_2)}\frac{\xi_t(z)-\xi_0(z)}{z^{k+1}}
=\frac{1}{2\pi i}\int_{C(0,|t|/\varepsilon_2)}\frac{\xi_t(z)-\xi_0(z)}{z^{k+1}}.$$
Estimates \eqref{eq:neck-estimate2} and \eqref{eq:neck-estimate3} give us respectively:
\begin{equation}
\label{eq:neck-estimate4}
\|A_k(t)\|\leq C\frac{|t|}{\varepsilon_2^{k+1}}\quad\mbox{ and }\quad
\|A_k(t)\|\leq C\frac{\varepsilon_2^{k+1}}{|t|^{k+1}}.\end{equation}
Then we have the following estimates:
\begin{eqnarray*}
\lefteqn{\int_0^1\left\|[\xi_t(\alpha_t(s))-\xi_0(\alpha_t(s))]\alpha'_t(s)\right\|\,ds
\leq\sum_{k\in\Z}\int_0^1\|A_k(t)\|\,|\alpha_t(s)|^k\,|\alpha'_t(s)|\,ds}\\
&\leq&\sum_{k\in\Z}\int_0^1 \|A_k(t)\| \,|t|^{(k+1)s/2}\,|\log|t|\,|\,ds
\quad\mbox{ using \eqref{eq:neck-alphau}}\\
&=&\|A_{-1}(t)\|\,|\log|t|\,|+\sum_{k\neq -1}\frac{2}{k+1}\|A_k(t)\|\left(
1-|t|^{(k+1)/2}\right)\\
&\leq &\|A_{-1}(t)\|\,|\log|t|\,|+2\sum_{k\geq 0}\|A_k(t)\|+
2\sum_{k\leq -2}\|A_k(t)\|\, |t|^{\frac{k+1}{2}}\\
&\leq& C|t\log|t|\,|+C\sum_{k\geq 0}\frac{|t|}{\varepsilon_2^{k+1}}
+C\sum_{k\leq -2}\left(\frac{\varepsilon_2}{|t|^{1/2}}\right)^{k+1}
\quad\mbox{ using \eqref{eq:neck-estimate4}}\\
&\leq& C|t\log|t|\,|+C|t|+C|t|^{\frac{1}{2}}.
\end{eqnarray*}
\cqfd

Returning to the proof of Theorem \ref{theorem:neck}, let $\Phi_0$ be the solution of
$d\Phi_0=\Phi_0\xi_0$ in $D(0,\varepsilon)$ with initial condition
$\Phi_0(1)=I_n$.
Let $Y_t(s)$ be the solution on $[0,1]$ of the Cauchy Problem
$$\left\{\begin{array}{l}
Y'_t(s)=Y_t(s)\xi_t(\alpha_t(s))\alpha'_t(s)\\
Y_t(0)=I_2.\end{array}\right.$$
By definition, $\cal{P}(\xi_t,\alpha_t)=Y_t(1)$.
Define
$$Z_t(s)=Y_t(s)-\Phi_0(\alpha_t(s)).$$
Then
\begin{eqnarray*}
Z'_t(s)&=&Y_t(s)\xi_t(\alpha_t(s))\alpha'_t(s)-\Phi_0(\alpha_t(s))\xi_0(\alpha_t(s))\alpha'_t(s)\\
&=&Z_t(s)\xi_t(\alpha_t(s))\alpha'_t(s)+\Phi_0(\alpha_t(s))[\xi_t(\alpha_t(s))-\xi_0(\alpha_t(s))]\alpha'_t(s).
\end{eqnarray*}
Hence
\begin{eqnarray*}
\lefteqn{\|Z_t(s)\|=\left\|\int_0^s Z'_t(x)dx\right\|
\leq\int_0^s\|Z_t(x)\|\,\|\xi_t(\alpha_t(x))\alpha'_t(x)\|dx}\\
&+&\int_0^s\|\Phi_0(\alpha_t(x))\|\,\|\left[\xi_t(\alpha_t(x))-\xi_0(\alpha_t(x))\right]\alpha'_t(x)\|dx
\end{eqnarray*}
By Gr\"onwall inequality:
$$\|Z_t(1)\|\leq \int_0^1\|\Phi_0(\alpha_t(s))\|\,\|\left[\xi_t(\alpha_t(s))-\xi_0(\alpha_t(s))\right]
\alpha'_t(s)\|ds\times\exp\left(\int_0^1\|\xi_t(\alpha_t(s))\alpha'_t(s)\|ds\right).$$
Using Lemma \ref{lemma:neck}, uniform bounds for $\Phi_0$ and $\xi_0$ in $D(0,1)$ and
the length of $\alpha_t$, we obtain
$$\|\cal{P}(\xi_t,\alpha_t)-\Phi_0(\alpha_t(1))\|=\|Z_t(1)\|\leq C|t|^{1/2}.$$
Since $\Phi_0$ is holomorphic in $D(0,1)$,
$$\|\Phi_0(\alpha_t(1))-\Phi_0(0)\|\leq C|\alpha_t(1)|=C|t|^{1/2}$$
Hence
\begin{equation}
\label{eq:neck-estimate5}
\|\cal{P}(\xi_t,\alpha_t)-\Phi_0(0)\|\leq C|t|^{1/2}.
\end{equation}
Let $\wh{\Phi}_0$ be the solution of $d\wh{\Phi}_0=\wh{\Phi}_0\wh{\xi}_0$ with initial condition
$\wh{\Phi}_0(1)=I_n$.
By the same argument, we have
\begin{equation}
\label{eq:neck-estimate6}
\|\cal{P}(\xi_t,\wh{\alpha}_t)-\wh{\Phi}_0(0)||\leq C|t|^{1/2}.
\end{equation}
By Equations \eqref{eq:neck-estimate5} and \eqref{eq:neck-estimate6}:
$$\|\cal{P}(\xi_t,\beta_t)-\Phi_0(0)\wh{\Phi}_0(0)^{-1}||=
\|\cal{P}(\xi_t,\alpha_t)\cal{P}(\xi_t,\wh{\alpha}_t)^{-1}-\Phi_0(0)\wh{\Phi}_0(0)^{-1}||
\leq C|t|^{1/2}.$$
Since $\cal{P}(\xi_t,\gamma)=I_2+O(t)$, we finally obtain
$$\left\|F(t)-\Phi_0(0)\wh{\Phi}_0(0)^{-1}\right\|\leq C|t|^{1/2}.$$
By Riemann Extension Theorem, $F$ extends holomorphically at $t=0$, and
$$F(0)=\Phi_0(0)\wh{\Phi}_0(0)^{-1}=\cal{P}(\xi_0,1,0)\cal{P}(\wh{\xi}_0,0,1).$$
Finally, to prove Point (3), assume that $t>0$ and write
$$\cal{P}(\xi_t,\beta_t)=\exp\left(\frac{t\log t}{2\pi i}t^{-1}\log \cal{P}(\xi_t,\gamma)\right)F(t).$$
Since $\cal{P}(\xi_0,\gamma)=I_2$, $t^{-1}\log\cal{P}(\xi_t,\gamma)$ extends holomorphically at $t=0$
with value $\frac{\partial}{\partial t}\cal{P}(\xi_t,\gamma)\mid_{t=0}$ and Point (3) follows.
\cqfd

\section{Differentiability of smooth functions of $t$ and $t\log t$}
\label{appendix:tlogt}
\begin{proposition}
\label{proposition:tlogt}
Let $E$ be a finite dimensional space and $g(t,s,z)$ be a smooth function from a neighborhood of $(0,0,z_0)$ in $\R^2\times E$ to a normed space $F$.
Define
$$f(t,z)=\left\{\begin{array}{ll}
g(t,t\log|t|,z)&\mbox{ if $t\neq 0$}\\
g(0,0,z)&\mbox{ if $t=0$.}\end{array}\right.$$
Assume that $g(0,s,z)$ only depends on $z$.
Then $f$ is of class $C^1$ and
$$df(0,z)=\frac{\partial g}{\partial t}(0,0,z)dt+d_zg(0,0,z).$$
\end{proposition}
Proof: $f$ is clearly continuous. For $t\neq 0$, we have by the chain rule:
$$df(t,z)=\frac{\partial g}{\partial t}(t,t\log |t|,z)dt+
\frac{\partial g}{\partial s}(t,t\log|t|,z)(1+\log|t|)dt+
d_zg(t,t\log|t|,z).$$
From the hypothesis, $\frac{\partial g}{\partial s}(0,s,z)=0$ so
$$\|\frac{\partial g}{\partial s}(t,s,z)\|=
\|\frac{\partial g}{\partial s}(t,s,z)-\frac{\partial g}{\partial s}(0,s,z)\|=O(t).$$
Hence
$$\lim_{(t,z)\to(0,z_0)}df(t,z)=\frac{\partial g}{\partial t}(0,0,z_0)dt+d_zg(0,0,z_0).$$
It follows (using the Mean Value Inequality) that $f$ is differentiable at $(0,z_0)$ and that it is
of classe $C^1$.
\cqfd


\begin{thebibliography}{9}
\bibitem{bobenko-heller-schmitt}
A. Bobenko, S. Heller, N. Schmitt:
{\em Minimal $n$-noids in hyperbolic and Anti-de Sitter 3-space.}
arXiv:1902.07992 (2019).
\bibitem{chae}
S. B. Chae. {\em Holomorphy and calculus in normed spaces.}
Monographs and textbooks in pure and applied mathematics, vol. 92 (1985).
\bibitem{dorfmeister-haak}
J. Dorfmeister, G. Haak:
{\em Constant mean curvature surfaces with periodic metric.}
Pacific Journal of Mathematics 182 (1998), 229--287.
\bibitem{dorfmeister-pedit-wu}
J. Dorfmeister, F. Pedit, H. Wu:
{\em Weierstrass type representation of
harmonic maps into symmetric spaces.}
Communications in Analysis and Geometry 6 (1998), 633-668.
\bibitem{dorfmeister-wu}
J. Dorfmeister, H. Wu:
{\em Construction of constant mean curvature n-noids
from holomorphic potentials.}
Mathematische Zeitschrift 258 (2008), 773--803.
\bibitem{fay} J.D. Fay: {\em Theta Functions on Riemann Surfaces}. Lecture Notes 
in Mathematics 352 (1973).
%\bibitem{forster}
%O. Forster:
%{\em Lectures on Riemann surfaces.}
%Graduate texts in Mathematics, Springer Verlag.
%\bibitem{fujimori-kobayashi-rossman}
%S. Fujimori, S. Kobayashi, W. Rossman:
%{\em Loop group methods for constant mean curvature surfaces.}
%Rokko Lectures in Mathematics 17, Koebe University.
%arXiv:math/0602570.
%\bibitem{gerding-pedit-schmitt}
%A. Gerding, F. Pedit, N. Schmitt:
%{\em Constant mean curvature surfaces: an integrable systems perspective.}
%Harmonic maps and differential geometry,
%Contemp. Math. 542 (2011), Amer. Math. Soc., 7--39.
\bibitem{grochenig}
K. Gr\"ochenig:
{\em Weight functions in time-frequency analysis.}
Fields Institute Communications 52, (2007), 343--366.
%\bibitem{karsten-kusner-sullivan}
%K. Gro\ss e-Brauckmann, R. Kusner, J. Sullivan:
%{\em Triunduloids: embedded constant mean curvature surfaces with three ends and genus zero.}
%J. Reine Angew. Math. 564 (2003), 35--61.
%\bibitem{karsten-kusner-sullivan2}
%K. Gro\ss e-Brauckmann, R. Kusner, J. Sullivan:
%{\em Coplanar constant mean curvature surfaces.}
%Comm. Anal. Geom. 15 (2007), no. 5, 985--1023.
%\bibitem{heller-heller-schmitt}
%L. Heller, S. Heller, N. Schmitt:
%{\em Navigating the Space of Symmetric CMC Surfaces.}
%arXiv:1501.01929.
\bibitem{heller1}
S. Heller:
{\em Higher genus minimal surfaces in $\S^3$ and stable bundles.}
J. Reine Angew. Math. 685 (2013), 105--122. 
\bibitem{heller2}
S. Heller:
{\em Lawson's genus two surface and meromorphic connections.}
Mathematische Zeitschrift 274 (2013), 745--760.
\bibitem{heller3}
S. Heller:
{\em A spectral curve approach to Lawson symmetric CMC surfaces of genus 2.}
Math. Annalen 360, Issue 3 (2014), 607--652.
\bibitem{lawson}
L. Heller, S. Heller, M. Traizet:
{\em Area estimates for high genus Lawson surfaces via DPW.}
arXiv 1907.07139 (2019).
\bibitem{kapouleas}
N. Kapouleas:
{\em Complete constant mean curvature surfaces in euclidean three-space.}
Annals of Mathematics 131 (1990), 239--330.
\bibitem{kilian-kobayashi-rossman-schmitt}
M. Kilian, S. Kobayashi, W. Rossman, N. Schmitt:
{\em Constant mean curvature surfaces of any positive genus.}
J. London Math. Soc. 72 (2005), 258--272.
\bibitem{kilian-mcintosh-schmitt}
M. Kilian, I. McIntosh, N. Schmitt:
{\em New constant mean curvature surfaces.}
Experiment. Math. 9 (2000), 595--611.
\bibitem{kilian-rossman-schmitt}
M. Kilian, W. Rossman, N. Schmitt:
{\em Delaunay ends of constant mean curvature surfaces.}
Compositio Mathematica 144 (2008), 186--220.
%\bibitem{korevaar-kusner-solomon}
%N. Korevaar, R. Kusner, B. Solomon:
%{\em The structure of complete embedded surfaces with constant mean curvature.}
%J. Diff. Geom. 30 (1989), 465--503.
%\bibitem{kusner-mazzeo-pollack}
%R. Kusner, R. Mazzeo, D. Pollack:
%{\em The moduli space of complete embedded constant mean curvature surfaces.}
%GAFA 6, Issue 1 (1996), 120--137.
\bibitem{masur} H. Masur: {\em The extension of the Weil-Petersson metric to the 
boundary of Teichmuller space}. Duke Math. J. 43  (1976), 623--635.
\bibitem{mazzeo-pacard}
R. Mazzeo, F. Pacard:
{\em Constant mean curvature surfaces with Delaunay ends.}
Comm. Anal. and Geom. 9, no. 1 (2001), 169--237.
\bibitem{presley-segal}
A. Pressley, G. Segal: {\em Loop Groups}. Oxford Mathematical Monographs,
Oxford University Press, Oxford, New York (1988).
\bibitem{raujouan}
T. Raujouan:
{\em On Delaunay ends in the DPW method.}
To appear in Indiana Univ. Math. J. %arXiv:1710.00768.
\bibitem{raujouan2}
T. Raujouan:
{\em Constant mean curvature n-noids in hyperbolic space.}
 arXiv:1905.09096 (2019).
%\bibitem{schmitt}
%N. Schmitt:
%{\em Constant mean curvature $n$-noids with platonic symmetries.}
%arXiv:math/0702469.
\bibitem{schmitt-kilian-kobayashi-rossman}
N. Schmitt, M. Kilian, S. Kobayashi, W. Rossman:
{\em Unitarization of monodromy
representations and constant mean curvature trinoids in 3-dimensional space
forms.}
Journal of the London Mathematical Society 75 (2007), 563--581.
%\bibitem{taylor} M. Taylor: {\em Introduction to Differential Equations.}
%Pure and Applied Undergraduate Texts 14, American Math. Soc. (2011).
\bibitem{teschl} G. Teschl: Ordinary differential equations and dynamical systems.
{\em Graduate Studies in Mathematics} 140, American Math. Soc. (2010).
\bibitem{nosym}
M. Traizet: An embedded minimal surface with no symmetries.
{\em J. Differential Geometry}, 60(1) (2002), 103--153.
\bibitem{triply} M. Traizet:
{\em On the genus of triply periodic minimal surfaces.}
J. Diff. Geom. 79 (2008), 243--275.
\bibitem{crelle} M. Traizet:
{\em Opening infinitely many nodes.}
J. reine angew. Math. 684 (2013), 165--186.
\bibitem{bryant}
M. Traizet: 
{\em Opening nodes on horosphere packings.}
Trans. Amer. Math. Soc. 368 (2016), 5701--5725.
\bibitem{nnoids}
M. Traizet:
{\em Construction of constant mean curvature $n$-noids using
the DPW method.}
To appear in J. reine angew. Math.
\bibitem{minoids}
M. Traizet:
{\em Gluing Delaunay ends to minimal n-noids using the DPW method.}
Mathematische Annalen, 377(3) (2020), 1481--1508
\bibitem{opening-nodes}
M. Traizet:
{\em Opening nodes and the DPW method.}
arXiv:1808.01366  (2018).
\end{thebibliography}
\end{document}